\documentclass[11pt,reqno]{amsproc}
\numberwithin{equation}{section}
\usepackage[table,xcdraw]{xcolor}
\usepackage{color}
\usepackage{mathrsfs}
\usepackage{fullpage}
\usepackage{caption}
\usepackage{amsmath}
\usepackage{graphicx}
\usepackage{subfigure}
\usepackage{float}
\DeclareGraphicsExtensions{.eps}
\usepackage{ucs}
\usepackage{latexsym}
\usepackage[utf8x]{inputenc}
\usepackage{epstopdf}
\usepackage{mathtools}
\usepackage[semicolon,square,authoryear]{natbib}
\usepackage[debug=false, colorlinks=true, pdfstartview=FitV, 
linkcolor=blue, citecolor=blue, urlcolor=blue]{hyperref}
\usepackage{algpseudocode}
\usepackage{algorithm}
\usepackage[most]{tcolorbox}
\usepackage{stmaryrd}

\newtheorem{remark}{Remark}[section]
\newtheorem{proposition}{Proposition}[section]

\usepackage{subfigure}
\usepackage{multirow}
\usepackage{makecell, boldline}
\usepackage{url}
\usepackage{listings}
\usepackage{soul}
\usepackage[unnumbered,framed]{mcode}

\usepackage{morefloats}
\newlength{\drop}
\definecolor{amethyst}{rgb}{0.6, 0.4, 0.8}
\definecolor{burgundy}{rgb}{0.5, 0.0, 0.13}
\newcommand\bfx{{\bf x}}

\title{\textbf{A vertex scheme for two-phase flow in heterogeneous media}}
\author{
	\textbf{M.~S.~Joshaghani},
	\textbf{V.~Girault},
	and
	\textbf{B.~Riviere\footnote{Supported by NSF-DMS 1913291.}}
	\\
	{\small 
\textbf{Correspondence to:}~\textsf{m.sarraf.j@rice.edu}}}

\newsavebox{\measurebox}
\begin{document}
	
	\date{\today}
	
	\begin{titlepage}
		\drop=0.1\textheight
		\centering
		\vspace*{\baselineskip}
		\rule{\textwidth}{1.6pt}\vspace*{-\baselineskip}\vspace*{2pt}
		\rule{\textwidth}{0.4pt}\\[0.25\baselineskip]
		{\Large \textbf{\color{burgundy}
				A vertex scheme for two-phase flow in heterogeneous media}}

		\rule{\textwidth}{0.4pt}\vspace*{-\baselineskip}\vspace{2pt}
		\rule{\textwidth}{1.6pt}\\[0.25\baselineskip]
		\scshape
		An e-print of the paper will be made available on arXiv. \par 
		\vspace*{0.3\baselineskip}
		Authored by \\[0.3\baselineskip]
		
		{\Large M.~S.~Joshaghani\par}
		{\itshape Postdoctoral Research Associate, Rice University, Houston, Texas 77005 \\ 
        \textbf{phone:} +1-281-781-5331, \textbf{e-mail:} m.sarraf.j@rice.edu}\\[0.25\baselineskip]
		
		{\Large V.~Girault\par}
		{\itshape Professor Emeritus, Laboratoire Jacques-Louis Lions \\
			University Pierre et Marie Curie, France}\\

		{\Large B.~Riviere\par}
		{\itshape Noah Harding Chair and Professor of Computational and Applied Mathematics  \\
			Rice University, Houston, Texas 77005} \\ 
	\vspace{1cm}
		\begin{figure}[h]
			\centering
			\subfigure[Highly-varying permeability field ($\mbox{m}^2$) in logscale \label{Fig:mainPageL}]{
				\includegraphics[clip,scale=0.5,trim=0 0cm 0.8cm 0]{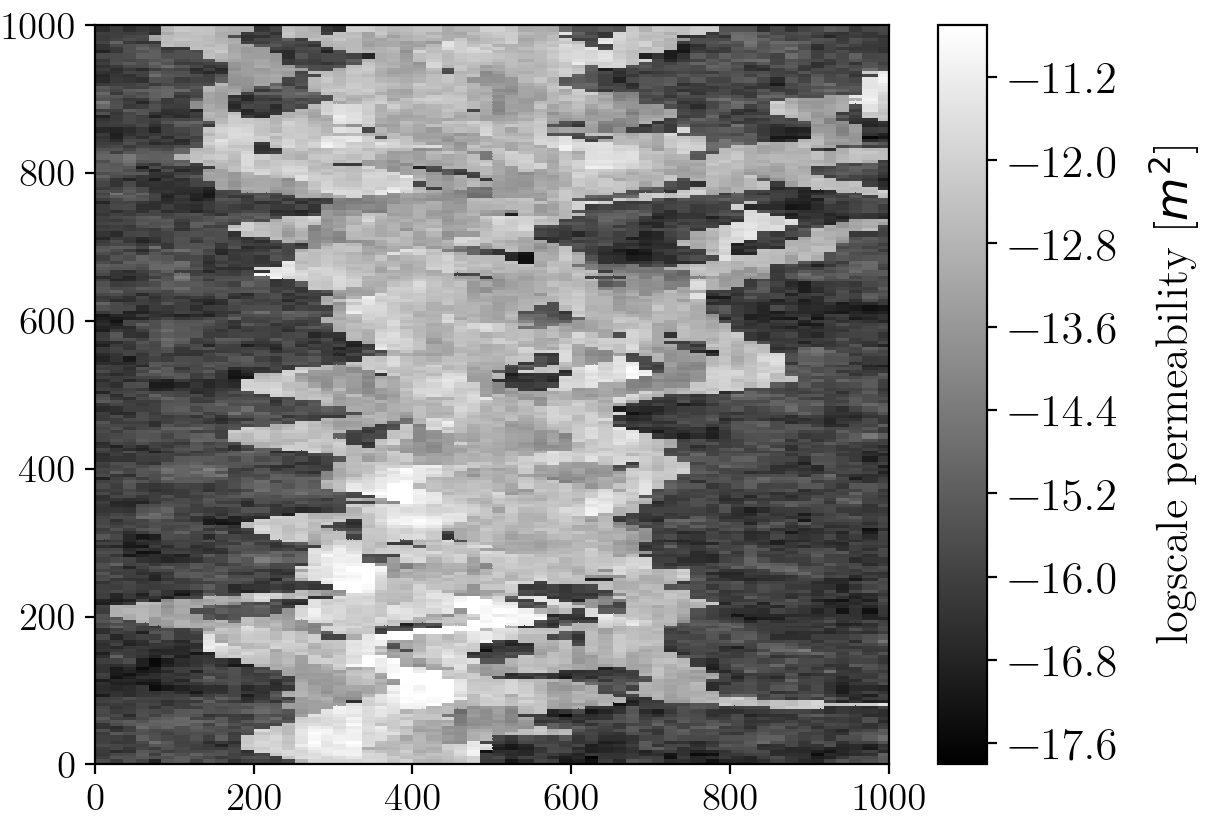}}
			\hspace{1cm}
			\subfigure[Wetting phase saturation; $t=0.725$ days  \label{Fig:mainPageR}]{
				\includegraphics[clip,scale=0.13,trim=0 3.5cm 0cm 0cm
				]{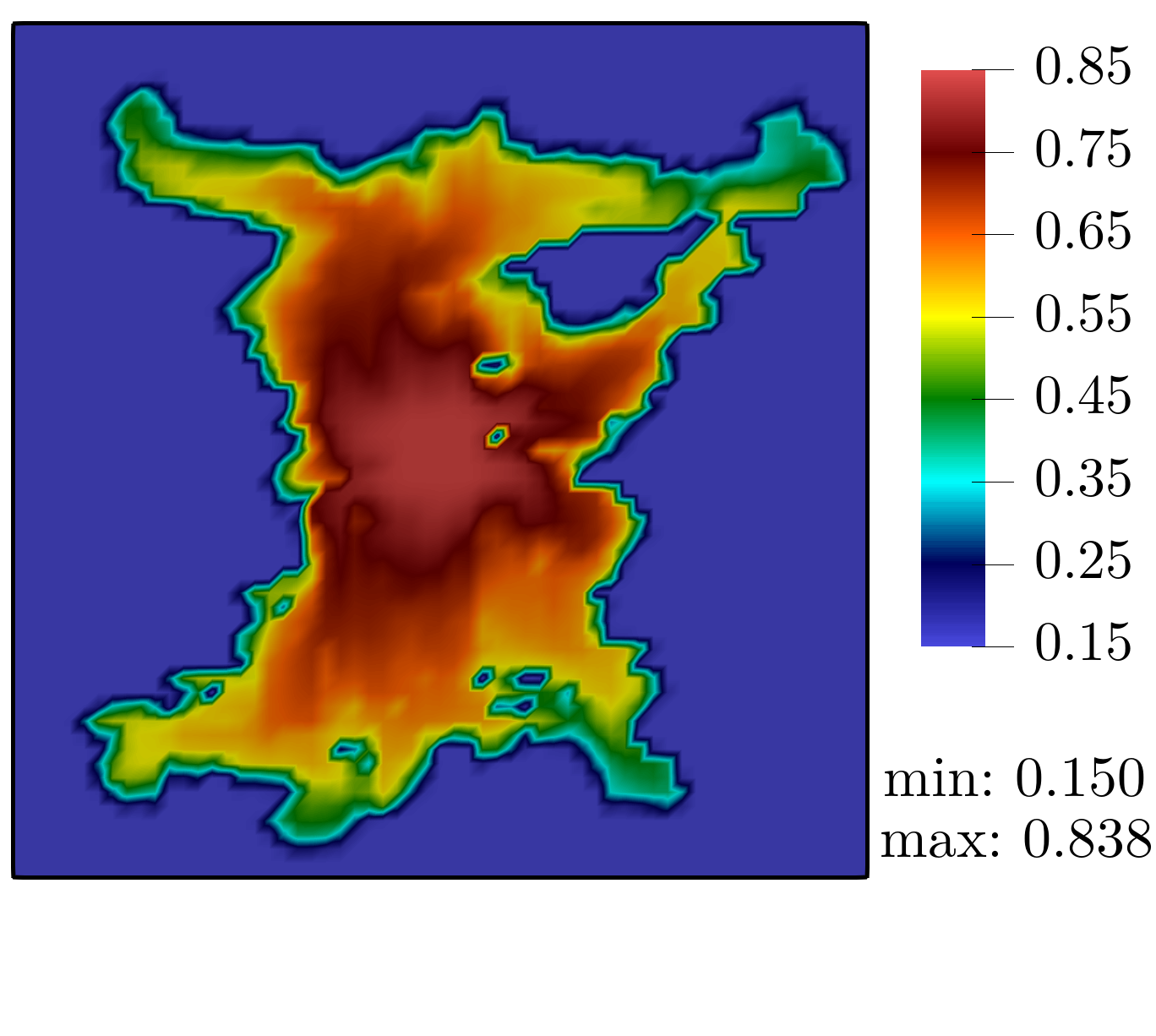}}

			\emph{{\small{
						Left figure shows a realistic discontinuous  permeability field in a 
                        $1000~\mbox{m} \times 1000~\mbox{m}$ domain.
						Right figure shows the saturation profile obtained under the proposed 
						finite element scheme. 
		                 We observe that: (i) the permeability field determines the pattern fluid flows through porous media
                         (ii) the proposed scheme exhibits satisfactory results with respect to maximum principle. 
                         This means that saturation solution always remains between $s_{rw}$ and $1-s_{ro}$ (between
                         $0.15$ and $0.85$ in this problem).
				      }}
					}
					\captionsetup{labelformat=empty}
				\vspace{-.9cm}
				\caption{}
		\end{figure}
		\vfill
		{\scshape 2021} \\
		{\small Computational Modeling of Porous Media (COMP-M) Group} \par
	\end{titlepage}
\setcounter{figure}{0}
\begin{abstract}
    This paper presents the numerical solution of immiscible two-phase flows in porous media, 
    obtained by a first-order finite element method equipped with mass-lumping and flux upwinding. 
    The unknowns are the physical phase pressure and phase saturation. 
    %
    Our numerical experiments confirm that the method converges optimally for manufactured solutions. 
    For both structured and unstructured meshes, we observe the high-accuracy wetting saturation profile 
    that ensures minimal numerical diffusion at the front. Performing several examples of quarter-five 
    spot problems in two and three dimensions, we show that the method can easily handle heterogeneities 
    in the permeability field. 
    Two distinct features that make the method appealing to reservoir simulators are: (i)  maximum principle
is satisfied, and (ii) mass balance is locally conserved.
    
    %

\end{abstract}
\maketitle
\section{INTRODUCTION}
\label{S1_Introduction}

A first-order finite element method is proposed to simulate two-phase flows in heterogeneous porous media. 
The method is defined for unstructured meshes made of simplices for two-dimensional or three-dimensional domains. Mass-lumping and upwinding techniques are 
employed to discretize the mass and stiffness matrices.  The proposed method solves for primary unknowns that are physical quantities,
namely the wetting phase pressure and the wetting phase saturation. Flows are driven by injection and production wells, represented by source and sink functions.
The method has recently been analyzed in the papers by \cite{GiraultRiviereCappaneraI,GiraultRiviereCappaneraII}. The fact that the relative permeabilities vanish when evaluated at the residual saturation values and that the capillary pressure has an unbounded derivative make the proofs for the well-posedness and convergence analysis of the scheme technical and complicated. In this current work, we extend the scheme to heterogeneous porous media for which the permeability field varies over several orders of magnitude across the domain. 
Several simulations of incompressible two-phase flow in two-dimensional and three-dimensional domains show the accuracy and robustness of the proposed mass-lumped upwinded finite element method. In this paper, our proposed scheme will be referred to as the ``vertex scheme'' because mass-lumping reduces the finite element integrals to quantities involving the values of the solution at the vertices.
Thanks to the use of mass-lumping and upwinding, the maximum principle is proved and observed in our computations.  The numerical solution of the saturation does not exhibit local oscillations near the front, which indicates the monotonicity of the scheme. 

The numerical modeling of incompressible two-phase flow in porous media has been widely studied in the literature. Only a small fraction of the proposed schemes with physical primary unknowns, has been theoretically analyzed. 
Besides our proposed finite element method \citep{GiraultRiviereCappaneraI,GiraultRiviereCappaneraII}, a cell-centered finite volume method has been analyzed 
in \citep{Eymard2003}. Many of the convergence works in the literature impose 
unrealistic constraints on the input data, in particular the relative permeabilities are assumed to be bounded below by positive 
constants and the derivative of the capillary pressure is assumed to be bounded \citep{ohlberger,epshteyn2009analysis,EymardGuichard}. 
If  non-physical primary unknowns are chosen, such as the global pressure introduced by \citet{chavent1986mathematical}, 
the degeneracy of the relative permeabilities can be circumvented (see \citep{douglas1983finite,Michel03,chen2001degenerate}).
The idea of using mass-lumping  has been proposed for the numerical solution of various partial differential equations 
(see for instance \citep{chen1985lumped,cohen2001higher}) as well as in the application of porous media 
(see \citep{Forsyth91,abriola1993mass}). Mass-lumping within the finite element method yields a diagonal mass matrix, 
which is a desirable feature for reducing the computational cost. 
Upwinding helps reducing the numerical oscillations near fronts in convection-dominated problems. 

An outline of the paper follows. In Section~\ref{Sec:S1_GEs}, the two-phase flow model is presented with wetting phase saturation 
and pressure as primary unknowns. The vertex scheme is defined in Section~\ref{Sec:S2_CG} and the resulting linearized system 
is described in Section~\ref{Sec:Solver}. Numerical simulations are shown in Section~\ref{Sec:numer} and are followed by conclusions.

\section{GOVERNING EQUATIONS}
\label{Sec:S1_GEs}
The incompressible two-phase flow model in a porous medium $\Omega\subset\mathbb{R}^d$, $d=2,3$, 
over a time interval $(0,T)$ is characterized by the following coupled equations:
\begin{subequations}
\begin{alignat}{3}
    \label{Eqn:GE_Wphase}
    \partial_{t} (\phi s) -\nabla \cdot (\eta_w(s) K \nabla p) &=
        f_w(s_{\mathrm{in}})\bar{q}-f_w(s) \underline{q},    
    \quad &&\mathrm{in}\; \Omega\times (0,T), \\
    \label{Eqn:GE_pressure}
   -\partial_{t} (\phi s) -\nabla \cdot \left(\eta_{o}(s)  K (
          \nabla p_c(s) + \nabla p)\right) &= 
        f_o(s_{\mathrm{in}})\bar{q}-f_o(s) \underline{q},    
    \quad && \mathrm{in}\; \Omega \times (0,T),\\
     \eta_{w}(s) K \nabla p \cdot \mathbf{n} &= 0, 
    && \mathrm{on} \; \partial \Omega \times (0,T),\\
     \eta_{o}(s) K \nabla p_o \cdot \mathbf{n} &= 0, 
    && \mathrm{on} \; \partial \Omega\times (0,T).
\end{alignat}
\end{subequations}
The primary unknowns are the wetting phase pressure, $p$, and wetting phase saturation, $s$.
The secondary unknowns, namely the non-wetting phase pressure and saturation, $(p_o,s_o)$, are recovered by using the relations:
\[
p_o = p_c(s)+p, \quad s_o = 1-s.
\]
The porosity and the permeability of the medium are denoted by $\phi$ and $K$ respectively. The mobilities, $\eta_\alpha$, are the
ratios of the relative permeabilities, $k_{r\alpha}$, to the phase viscosities, $\mu_\alpha>0$, for $\alpha=o,w$:
\begin{align}
    \eta_{\alpha}(s)=\frac{k_{r\alpha}(s)}{\mu_{\alpha}},
    \quad \alpha = w,o.
\end{align}
Both capillary pressure, $p_c$, and relative permeabilities are functions of the saturation 
(see \citep{brooks1964hydrau,van1980closed}).
In this work, the Brooks-Corey model is used.  The residual saturations, $s_{rw}, s_{ro}$, are constants in the interval $[0,1]$.
\begin{align}
    k_{rw}(s) &=\bar{s}^{\frac{2+3\theta}{\theta}},
    \quad
    k_{ro}(s) = (1-\bar{s})^2(1-\bar{s}^{\frac{2+3\theta}{\theta}}), \quad
    \bar{s} = \frac{s-s_{rw}}{1-s_{rw}-s_{ro}}, \label{eq:krwo} \\
    p_c(s) &= 
    \begin{cases}
        p_d \bar{s}^{\frac{-1}{\theta}} \quad \mbox{if } \bar{s}>R \\
        p_d R^{\frac{-1}{\theta}}-\frac{p_d}{\theta}R^{(-1-\frac{1}{\theta})} (\bar{s}-R) \quad \mbox{otherwise}.
    \end{cases}
\label{eq:pc}
\end{align}

This model introduces a  parameter $\theta \in [0.2,3.0]$, which characterizes the inhomogeneity of the medium. 
The entry pressure, $p_d$, is a constant pressure corresponding to the capillary pressure required
to displace the fluid from the largest pore.
The fractional flows of each phase are related to the mobilities as follows:
\begin{align}
    f_w(s) = \frac{\eta_w(s)}{\eta_w(s)+\eta_o(s)}, \quad
    f_o(s) = 1 - f_w(s).
\end{align}
Flow rates at the injection and production wells,
$\bar{q}$, and $\underline{q}$ satisfy:
\begin{align}
    \bar{q}\ge 0, \quad
    \underline{q} \ge 0, \quad 
    \int_{\Omega} \bar{q} = 
    \int_{\Omega} \underline{q},
\end{align}
and the saturation at the injection wells is set equal to a constant value $s_\mathrm{in}$.
Finally the model is completed by the initial condition:
\begin{equation}
    \label{Eqn:GE_IC}
     s = s^{0},\quad \mathrm{in} \; \Omega.
\end{equation}

\section{NUMERICAL SCHEME}
\label{Sec:S2_CG}
The domain $\Omega$ is partitioned into triangular elements in 2D and tetrahedral elements in 3D. Let $h$ denote
the maximum diameter of each element in the mesh $\mathcal{E}_h$.  
Let $\tau$ denote the time step size and let $p_h^n, s_h^n$ denote the discrete pressure and saturation
respectively at time $t^n = n\tau$.  They belong to the finite element space $X_h$ of order one:
\[
X_h = \{ v_h \in \mathcal{C}^0(\bar\Omega): \, \forall E\in\mathcal{E}_h, \, v_h|_E \in \mathbb{P}_1(E)\}.
\]
Let $M$ be the dimension of $X_h$; it is the number of nodes (i.e. vertices) of $\mathcal{E}_h$.  Let $\Phi_i$ be the
Lagrange basis function, that is piecewise linear and takes the value $1$ at node $i$ and $0$ at all the other nodes.
We write
\begin{equation}\label{eq:basislin}
p_h^n(\bfx) = \sum_{i=1}^M P_i^n \Phi_i(\bfx), \quad s_h^n(\bfx) = \sum_{i=1}^M S_i^n \Phi_i(\bfx), \quad \bfx\in\Omega.
\end{equation}
We now define coefficients that arise from the mass-lumping technique.
For a fixed node $i$, the macro-element  $\Delta_i$ is defined as the union of elements sharing the node $i$.
\[
c_{ij,E} = \int_E \vert \nabla \Phi_i \cdot \nabla\Phi_j \vert, \quad
c_{ij} = \sum_{E\in\Delta_i\cap\Delta_j} c_{ij,E}, \quad
c_{ij}(K) = \sum_{E\in\Delta_i\cap\Delta_j} K_E \, c_{ij,E}, \quad \forall 1\leq i,j\leq M,
\]
\[
m_i = \frac{\vert \Delta_i\vert}{d+1}, \quad \forall 1\leq i, j\leq M.
\]
We assume here that the permeability $K$ is piecewise constant and we denote by $K_E$ the constant value that is the
restriction of $K$ on the element $E$. Clearly, if $K$ is constant everywhere, then $c_{ij}(K) = K c_{ij}$.

We first introduce the nonlinear scheme, that is written
with respect to the nodal values of the numerical pressure and saturation. 
For $n\geq 1$, given $s_h^{n-1}\in X_h$, find $(p_h^{n}, s_h^{n}) \in X_h\times X_h$ satisfying \eqref{eq:basislin} and 
\begin{align}
    \label{eq:or1}
    m_i\phi\frac{S_i^{n}-S_i^{n-1}}{\tau}
    -\sum_{j=1}^M c_{ij}(K) \eta_{w}(S_{w,ij}^{n}) (P_j^{n}-P_i^{n}) \nonumber \\
    = m_i (f_w(s_{\mathrm{in}})\bar{q}_i-f_w(S_i^{n-1})\underline{q}_i ), \quad 1\le i\le M-1,
\end{align}
\begin{align}
    \label{eq:or2}
    -m_i \phi \frac{S_i^{n}-S_i^{n-1}}{\tau}
    -\sum_{j=1}^M c_{ij}(K) \eta_{o}(S_{o,ij}^{n}) (P_j^{n}-P_i^{n}) \nonumber \\
    -\sum_{j=1}^M c_{ij}(K) \eta_{o}(S_{o,ij}^{n})\left(
        p_c(S_j^{n-1})  + p_c'(S_j^{n-1})
     (S_j^{n}-S_j^{n-1})
    -p_c(S_i^{n-1}) - p_c'(S_i^{n-1})
    (S_i^{n}-S_i^{n-1})
    \right) \nonumber \\
    = m_i (f_o(s_{\mathrm{in}})\bar{q}_i-f_o(S_i^{n-1})\underline{q}_i), \quad 1\le i\le M, 
\end{align}
\begin{align}
    \sum_{i=1}^M m_i \; P_i^{n} = 0. \label{eq:or3}
\end{align}
The values $S_{w,ij}^{n}$ and $S_{o,ij}^n$ are upwind values, i.e. they are nodal values of the saturation at either
node $i$ or node $j$, that are made precise in the linearized scheme below.
In the case of constant permeability, well-posedness and convergence of the nonlinear scheme are proved 
in~\citep{GiraultRiviereCappaneraI,GiraultRiviereCappaneraII}. 
\begin{proposition}
Let $(s,p)$ be a weak solution to problem \eqref{Eqn:GE_Wphase}-\eqref{Eqn:GE_pressure}. Assume that the porosity
$\phi$ and permeability $K$ are positive constants. 
As the mesh size $h$ and time step size $\tau$ tend to zero, the discrete saturation satisfying
\eqref{eq:or1}-\eqref{eq:or3} converges, up to a subsequence, strongly to $s$ 
in the $L^2$ norm and the discrete pressure converges, up to a subsequence, weakly to $p$. 
In addition, the saturation satisfies the maximum principle:
\begin{equation}\label{eq:maxsat}
s_{rw} \leq s_h^n(\bfx) \leq 1-s_{ro}, \quad \forall \bfx\in \Omega.
\end{equation}
\label{prop:conv}
\end{proposition}
%

\underline{Linearized Scheme:} We linearize the equations \eqref{eq:or1}-\eqref{eq:or2} by using a fixed point iteration and 
approximating the capillary
pressure by a first-order Taylor expansion:
\[
p_c(S_i^{n}) \approx p_c(S_i^{n-1}) + p_c'(S_i^{n-1}) (S_i^{n}-S_i^{n-1}).
\]
At each time step $t^n$, we will solve for 
a sequence of nodal values $(P_i^{n,k}, S_i^{n,k})$ where the superscript $k$ denotes the fixed-point iteration number.
\begin{align}
    \label{eq:scheme1}
    m_i\phi\frac{S_i^{n,k}-S_i^{n-1}}{\tau}
    -\sum_{j=1}^M c_{ij}(K) \eta_{w}(S_{w,ij}^{n,k-1}) (P_j^{n,k}-P_i^{n,k}) \nonumber \\
    = m_i (f_w(s_{\mathrm{in}})\bar{q}_i-f_w(S_i^{n-1})\underline{q}_i ), \quad 1\le i\le M-1,
\end{align}
\begin{align}
    \label{eq:scheme2}
    -m_i \phi \frac{S_i^{n,k}-S_i^{n-1}}{\tau}
    -\sum_{j=1}^M c_{ij}(K) \eta_{o}(S_{o,ij}^{n,k-1}) (P_j^{n,k}-P_i^{n,k}) \nonumber \\
    -\sum_{j=1}^M c_{ij}(K) \eta_{o}(S_{o,ij}^{n,k-1})\left(
        p_c(S_j^{n-1})  + p_c'(S_j^{n-1})
     (S_j^{n,k}-S_j^{n-1})
    -p_c(S_i^{n-1}) - p_c'(S_i^{n-1})
    (S_i^{n,k}-S_i^{n-1})
    \right) \nonumber \\
    = m_i (f_o(s_{\mathrm{in}})\bar{q}_i-f_o(S_i^{n-1})\underline{q}_i), \quad 1\le i\le M, 
\end{align}
\begin{align}
    \sum_{i=1}^M m_i P_i^{n,k} = 0. \label{eq:scheme3}
\end{align}
We now make precise the choice of the upwind values, $S_{w,ij}^{n,k}$ and $S_{o,ij}^{n,k}$:
\begin{align}
    S_{w,ij}^{n,k-1} = 
\begin{cases}
  S_i^{n,k-1} \quad \mathrm{if} \quad P_i^{n,k-1} > P_j^{n,k-1} \\
  S_j^{n,k-1} \quad \mathrm{if} \quad P_i^{n,k-1} < P_j^{n,k-1} \\
  \max(S_i^{n,k-1}, S_j^{n,k-1}) \quad  \mathrm{if} \quad P_i^{n,k-1} = P_j^{n,k-1}, 
\end{cases}
\end{align}
\begin{align}
    \label{Eqn:upwind_sat}
    S_{o,ij}^{n,k-1} = 
\begin{cases}
    S_i^{n,k-1} \quad \mathrm{if} \quad p_c(S_{i}^{n,k-1}) + P_i^{n,k-1} > p_c(S_{j}^{n,k-1}) + P_j^{n,k-1} \\
    S_j^{n,k-1} \quad \mathrm{if} \quad p_c(S_{i}^{n,k-1}) + P_i^{n,k-1} < p_c(S_{j}^{n,k-1}) + P_j^{n,k-1} \\
  \min(S_i^{n,k-1}, S_j^{n,k-1}) \quad  \mathrm{if}  \quad p_c(S_{i}^{n,k-1}) + P_i^{n,k-1} = p_c(S_{j}^{n,k-1}) + P_j^{n,k-1}.
\end{cases}
\end{align}

We initialize the iterates with the values at the previous time-step:
\[
P_i^{n,0} = P_i^{n-1}, \quad S_i^{n,0} = S_i^{n-1}.
\]
Convergence is obtained when the difference between two iterates for both discrete pressure and saturation is small (less than $10^{-5}$)  {in the $L^{\infty}$ norm.  The nodal values of
the saturation and pressure at time $t^n$ are the nodal values of the converged iterates.
Since the finite element solutions $s_h^n$ and $p_h^n$ uniquely depend on the nodal values, they can be evaluated at any point
in the domain.

Finally, to start the algorithm, we choose for $s_h^0$ the Lagrange interpolant of the saturation $s_0$ and
for $p_h^0$ a constant value so that $S_{w,ij}^{0,0}$ and $S_{o,ij}^{0,0}$ are well defined. 

\begin{remark}
Equation~\eqref{eq:scheme1} is valid for $i=M$; this can be obtained by adding \eqref{eq:scheme1} and \eqref{eq:scheme2} and by using
\eqref{eq:scheme3}.
\end{remark}

\section{SOLVER METHODOLOGY} 
\label{Sec:Solver}
The fully discrete formulations \eqref{eq:scheme1}-\eqref{eq:scheme3} yield a $2\times 2$ block linear system of the form:
\begin{align}
    \label{Eqn:kxf}
    \begin{pmatrix}
        \mathbf{K}_{ss} & \mathbf{K}_{sp} \\
        \mathbf{K}_{ps} & \mathbf{K}_{pp}
    \end{pmatrix}
    \begin{pmatrix}
        \mathbf{s} \\
        \mathbf{p}
    \end{pmatrix}
    =
    \begin{pmatrix}
        \mathbf{f}_s \\
        \mathbf{f}_p
    \end{pmatrix}
\end{align}
where each block is of size $M\times M$ and has entries that depend on the time step and the Picard iterate.  
Because of the local support of the basis functions, the sums over all the nodes in \eqref{eq:scheme1} and \eqref{eq:scheme2} reduce to sums over a small set of nodes, which leads to sparse matrices.  To be precise, let $\mathcal{N}(i)$ be the set of indices of all
nodes in the macro-element $\Delta_i$. \\
The block $\mathbf{K}_{ss}$ is a diagonal matrix:
\[
(\mathbf{K}_{ss})_{ii} = \frac{m_i\phi}{\tau}, \quad 1\leq i\leq M-1,\quad
(\mathbf{K}_{ss})_{MM} = 0.
\]
The non-zero entries in the block $\mathbf{K}_{sp}$ are:
\[
(\mathbf{K}_{sp})_{ij} = -c_{ij}(K) \eta_w(S_{w,ij}^{n,k-1}), \quad 1\leq i\leq M-1, \, j \in \mathcal{N}(i),  \, j\neq i, \quad
(\mathbf{K}_{sp})_{ij} = m_j, \quad i=M, \, 1\leq j\leq M.
\]
The non-zero entries in the block $\mathbf{K}_{ps}$ are:
\[
(\mathbf{K}_{ps})_{ii} = -\frac{m_i\phi}{\tau}, \quad 1\leq i\leq M, \quad
(\mathbf{K}_{ps})_{ij} = - c_{ij}(K) \eta_o(S_{o,ij}^{n,k-1}) p_c'(S_j^{n-1}), \quad 1\leq i\leq M, \, j \in \mathcal{N}(i),  \, j\neq i.
\]
The non-zero entries in the block $\mathbf{K}_{pp}$ are:
\[
(\mathbf{K}_{pp})_{ij} = - c_{ij}(K) \eta_o(S_{o,ij}^{n,k-1}), \quad 1\leq i\leq M, \, j \in \mathcal{N}(i),  \, j\neq i.
\]
For completeness, we display the entries of the right-hand side vectors $\mathbf{f}_s$ and $\mathbf{f}_p$.
\[
(\mathbf{f}_s)_i = \frac{m_i \phi}{\tau} + m_i (f_w(s_{\mathrm{in}})\bar{q}_i-f_w(S_i^{n-1})\underline{q}_i), \quad 1\leq i \leq M-1, \quad
(\mathbf{f}_s)_M = 0,
\]
\begin{multline}
(\mathbf{f}_p)_i = -\frac{m_i \phi}{\tau} + m_i (f_o(s_{\mathrm{in}})\bar{q}_i-f_o(S_i^{n-1})\underline{q}_i)
\\
+\sum_{j\in\mathcal{N}(i)} c_{ij} \eta_{o}(S_{o,ij}^{n,k-1}) \left( p_c(S_j^{n-1})  - p_c'(S_j^{n-1}) S_j^{n-1}
    -p_c(S_i^{n-1}) + p_c'(S_i^{n-1}) S_i^{n-1}\right).
\end{multline}

It is worth noting that the construction of the global matrix $\mathbf{K}$ is done by assembling local matrices,
as this is usually done in the finite element framework.  
For example we describe the procedure for assembling the block $\mathbf{K}_{ps}$ in Algorithm \ref{Alg:Vertexassembler};
the other blocks $\mathbf{K}_{ss}$, $\mathbf{K}_{sp}$, and $\mathbf{K}_{pp}$ are handled similarly. 
Let $\mathbf{C}_{loc}^E$ be the local matrix associated with the $c_{ij}$ coefficients restricted to an element $E$. 
\[
(\mathbf{C}_{loc}^E)_{i_{loc},j_{loc}} = \int_E \vert \nabla \Psi_{i_{loc}} \cdot \nabla \Psi_{j_{loc}}\vert, \quad
\forall 1\leq i_{loc}, j_{loc}\leq d+1, 
\]
where the functions $\Psi_{i_{loc}}$ are linear polynomials on $E$ that correspond to the restriction of a global
basis $\Phi_i$ on $E$ for the node $i$ with local number equal to $i_{loc}$.

\begin{remark}
    In the case of a two-dimensional domain partitioned into a structured mesh of right-triangular elements of size $h$, the
local matrix $\mathbf{C}_{loc}^E$ is the same constant matrix for all elements $E$. Taking the local numbering counterclockwise 
and start from the right-angle node, $\mathbf{C}_{loc}^E$ reads as follows:
    \begin{align}
        \mathbf{C}_{loc}^E = 
        \begin{pmatrix}
            0.5 & 0.5 & 0.5\\
            0.5 & 0.5 & 0.0\\
            0.5 & 0.0 & 0.5
        \end{pmatrix}, \quad \forall E\in\mathcal{E}_h.
    \end{align}
  However, for unstructured meshes, the entries of the local matrix will depend on the element. 
\end{remark}

%
\begin{algorithm}
    \caption{Assembly for $\mathbf{K}_{ps}$ matrix at time step $t^{n}$ and Picard's iteration ${k}$.}
	\label{Alg:Vertexassembler}
	\begin{algorithmic}[0]
        \State Input: $(S_i^{n-1})_i, \, (S_i^{n,k-1})_i$, $(P_i^{n,k-1})_i$.
        \State Initialize $\mathbf{K}_{ps}$ to the zero matrix.
	\For{each element $E$ in mesh} 
        \State construct  $\mathbf{C}_{loc}^E$ \Comment For structured mesh: $\mathbf{C}_{loc}^E$ is constant over all elements
		%
        %
        \State evaluate $K_E = K|_E$ \Comment Value of permeability in element $E$.
        \For{$i_{loc} = \{1,2, \cdots, d+1\}$ } \Comment There are $(d+1)$  degrees-of-freedom per element
        \State $i = \mathrm{glodofs}(i_{loc})$ \Comment Global number of local node
        \State $\mathbf{K}_{ps}(i,i) = -\frac{\phi \vert E\vert}{(d+1)\tau}$ \Comment Mass-lumping operation
        \For{$j_{loc} = \{ 1,2,\cdots, d+1 \}\setminus \{i_{loc}\}$}
        \State $j = \mathrm{glodofs}(j_{loc})$ \Comment Global number of local node
        %
        %
        \If {$p_c(S_i^{n,k-1})+P_i^{n,k-1}>p_c(S_j^{n,k-1})+P_j^{n,k-1}$} \Comment Upwinding operation
        \State Etaij = $\eta_o(S_i^{n,k-1})$
        \ElsIf {$p_c(S_i^{n,k-1})+P_i^{n,k-1}<p_c(S_j^{n,k-1})+P_j^{n,k-1}$}
	\State Etaij = $\eta_o(S_j^{n,k-1})$
		\Else
        \State Etaij = $\eta_o\big(\min(S_i^{n,k-1},S_j^{n,k-1})\big)$

		\EndIf
		%
        %
        \State $\mathbf{K}_{ps}(i,j) = -\mbox{Etaij} \times \mathbf{C}_{loc}^E(i_{loc},j_{loc}) \times K_E \times$  $\frac{\mathrm{d}p_c}{\mathrm{d}s}(S_j^{n-1})$ 
		\EndFor
		%
		\EndFor 
        \EndFor  
	\end{algorithmic}
\end{algorithm}

A Schur complement approach is used to factorize the matrix $\mathbf{K}$ 
following \citep{Mapakshi_JCP_2018,joshaghani2019composable} 
and the references within can be applied.  
Since the block $\mathbf{K}_{ss}$ is not invertible, we rewrite the system $\mathbf{K}$ as
\[
\mathbf{K}=
    \begin{pmatrix}
        \mathbf{A}_{11} & \mathbf{A}_{12} \\
        \mathbf{A}_{21} & \mathbf{A}_{22}
    \end{pmatrix},
\]
where $\mathbf{A}_{11} = \mathbf{K}_{ss}(1:M-1,1:M-1)$. In other words, we shifted the definition of the blocks
so that the first block $\mathbf{A}_{11}$ is of size $(M-1)\times (M-1)$ and it is now invertible.
We now write
\begin{align}
    \mathbf{K} = 
    \begin{pmatrix}
        \mathbf{I} & \mathbf{0} \\
        \mathbf{A}_{21}(\mathbf{A}_{11})^{-1} & \mathbf{I}
    \end{pmatrix}
    \begin{pmatrix}
        \mathbf{A}_{11} & \mathbf{0} \\
        \mathbf{0} & \mathbf{S}
    \end{pmatrix}
    \begin{pmatrix}
        \mathbf{I} & (\mathbf{A}_{11})^{-1}\mathbf{A}_{12} \\
        \mathbf{0} & \mathbf{I}
    \end{pmatrix},
\end{align}
where $\mathbf{I}$ is the identity matrix and 
\begin{align}
    \mathbf{S} = \mathbf{A}_{22} - \mathbf{A}_{21}(\mathbf{A}_{11})^{-1}\mathbf{A}_{12}
\end{align}
is the Schur complement. The inverse can therefore be written as:
\begin{align}
    \mathbf{K}^{-1} = 
    \begin{pmatrix}
        \mathbf{I} & -(\mathbf{A}_{11})^{-1}\mathbf{A}_{12} \\
        \mathbf{0} & \mathbf{I}
    \end{pmatrix}
    \begin{pmatrix}
        (\mathbf{A}_{11})^{-1} & \mathbf{0} \\
        \mathbf{0} & \mathbf{S}^{-1}
    \end{pmatrix}
    \begin{pmatrix}
        \mathbf{I} & \mathbf{0} \\
        -\mathbf{A}_{21}(\mathbf{A}_{11})^{-1} & \mathbf{I}
    \end{pmatrix}.
\end{align}
The task at hand is to find the inverse of $\mathbf{S}$.
Note that $\mathbf{A}_{11}$ is a diagonal mass matrix for the saturation equation and hence
it is straightforward to obtain the inverse. 
For the Schur complement block we employ the multigrid V-cycle on $\mathbf{S}$ from the HYPRE boomerAMG package 
\citep{falgout2002hypre}.
We expect this to work since the $\mathbf{S}$ block is spectrally equivalent to the Laplacian.
When the inverses are obtained, we rely on GMRES \citep{saad1986gmres} with relative tolerance of $1\times 10^{-8}$
to solve the entire block system.
It is found in \citep{Mapakshi_JCP_2018} that this methodology is computationally less expensive and more practical 
for large-scale computations. 
Solving the system of equations \eqref{Eqn:kxf} 
in fast and efficient way can be done through PETSc \citep{petsc_user_ref,petsc_web_page,Dalcin2011} 
and its composable solver capabilities \citep{brown2012composable}.
Appendix~\ref{Sec:Appendix} contains the  necessary PETSc command-line options for the described Schur complement approach.
All the numerical results are generated using FEniCS Project
\citep{LoggWells2010a,brooks1964hydrau}. Among the various components available
in FEniCS, we use the DOLFIN library \citep{LoggWellsEtAl2012a}
and the Unified From Language library \citep{Alnaes2012a}.
Simulations are conducted on a single socket Intel Core i7-7920HQ server node 
by utilizing a single MPI process.
Computer codes implementing the proposed computational framework can be found at 
\citep{zenodo/Vertex_based_method}.

%
\section{REPRESENTATIVE NUMERICAL RESULTS}
\label{Sec:numer}
\subsection{Analytical problem and $h-$convergence study}
We first perform an $h$-convergence study on two-dimensional structured triangular meshes of size $h$.
Consider a unit square to be the computational domain with the 
following expressions for the saturation and pressure fields:
\begin{subequations}
    \label{Eqn:analytical}
\begin{align}
&s(x,y,t) = 0.4+0.4xy+0.2\cos(t+x),\\
&p(x,y,t) = 2+x^2y-y^2+x^2\sin(y+t)-\frac13 \cos(t) + \frac13 \cos(t+1) -\frac{11}{6}.
\end{align}
\end{subequations}
We replace the source/sink terms (i.~e.,~wells flow rates) of equations
\eqref{Eqn:GE_Wphase}--\eqref{Eqn:GE_pressure} by functions denoted by $f_1$ and $f_2$, obtained
via the method of manufactured solutions..
%
Dirichlet boundary conditions are applied on $\partial \Omega$ on both saturation and 
pressure fields. The input parameters are:
\[
\phi=2, \, K=1, \, \mu_w=\mu_o=1,\, s_{rw}=s_{ro}=0, \,  k_{rw}(s)=s^2, \, k_{ro}(s)=(1-s)^2.
\]
The capillary pressure satisfies \eqref{eq:pc} with $\theta=2$, $p_d=50$,  and $R=0.05$.
Table \ref{tab:convergence} shows the errors in $L^2$ and $H^1$ norms evaluated at $T=1$ and the
corresponding convergence rates for saturation and pressure. 
The rates are optimal in the $H^1$ norm. The suboptimal rate in the $L^2$ norm
is expected as first order Taylor expansion is used for capillary pressure, 
and phase mobilities are evaluated through Picard's iterations.
%
The vertex scheme results in the theoretical convergence rate of one for both unknowns,
which confirms the correct behavior of the algorithm.
\begin{table}[]
\caption{Results of convergence test where the mesh size is denoted by $h$.
The time step $\tau$ is set to mesh size and $L^2$ and $H^1$ norms are computed at the
final time $T=1$}
\label{tab:convergence}
\scalebox{0.88}{
\begin{tabular}{lllllllllll}
\Xhline{2\arrayrulewidth}
 &
   &
   &
  \multicolumn{2}{l}{$||s_h-s||_{L^2(\Omega)}$} &
  \multicolumn{2}{l}{$||p_h-p||_{L^2(\Omega)}$} &
  \multicolumn{2}{l}{$||s_h-s||_{H^1(\Omega)}$} &
  \multicolumn{2}{l}{$||p_h-p||_{H^1(\Omega)}$} \\ \cline{4-11} 
\multirow{-2}{*}{$h$} &
\multirow{-2}{*}{$M$} &
  \multirow{-2}{*}{$\tau$} &
  Error &
  Rate &
  Error &
  Rate &
  Error &
  Rate &
  Error &
  Rate \\ \hline
1/4 &
  25 &
1/4 &
  \cellcolor[HTML]{C0C0C0} $9.430\times10^{-4}$ &
  \cellcolor[HTML]{C0C0C0} $-$ &
  $8.830\times10^{-3}$ &
  - &
  \cellcolor[HTML]{C0C0C0} $5.160\times10^{-3}$ &
  \cellcolor[HTML]{C0C0C0}- &
  $4.980\times10^{-2}$ &
  - \\
1/8 &
  81 &
1/8 &
  \cellcolor[HTML]{C0C0C0} $6.600\times10^{-4}$ &
  \cellcolor[HTML]{C0C0C0}0.515 &
  $4.740\times10^{-3}$ &
  0.899 &
  \cellcolor[HTML]{C0C0C0} $3.610\times10^{-3}$ &
  \cellcolor[HTML]{C0C0C0}0.514 &
  $2.610\times10^{-2}$ &
  0.934 \\
1/16 &
  289 &
1/16 &
  \cellcolor[HTML]{C0C0C0} $3.650\times10^{-4}$ &
  \cellcolor[HTML]{C0C0C0}0.853 &
  $2.370\times10^{-3}$ &
  1.000 &
  \cellcolor[HTML]{C0C0C0} $2.010\times10^{-3}$ &
  \cellcolor[HTML]{C0C0C0}0.846 &
  $1.300\times10^{-2}$ &
  1.003 \\
1/32 &
  1069 &
1/32 &
  \cellcolor[HTML]{C0C0C0} $1.890\times10^{-4}$ &
  \cellcolor[HTML]{C0C0C0}0.949 &
  $1.170\times10^{-3}$ &
  1.014 &
  \cellcolor[HTML]{C0C0C0} $1.040\times10^{-3}$ &
  \cellcolor[HTML]{C0C0C0}0.944 &
  $6.440\times10^{-3}$ &
  1.014 \\
1/64 &
  4225 &
1/64 &
  \cellcolor[HTML]{C0C0C0} $9.350\times10^{-5}$ &
  \cellcolor[HTML]{C0C0C0}1.018 &
  $5.500\times10^{-4}$ &
  1.094 &
  \cellcolor[HTML]{C0C0C0} $5.220\times10^{-4}$ &
  \cellcolor[HTML]{C0C0C0}1.000 &
  $3.270\times10^{-3}$ &
  0.975 \\ 
\Xhline{2\arrayrulewidth}
\end{tabular}
}
\end{table}

\subsection{Physical problems}
In this section, robustness of the proposed vertex scheme is assessed using standard two- and three-dimensional
test problems. Numerical responses of several five spot and quarter-five spot problems, with homogeneous and 
heterogeneous permeability fields are investigated.
We examine the element-wise mass balance property associated with the vertex scheme and 
also comment on capability of the scheme in satisfying the maximum principle.
Let water and oil be the wetting phase and non-wetting phase, respectively. For
all problems, the relative permeability and capillary data satisfy \eqref{eq:krwo}, \eqref{eq:pc} and we assume the following:
\begin{align}
    &s_{\mathrm{in}} = 0.85, \quad 
    \mu_{w} = 5 \times 10^{-4} \mbox{~kg/ms}, \quad
    \mu_{o} = 2 \times 10^{-3} \mbox{~kg/ms}, \quad \\
    &\phi = 0.2, \quad
    p_d = 5 \times 10^{3} \mbox{~Pa}, \quad 
    \theta = 3, \quad 
    s_{rw} = s_{ro} = 0.15, \\
    & s^{0}=0.15, \quad p^{0}=1\times{10}^{6} \mbox{~Pa}.
\end{align}
\subsubsection{Two-dimensional homogeneous medium}
\label{sub:2D_homogen}
We take a domain of $\Omega = [0,100]^2\;\mathrm{m^2}$ 
with mesh-size of $h=100/40$~m. No-flow boundary condition over $\partial \Omega$ is chosen for 
this problem (see Figure \ref{Fig:2D_BVP}) and flow is driven from the injection to the production wells by
introducing source and sink terms.
\begin{figure}
  \centering
  \sbox{\measurebox}{%
    \begin{minipage}[b]{.5\textwidth}
      \subfigure
	  [Schematic]
	  {\label{Fig:2D_BVP}\includegraphics[scale=0.395]{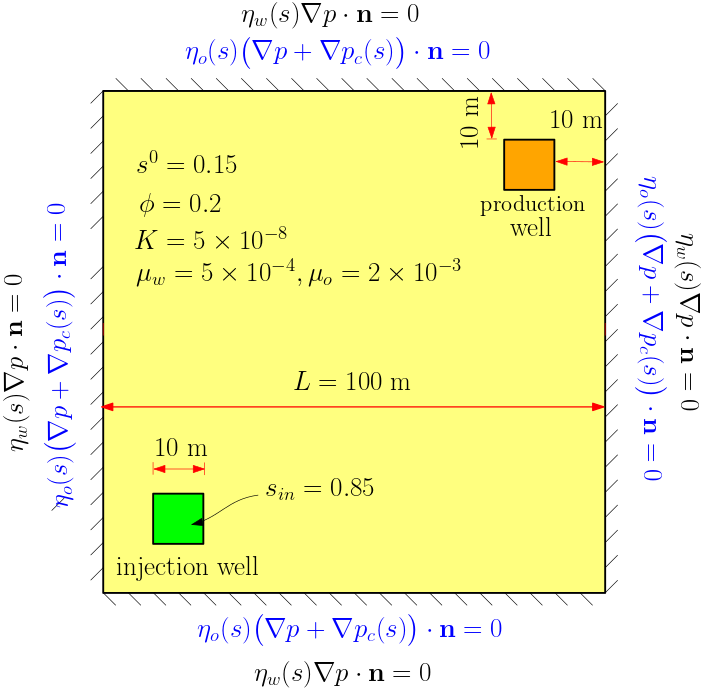}}
  \end{minipage}}
  \hspace{0.6cm}
  \usebox{\measurebox}\qquad
  \begin{minipage}[b][\ht\measurebox][s]{.4\textwidth}
    \centering
    \subfigure
	[Structured mesh]
    {\label{Fig:Mesh_structured}\includegraphics[scale=0.23,trim=0 0cm 0 0]{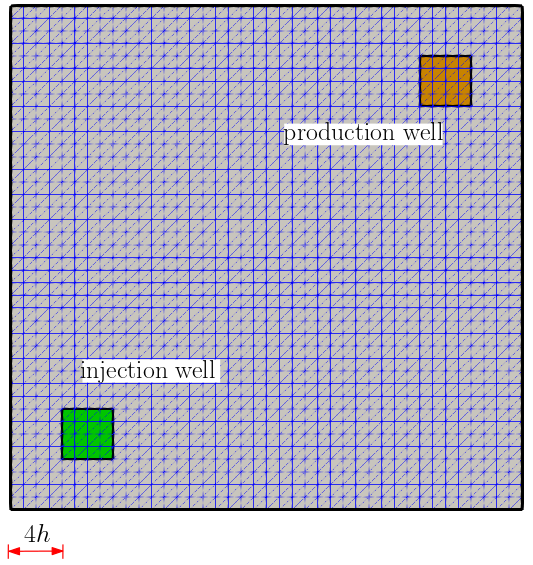}}		
	\vfill		
	\subfigure
	    [Unstructured mesh]
	    {\label{Fig:Mesh_unstructured}\includegraphics[scale=0.18]{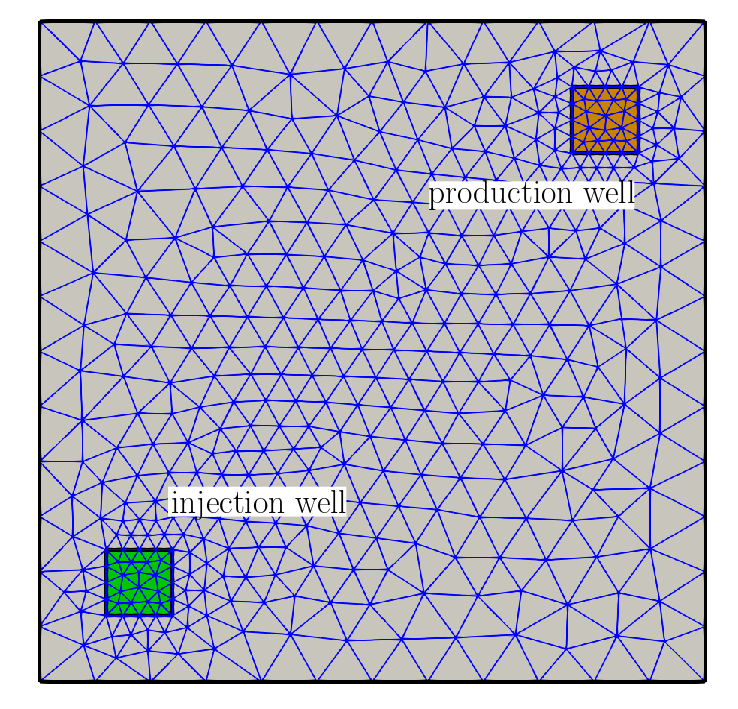}}
  \end{minipage}
  \caption{\textsf{Two-dimensional homogeneous medium:} 
    This figure provides a pictorial description
    of the boundary value problem and shows the
    typical meshes employed in our numerical
simulations. \label{Fig:2D_schematic}}
\end{figure}
The injection and production flow rates satisfy:
\begin{align}
    \int_{\Omega}\bar{q}=\int_{\Omega}\underline{q}=0.1,
\end{align}
where $\bar{q}$ is piecewise constant on $[10,20]\;\mbox{m}\times[10,20]\;\mbox{m}$
and $\bar{q}=0$ elsewhere and $\underline{q}$ is piecewise constant
on $[80,90]\;\mbox{m}\times[80,90]\;\mbox{m}$ and $\underline{q}=0$ elsewhere. 
We choose a constant permeability $K= 5 \times 10^{-8}~\mathrm{m}^2$. 
Domain is discretized with a triangular structured mesh and the time step is $\tau = 60$ s.
The final simulation time is $T=12000$ s, and we provide solutions snapshots at $t=1800$ s, $t=6000$ s, and $t=12000$ s. 
The saturation and pressure profiles obtained under vertex scheme are, respectively, displayed in Figures
\ref{Fig:sat_2Dhom_30}--\ref{Fig:sat_2Dhom_200} and \ref{Fig:p_2Dhom_30}--\ref{Fig:p_2Dhom_200}.
The wetting phase is injected at the lower left end of the domain, and displaces the non-wetting fluid to 
the upper right corner. Note that this problem is convection-dominated but
it is evident that numerical saturation  remain within physical bounds ($s_h^n\geq 0.15$ and $s_h^n\leq 0.85$) during simulation
and no undershoot and overshoot are observed.
It only takes $4$ to $5$ Picard's iterations at each time step for convergence of the vertex scheme. 
This is true for all two-dimensional test cases unless specified otherwise.
\begin{figure}
	\subfigure[$t = 1800$ s\label{Fig:sat_2Dhom_30}]{
		\includegraphics[clip,scale=0.18,trim=0 0cm 0cm 0]{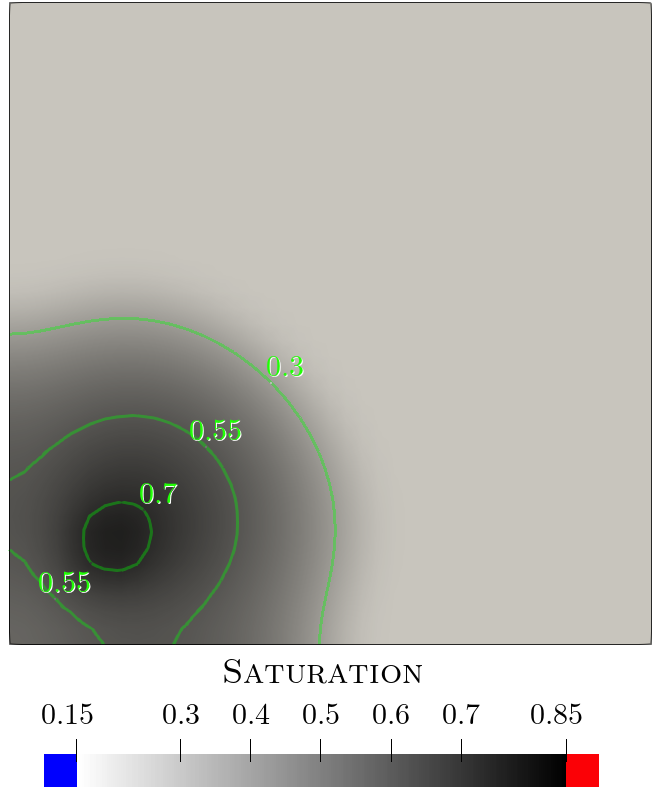}} 
        \hspace{.25cm}
	\subfigure[$t = 6000$ s \label{Fig:sat_2Dhom_100}]{
		\includegraphics[clip,scale=0.18,trim=0 0cm 0cm 0]{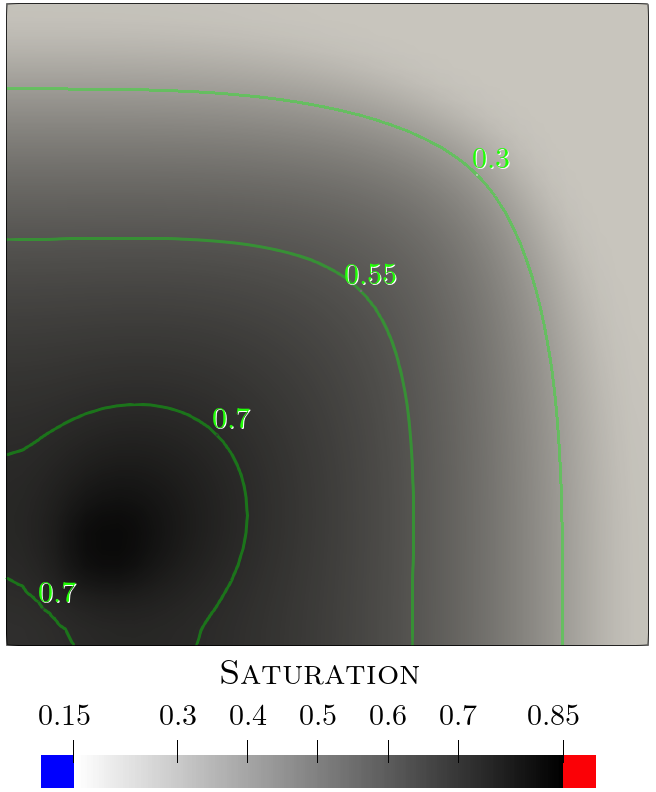}} 
        \hspace{.25cm}
    \subfigure[$t = 12000$ s\label{Fig:sat_2Dhom_200}]{
		\includegraphics[clip,scale=0.18,trim=0 0cm 0cm 0]{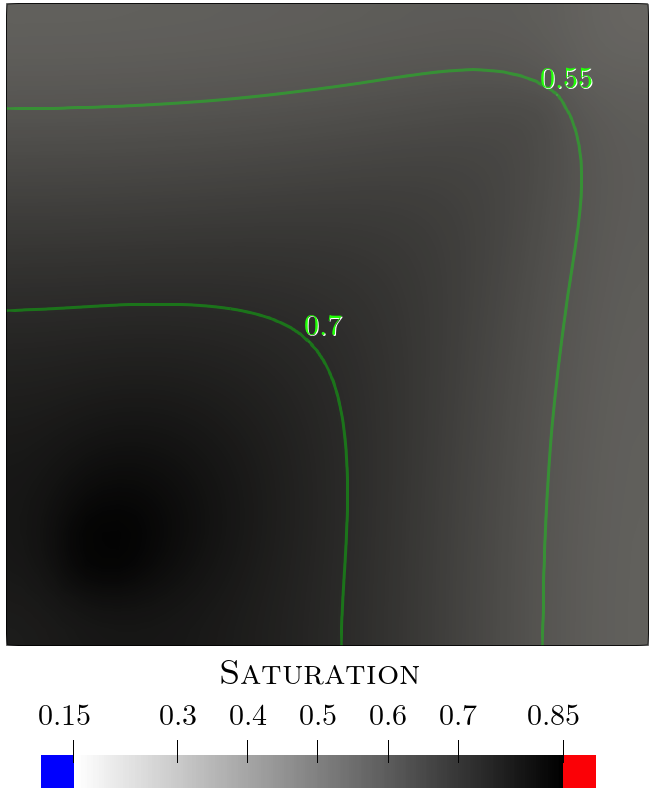}} \\
	\subfigure[$t = 1800$ s\label{Fig:p_2Dhom_30}]{
		\includegraphics[clip,scale=0.1,trim=0 0 0 10cm]{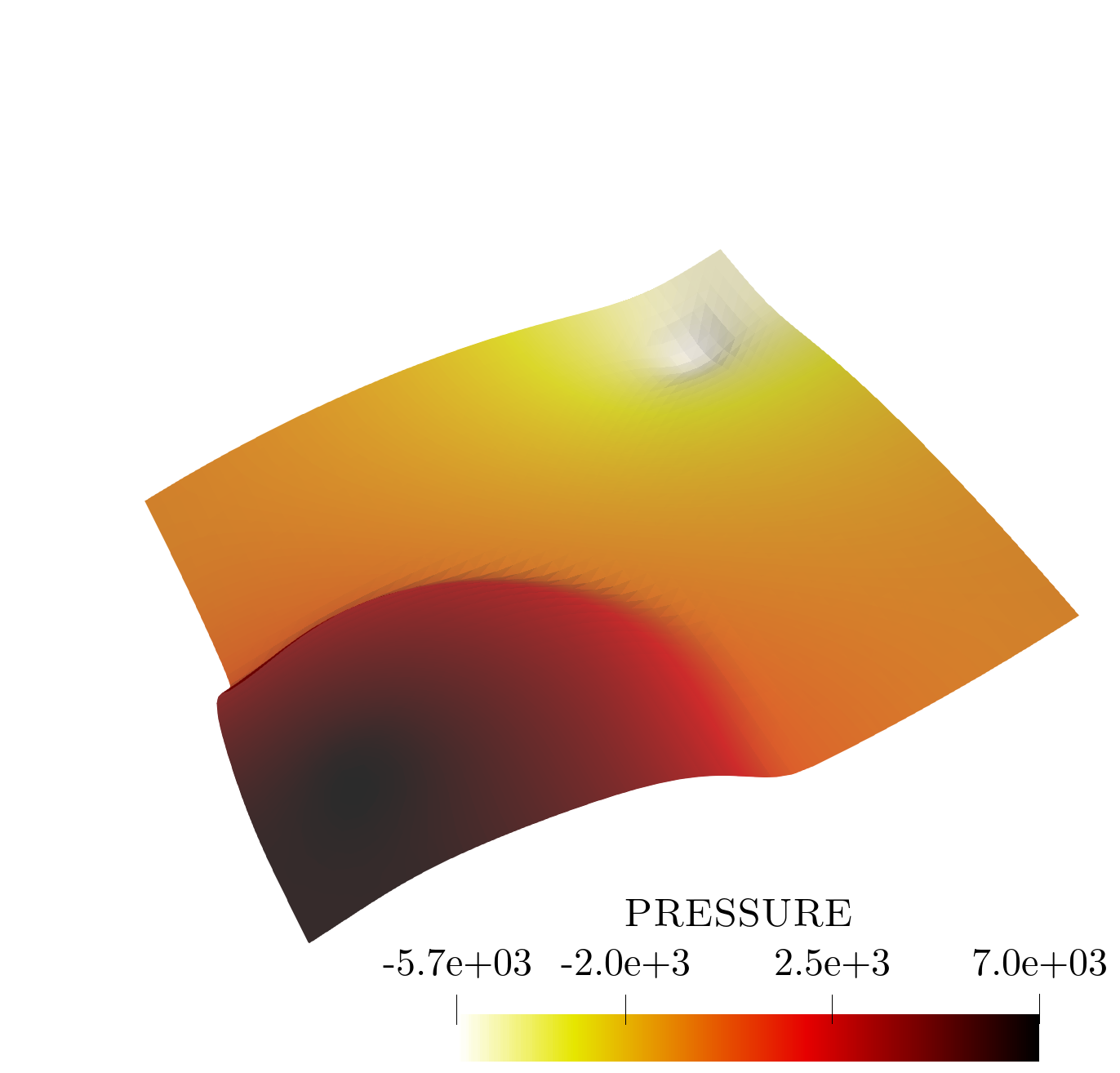}} 
        \hspace{.25cm}
	\subfigure[$t = 6000$ s \label{Fig:p_2Dhom_100}]{
		\includegraphics[clip,scale=0.1,trim=0 0cm 0cm 10cm]{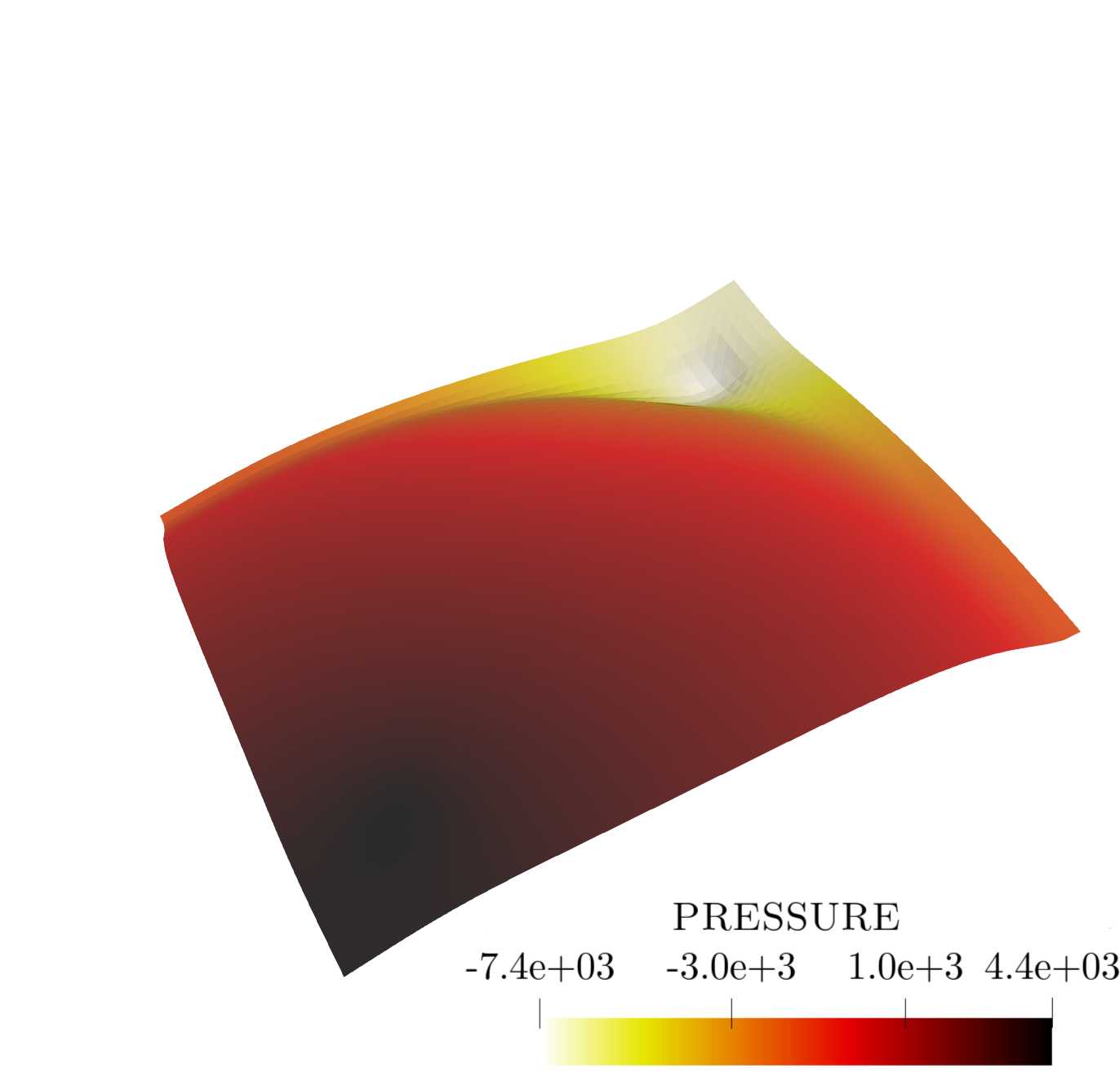}} 
        \hspace{.25cm}
    \subfigure[$t = 12000$ s \label{Fig:p_2Dhom_200}]{
		\includegraphics[clip,scale=0.1,trim=0 0cm 0cm 10cm]{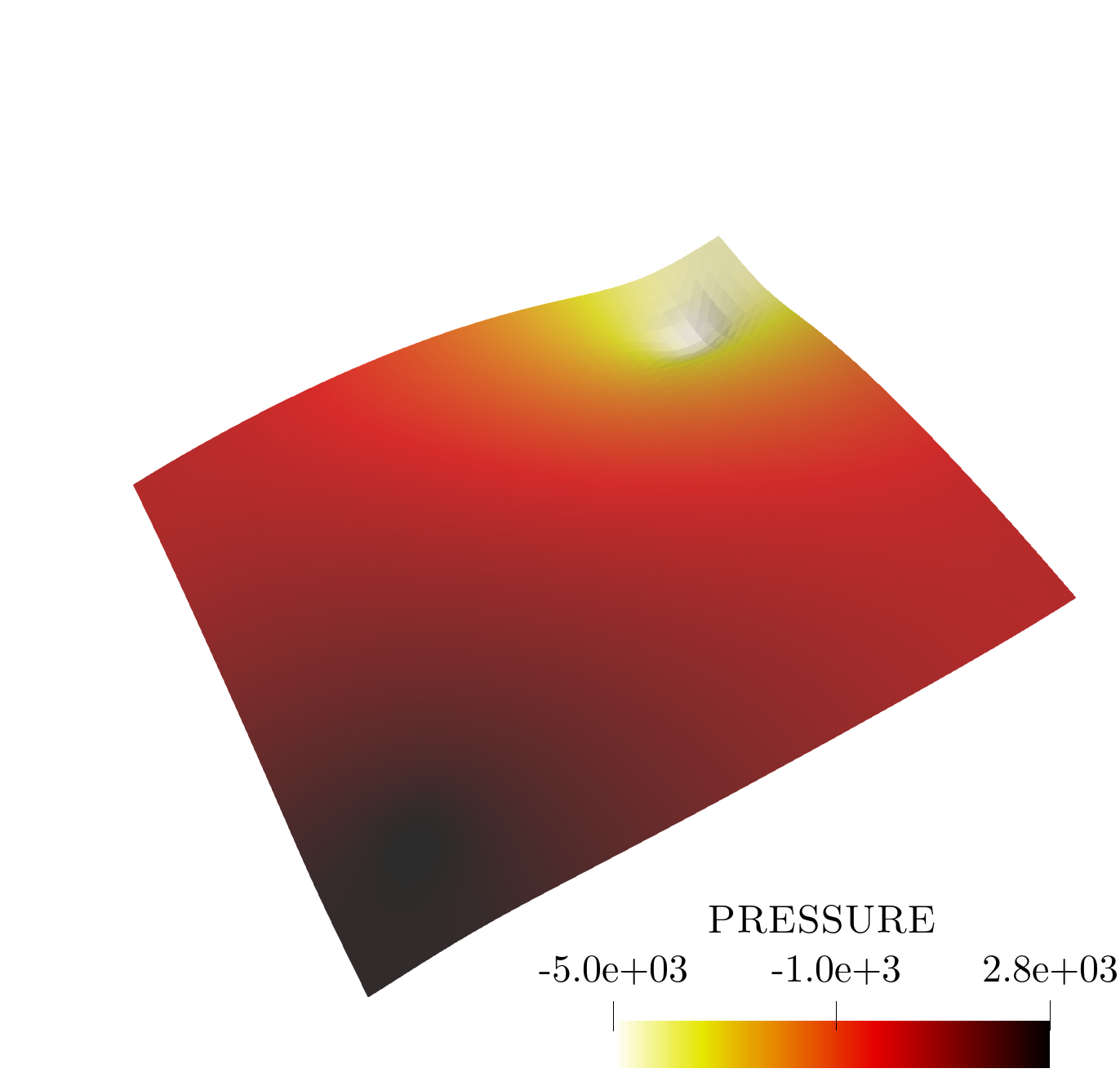}} 
        \caption{\textsf{Two-dimensional homogeneous medium:}
                This figure shows the saturation and pressure profiles obtained from the proposed vertex scheme  
                at three different time steps. Structured triangular mesh is used (see Figure~\ref{Fig:Mesh_structured}).
                The saturation field (top figures) should be between $0.15$ and $0.85$.
                This figure suggests that the proposed scheme is capable of
                providing maximum-principle satisfying results.
                To wit, no undershoots (blue-colored cell) or overshoots (red-colored cell) observed.
                Highest value and lowest value of pressure solutions (bottom figures)
                detected at the injection well and the production well, respectively. The pressure difference forces the
                wetting phase flow through the domain. Pressure differences subside as the front reaches the production well.
                Pressure solutions are warped for better visualization.
        \label{Fig:2D_homogen}}
\end{figure}

We compare saturation profiles obtained from the vertex scheme
with the solutions obtained from the fully implicit discontinuous 
Galerkin (DG) formulation developed by \citet{epshteyn2007fully}. 
For the chosen DG formulation polynomial order is set to $\mathbb{P}=1$, 
DG symmetry parameter is set to $\epsilon = +1$ (i.e., NIPG), and the penalty parameter is set to $\sigma = 0.1$.
Both DG formulation and vertex scheme are solved on the structured triangular mesh 
(as shown in Figure~\ref{Fig:Mesh_structured}) and the time step is $\tau=60$ s.

\begin{figure}
	\subfigure[Saturation profile  \label{Fig:vertex_vs_DG_sat}]{
    \includegraphics[clip,width=0.45\linewidth,trim=0 0 7.5cm 0]{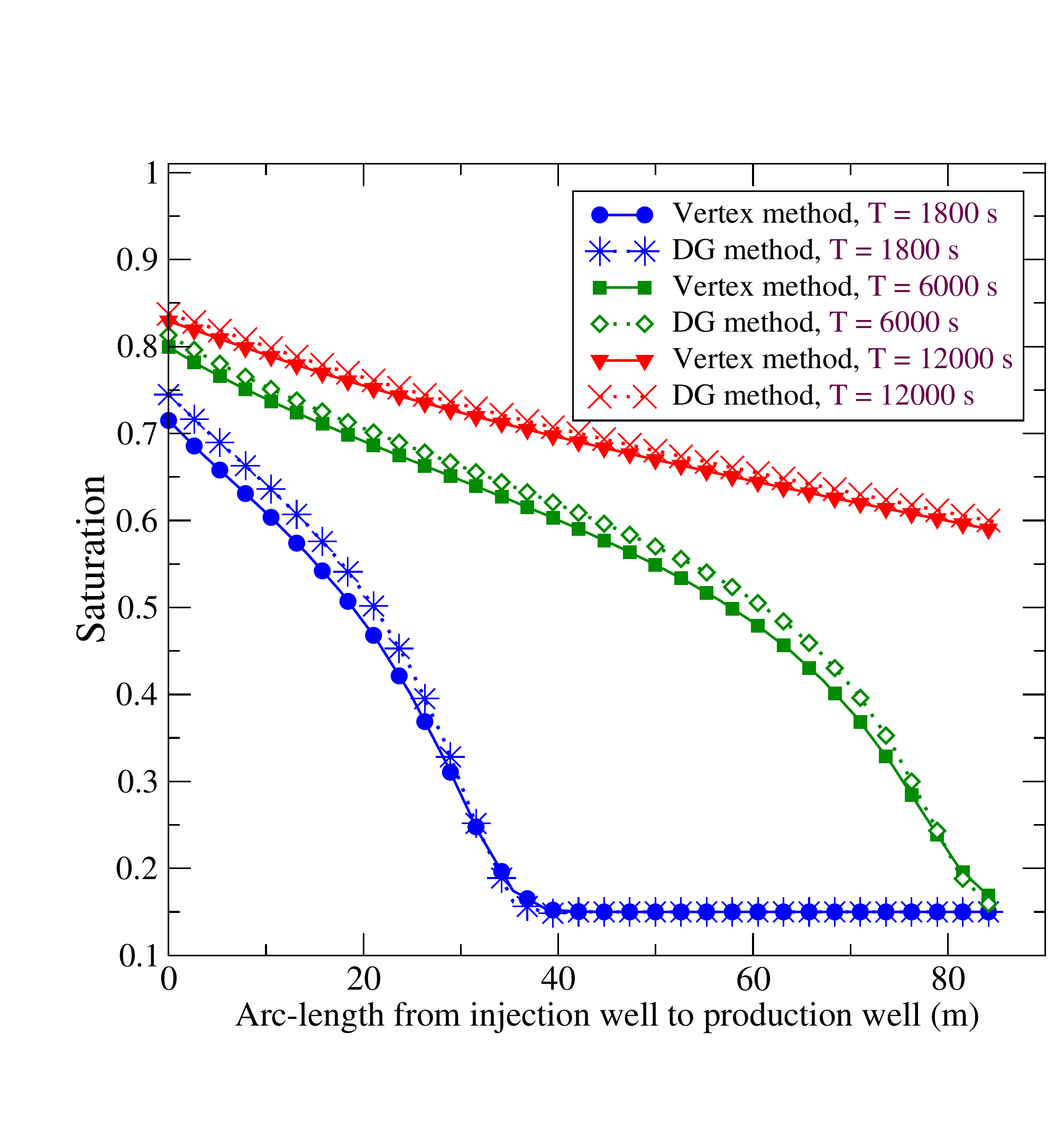}}
    \hspace{1cm}
	\subfigure[Pressure profile \label{Fig:vertex_vs_DG_pres}]{
    \includegraphics[clip,width=0.45\linewidth,trim=0 0 7.5cm 0]{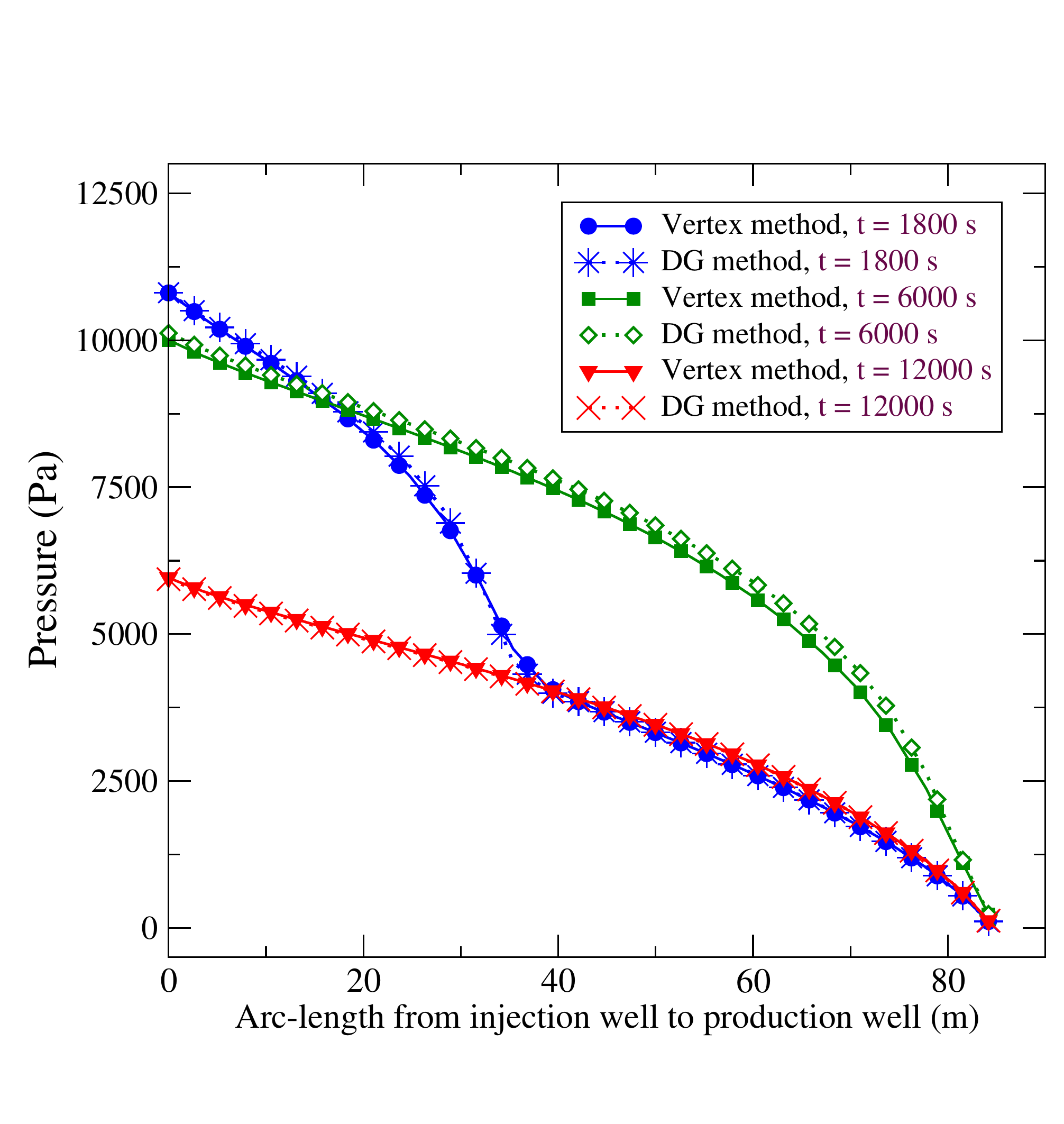}}\\
  %
  \caption{
      \textsf{Two-dimensional homogeneous medium:}
      This figure compares the saturation and pressure profiles obtained from the vertex scheme and a
      first order discontinuous Galerkin formulation (i.e., DG-NIPG with $\mathbb{P}=1$ and $\sigma=0.1$).
      Solutions are plotted along the diagonal line spanned from point $(20,20)$~m to point $(80,80)$~m.
      DG and vertex scheme reproduce similar responses throughout the simulation.
      Given that our solutions are only plotted on the diagonal line, the results from the left figure should be treated 
      with considerable caution.
      As highlighted in Figure~\ref{Fig:2D_homogen}, saturation solutions under the vertex scheme are always bounded by physical
      values. However, we are aware the DG method (and finite element methods in general) do not enjoy maximum principle.
      For example, in this experiment the lowest and highest value of saturation found under the DG-NIPG formulation 
      was $0.148$ and $1.185$, respectively. 
    \label{Fig:vertex_vs_DG}
  }
\end{figure}
For three representative time steps, the saturation 
and pressure profiles, along the diagonal $\{ (x,y):x=y \}$ from the injection well upto production well 
are illustrated in Figures~\ref{Fig:vertex_vs_DG_sat} 
and~\ref{Fig:vertex_vs_DG_pres}. \\
We observe that
the finite element solutions are accurate and in very good agreement with the DG solutions.
It can be seen that the saturation fronts, under both DG and proposed vertex scheme, 
propagates with the same speed. 
We recall that the proposed finite element scheme satisfies a maximum principle, as mentioned in Proposition~\ref{prop:conv}.
However, the DG approximations of the saturation are not guaranteed to satisfy \eqref{eq:maxsat}
and small undershoot (usually at the injection well) and overshoots 
(usually right after the saturation front) are observed for saturation profile.
Increasing DG polynomial order (in addition to sharpening front) is reported to relatively
reduce these unphysical violations \citep{epshteyn2009analysis}.
Even so, DG schemes still require external bound-preserving mechanisms
such as slope/flux limiting \citep{kuzmin2010vertex},
artificial viscosity \citep{van2002space},
or nodal-based optimization \citep{joshaghani2020modeling} to completely enforce maximum principle. 
A comprehensive survey of bound-preserving methods is described in \citep{zhang2011maximum}.

\subsubsection{Conservation of local mass balance}
Next, we investigate the local mass conservation property of the proposed scheme for the incompressible two-phase flow model. 
The local mass conservation of an element $\omega$ at each time step, is calculated as follows:
\begin{align}
    m(E) =  \int_{E} \frac{\phi (s_h^{n}-s_h^{n-1}) }{\tau} 
    - \int_{\partial E} \eta_{w}(s^{n}) K_E \nabla p^{n} \cdot \mathbf{n}_{E}
    - \int_{E} \left(f_w(s_{\mathrm{in}}) \bar{q} + f_w(s^{n}) \underline{q}\right).
\end{align}
A true locally mass conservative scheme should produce the zero value for $m(E)$ for each element $E$.
We compute the mass balance values for the problem described in Section~\ref{sub:2D_homogen}. 
In Figure~\ref{Fig:BoM}, the values of $m(E)$ are displayed at three representative time steps.
One can see that the magnitude of $m(E)$ in the domain (except at the wells' locations) is always less than
$10^{-5}$, which is the tolerance set for the Picard's iteration. Hence, the proposed scheme
is locally mass conservative.
We note that the source/sink models result in a higher mass balance error values (of the order of $10^{-3}$) on the elements
that form the support of the  injection and production wells.
\begin{figure}
	\subfigure[$t=1800$ s  \label{Fig:BoM_vertex_1800}]{
		\includegraphics[clip,scale=0.17,trim=0 0cm 12cm 0]{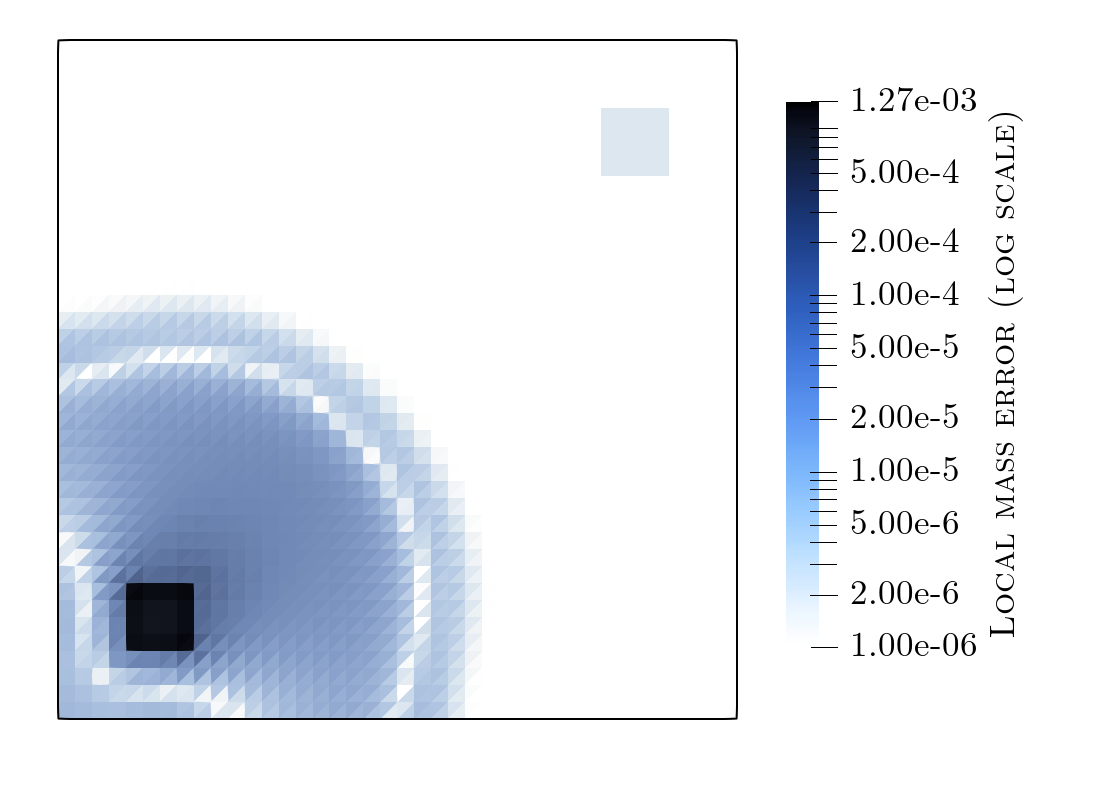}}
	\hspace{-0.35 cm}
	\subfigure[$t=6000$ s \label{Fig:BoM_vertex_6000}]{
		\includegraphics[clip,scale=0.17,trim=0 0cm 12cm 0]{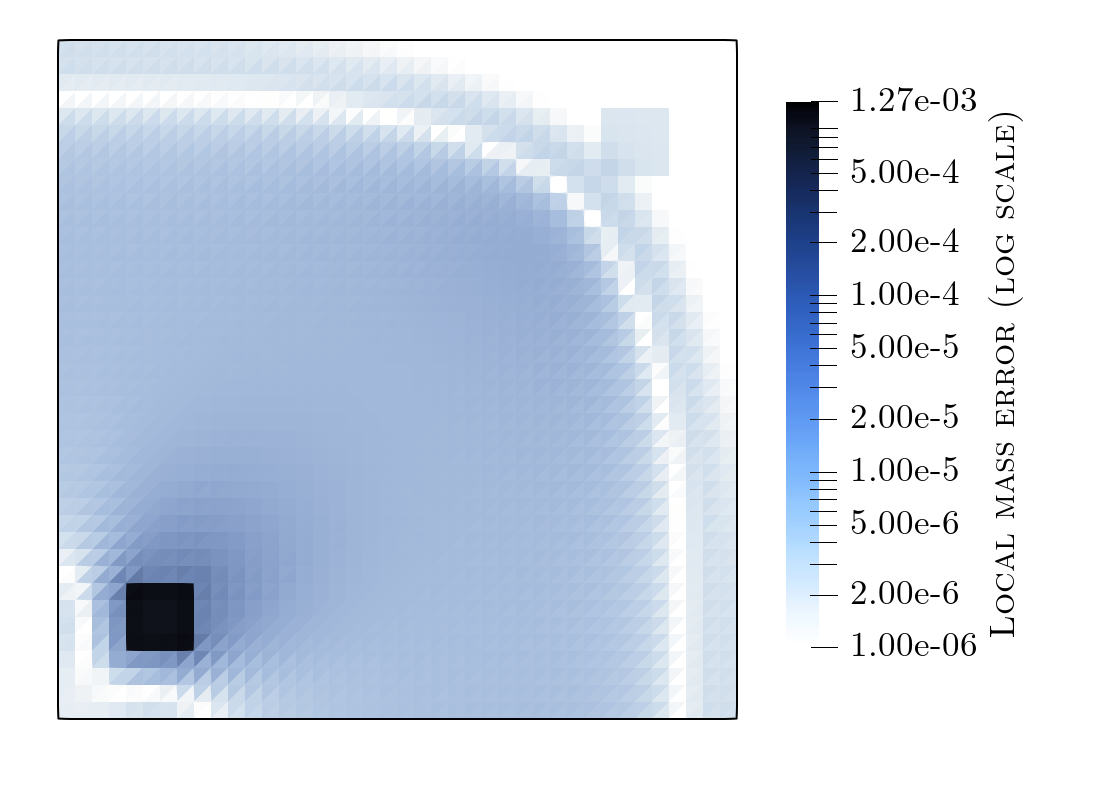}}
	\hspace{-0.35 cm}
	\subfigure[$t=12000$ s  \label{Fig:BoM_vertex_12000}]{
		\includegraphics[clip,scale=0.17,trim=0 0cm 12cm 0]{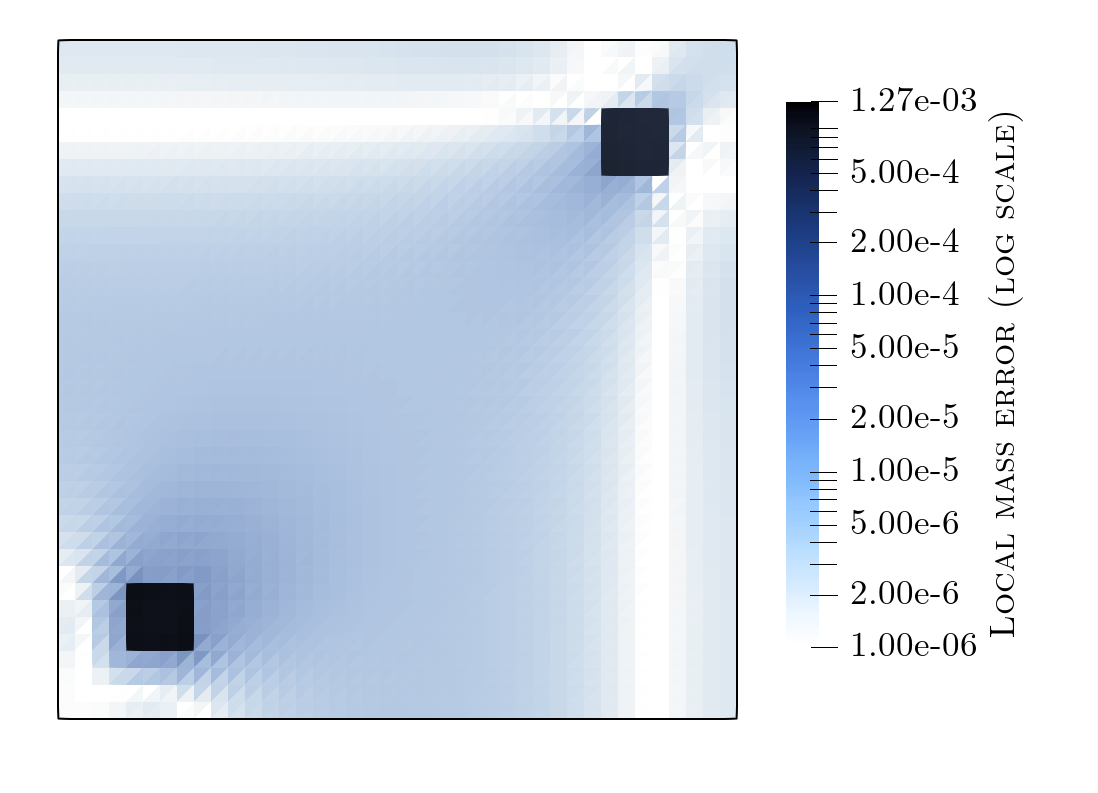}}
	\hspace{0.1 cm}
	\subfigure{
		\includegraphics[clip,scale=0.18,trim=1cm 3cm 2cm 0cm]{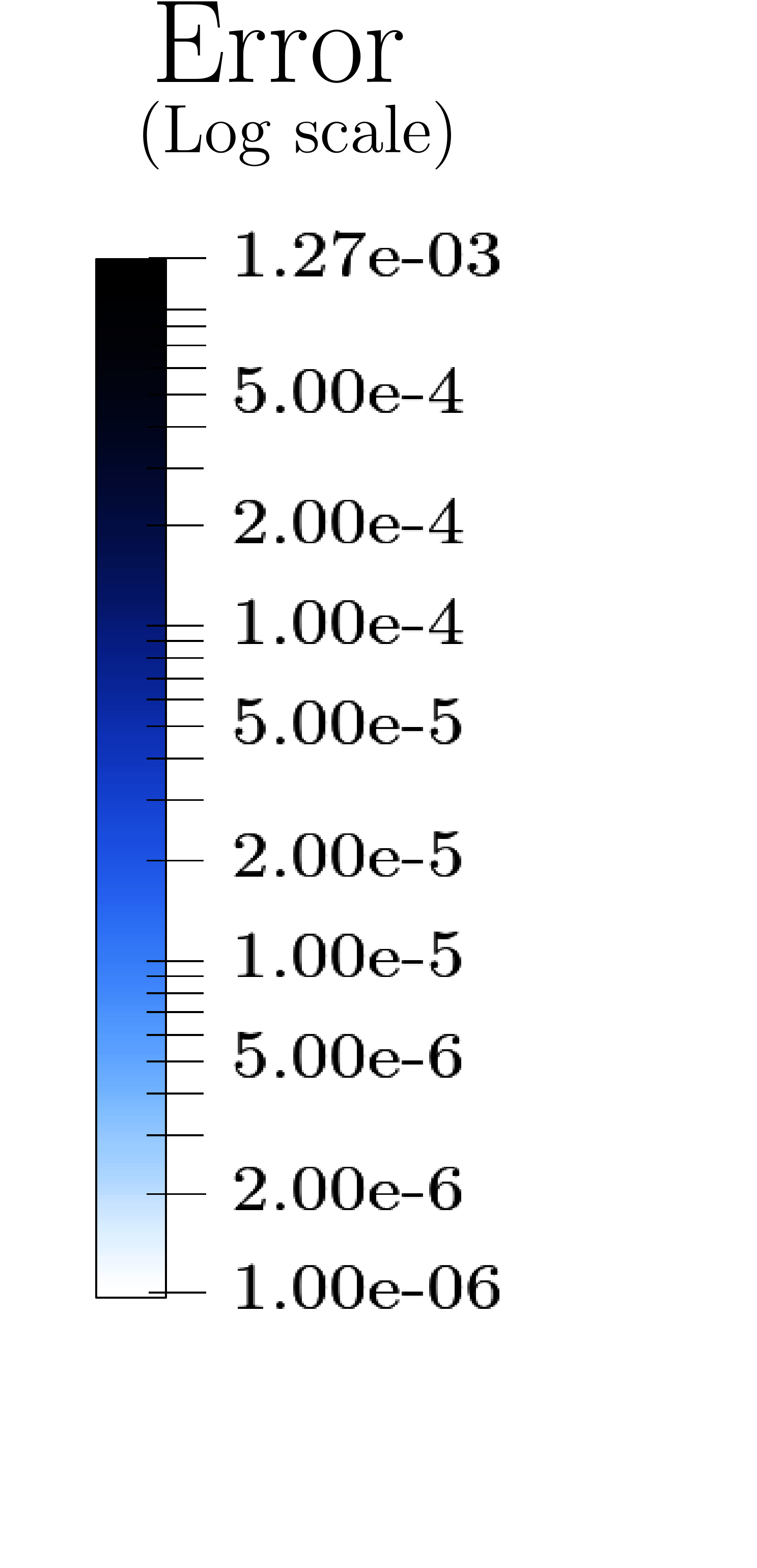}}
        \caption{
            \textsf{Local mass balance conservation:}
            This figure illustrates the local mass conservation properties of the vertex scheme 
            for two-phase incompressible flow problem on a homogeneous domain.
            The mass balance error remains small as time advances.
            This value is always less than $10^{-5}$ inside the reservoir (excluding elements allocated
            to wells).
            \label{Fig:BoM}
        }
\end{figure}

\subsubsection{Two-dimensional domain with unstructured mesh}
All parameters are the same as in Section~\ref{sub:2D_homogen}, 
except for the mesh that is triangular unstructured as depicted in Figure \ref{Fig:Mesh_unstructured}.
Figure \ref{Fig:2D_Unstructured} shows the saturation profiles at four different time steps.
We observe that the saturation remains bounded and no violations of maximum principle 
are observed throughout the simulation. This result also shows that the proposed finite element  
scheme handles unstructured meshes as expected. 
\begin{figure}
	\subfigure[$t=1500$ s  \label{Fig:BoM_vertex_1800}]{
		\includegraphics[clip,scale=0.125,trim=0 0cm 0 0]{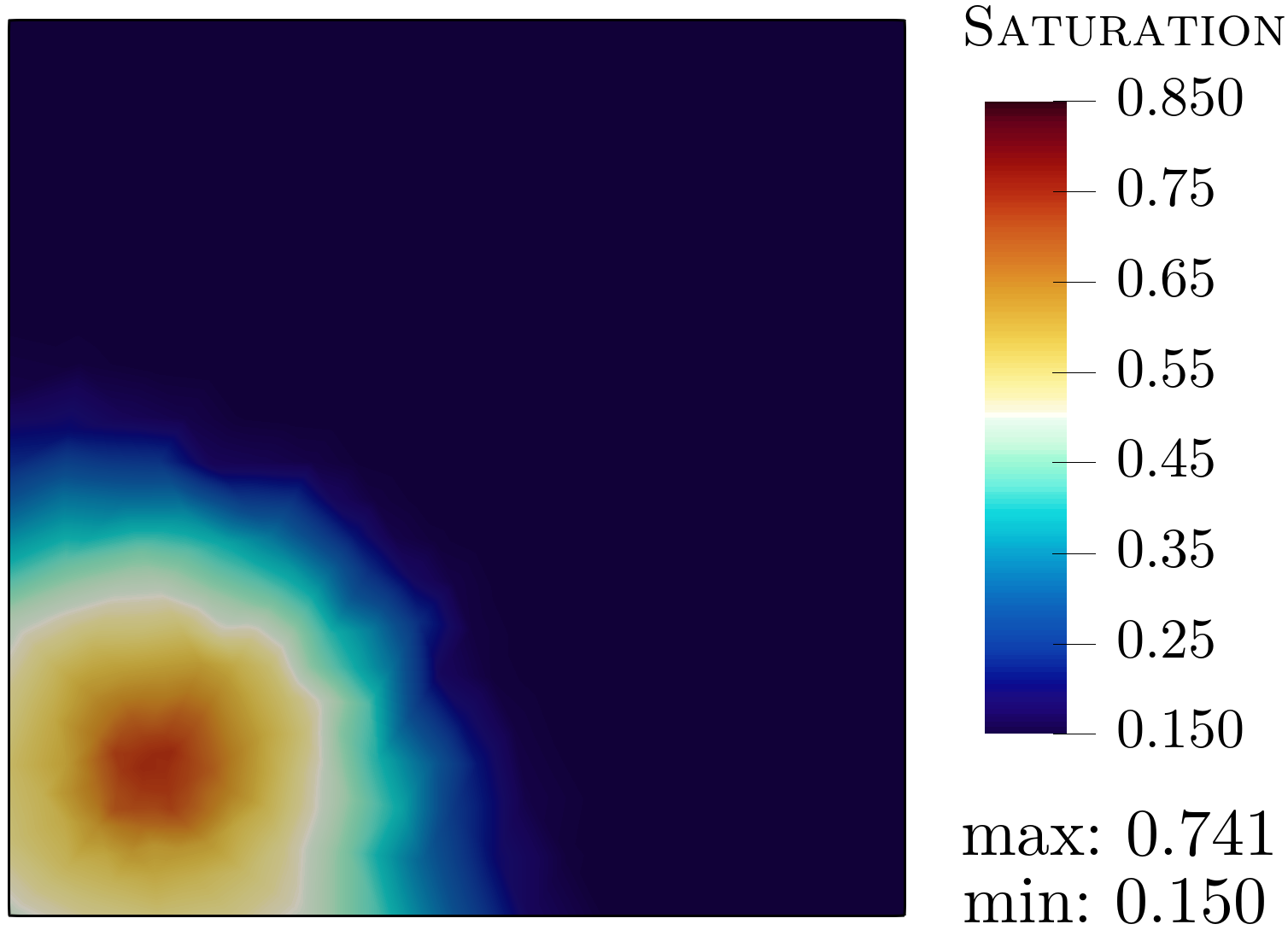}}
        \hspace{0.5cm}
	\subfigure[$t=4500$ s  \label{Fig:BoM_vertex_1800}]{
		\includegraphics[clip,scale=0.125,trim=0 0cm 0 0]{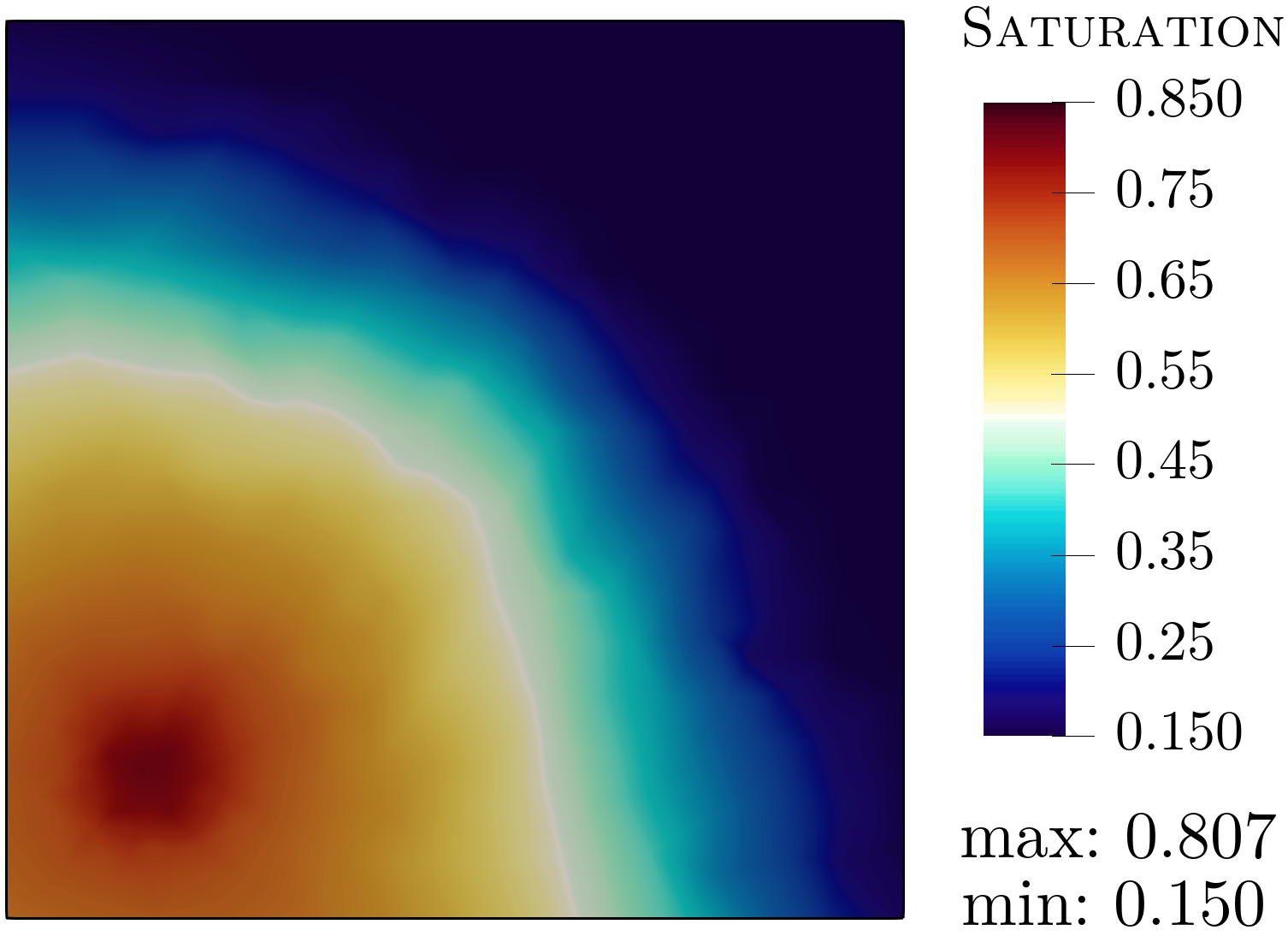}}
        \vspace{0.75cm} 
	\subfigure[$t=6900$ s  \label{Fig:BoM_vertex_1800}]{
		\includegraphics[clip,scale=0.125,trim=0 0cm 0 0]{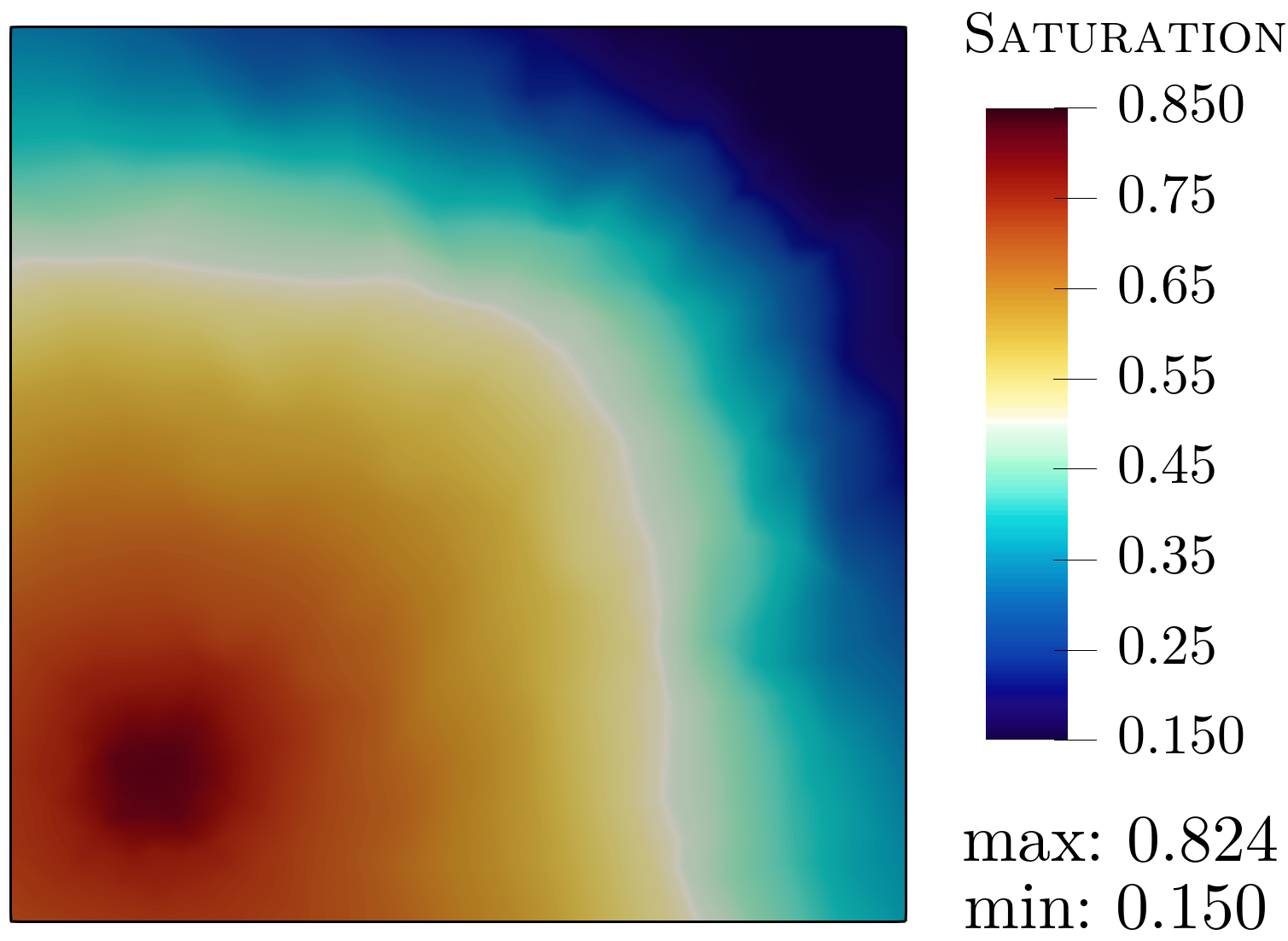}}
        \hspace{0.58cm}
	\subfigure[$t=12000$ s  \label{Fig:BoM_vertex_1800}]{
		\includegraphics[clip,scale=0.175,trim=0 0cm 0 0]{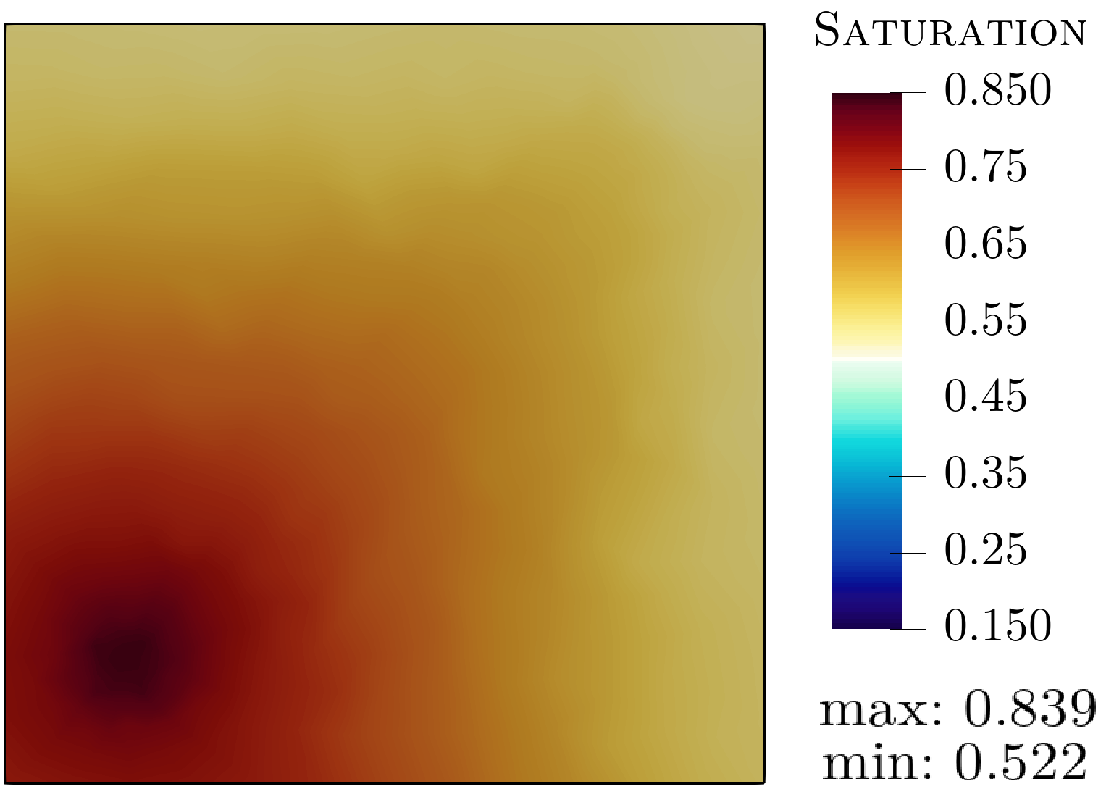}}
        \caption{
            \textsf{Two-dimensional homogeneous domain with unstructured mesh}:
              This figure shows saturation contour plots obtained  on an unstructured triangular mesh
              (Figure~\ref{Fig:Mesh_unstructured}).
              These results pinpoint that the proposed scheme can provide accurate and maximum-principle satisfying
              results on unstructured meshes.
        \label{Fig:2D_Unstructured}}
\end{figure}
\vspace{-\baselineskip}
\subsubsection{Two-dimensional porous media with low permeability block}
In this problem, the domain is $\Omega=[0,100]^2\;\mathrm{m^2}$, 
and as shown in Figure~\ref{Fig:heterogen_schematic} the permeability is 
$5\times 10^{-8}$ everywhere except inside the square inclusion of size $L=20$ m, where the permeability $k_\mathrm{In}$ is 10 
times smaller. The remaining parameters are the same as in Section~\ref{sub:2D_homogen}.
\begin{figure}
  \includegraphics[clip,width=0.6\linewidth]{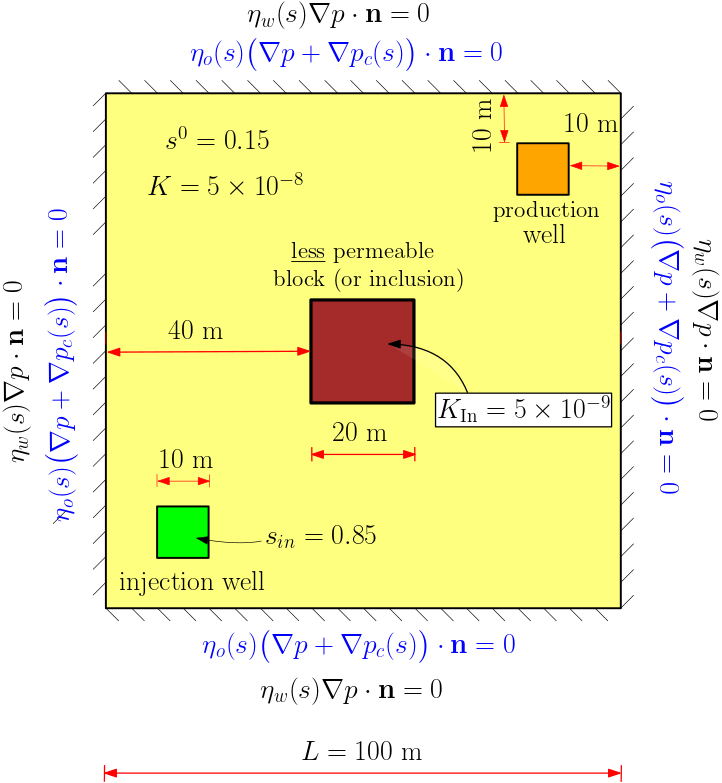}
  \caption{
      \textsf{Two-dimensional porous media with low permeability block:}
      This figure shows the representative computational domain and the boundary value problem. No flow boundary conditions
      are assigned along all boundaries.
  \label{Fig:heterogen_schematic}}
\end{figure}
The saturation solutions at different time steps are depicted in Figure~\ref{Fig:sat_heterogen}.
As expected, the wetting phase initially avoids the region of lower permeability, while still traveling towards the production well.
Toward the end of the simulation, we can observe that the wetting phase has started to penetrate the inclusion region. 
However, when we increase order of difference in permeabilities (e.g., $K/K_{\mathrm{In}}=10000$), the inclusion becomes 
impenetrable throughout the simulation. This trend is clearly shown in Figure~\ref{Fig:sat_two_case} and was reported in the literature 
for single-phase flow \citep{li2015numerical}, and two-phase flow \citep{fabien2020high}.
Figures \ref{Fig:sat_heterogen} and \ref{Fig:sat_two_case} also highlight that the vertex scheme (on a relatively coarse mesh) 
is capable of generating sharp saturation front in the domain with non-homogeneous permeability, 
in addition to completely suppressing undershoots and overshoots in the saturation profile.
\begin{figure}
	\subfigure[$t=1500$~s  \label{Fig:2D_het_1500}]{
		\includegraphics[clip,scale=0.135,trim=0 0cm 18cm 0]{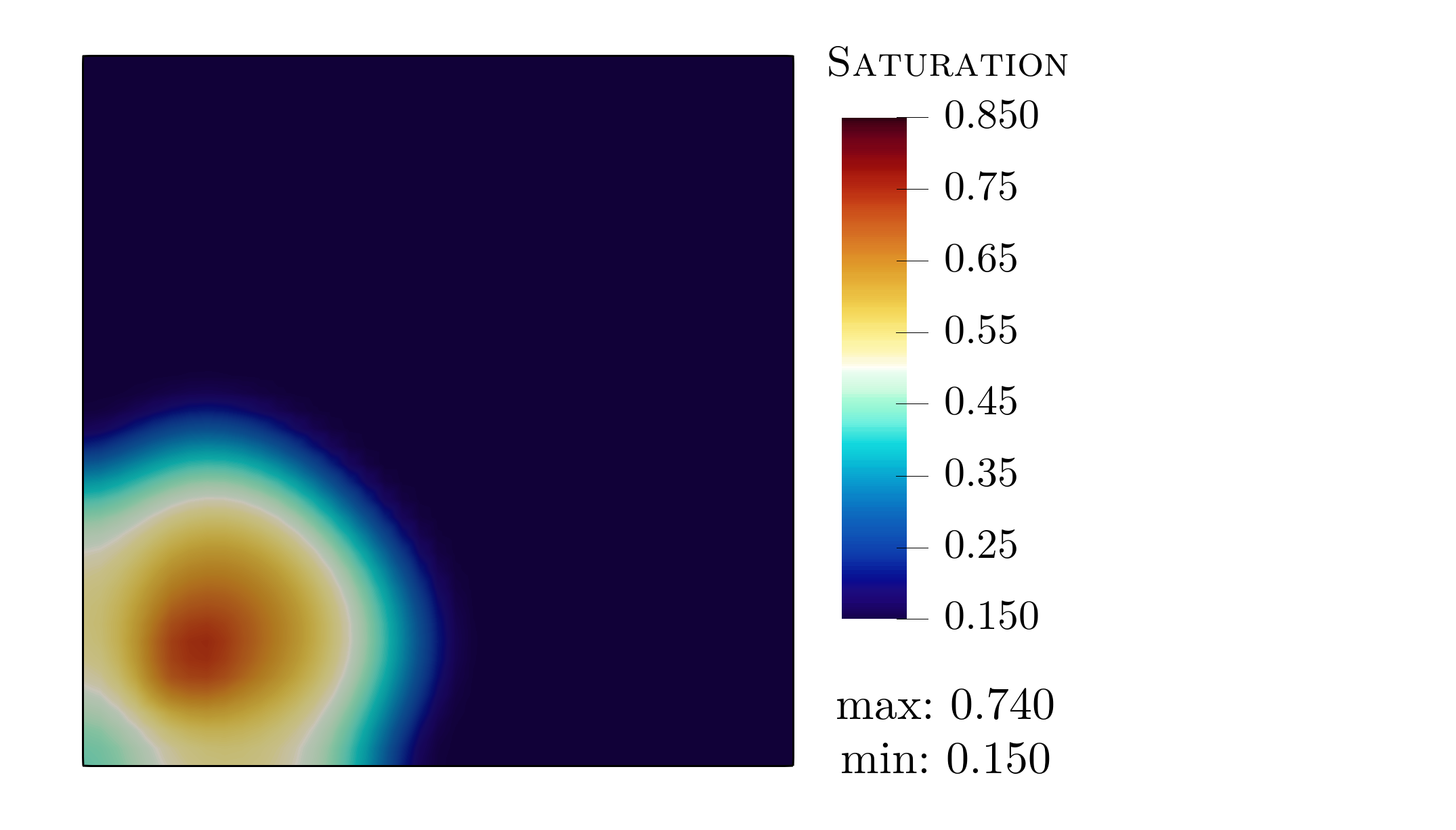}}
	\subfigure[$t=4500$~s  \label{Fig:2D_het_4500}]{
		\includegraphics[clip,scale=0.135,trim=0 0cm 18cm 0]{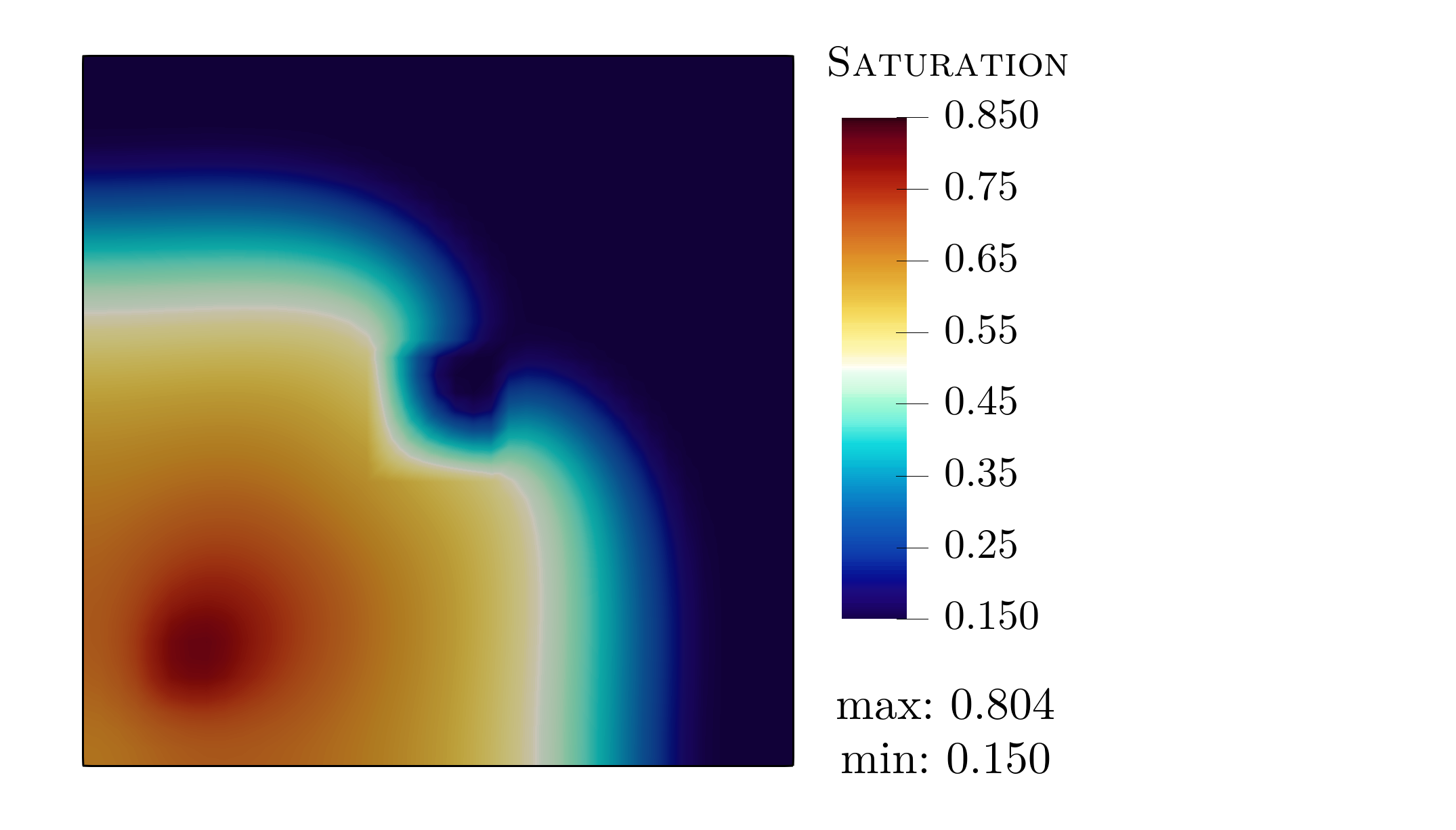}}
        \vspace{0.5cm}
	\subfigure[$t=6900$~s  \label{Fig:2D_het_6900}]{
		\includegraphics[clip,scale=0.135,trim=0 0cm 18cm 0]{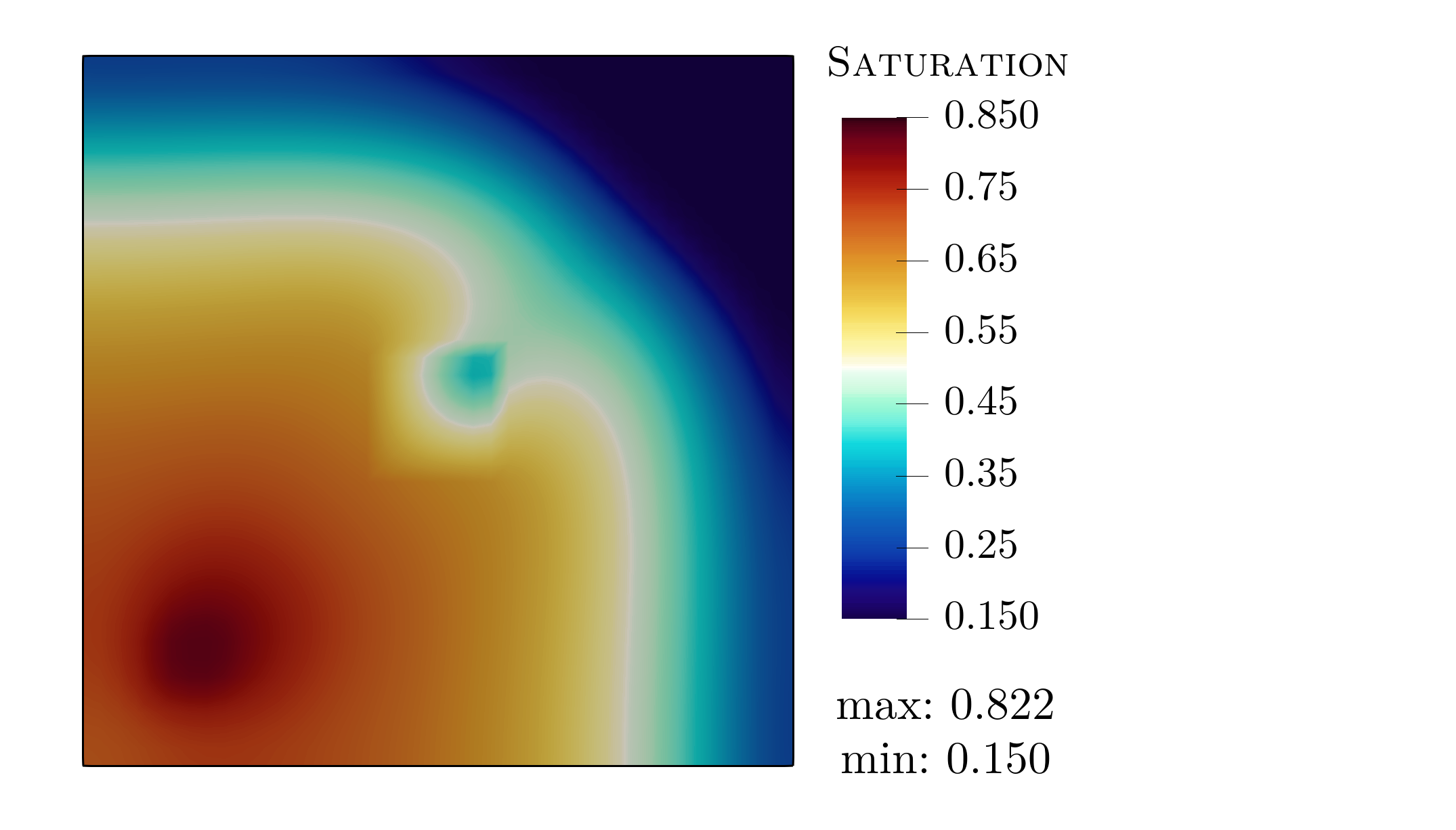}}
	\subfigure[$t=12000$~s  \label{Fig:2D_het_12000}]{
		\includegraphics[clip,scale=0.135,trim=0 0cm 18cm 0]{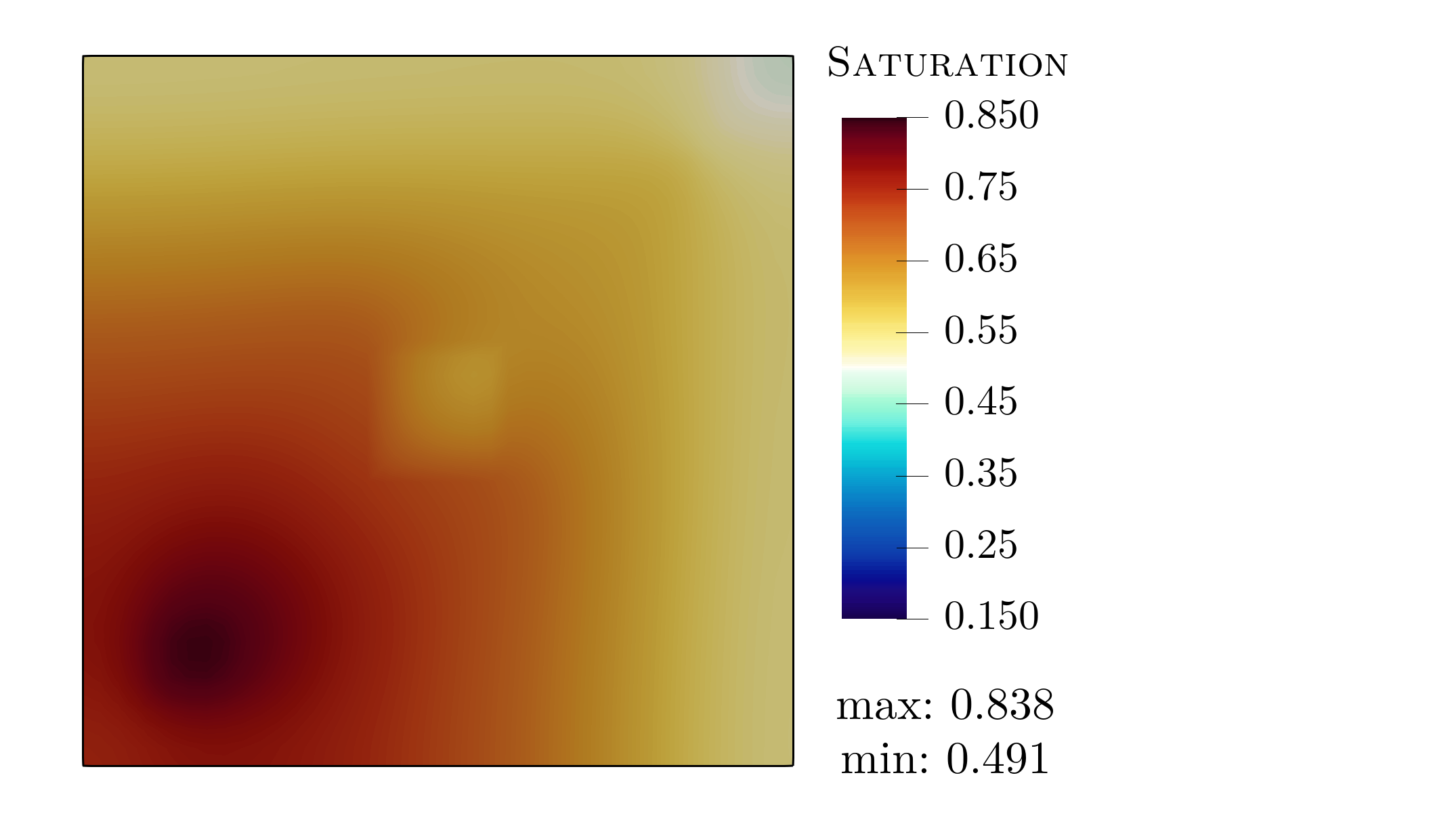}}
        \caption{
            \textsf{Two-dimensional medium with a low permeability block:}
            This figure shows the evolution of saturation solutions. Permeability of block is
            one order of magnitude smaller than the rest of domain. Structured triangular mesh
            (Figure~\ref{Fig:Mesh_unstructured}) was used for this problem.
            Vertex scheme exhibits the expected response, since the saturation avoids the region 
            of low permeability. We note that no violations with respect to maximum principle occurred
            and saturation field remains smooth and monotone even near the corners of the low permeability block.
        \label{Fig:sat_heterogen}}
\end{figure}
Snapshots of the pressure solution along the diagonal line are shown in Figure~\ref{Fig:p_heterogen}.
The less permeable region slightly undulates the curve in that region by  increasing the pressure drop. However,
similar to results obtained in homogeneous porous media (see Figure~\ref{Fig:vertex_vs_DG_pres}), 
the pressure difference drops as more wetting phase reaches the production well.
\begin{figure}
    \centering
    \includegraphics[width=0.8\linewidth]{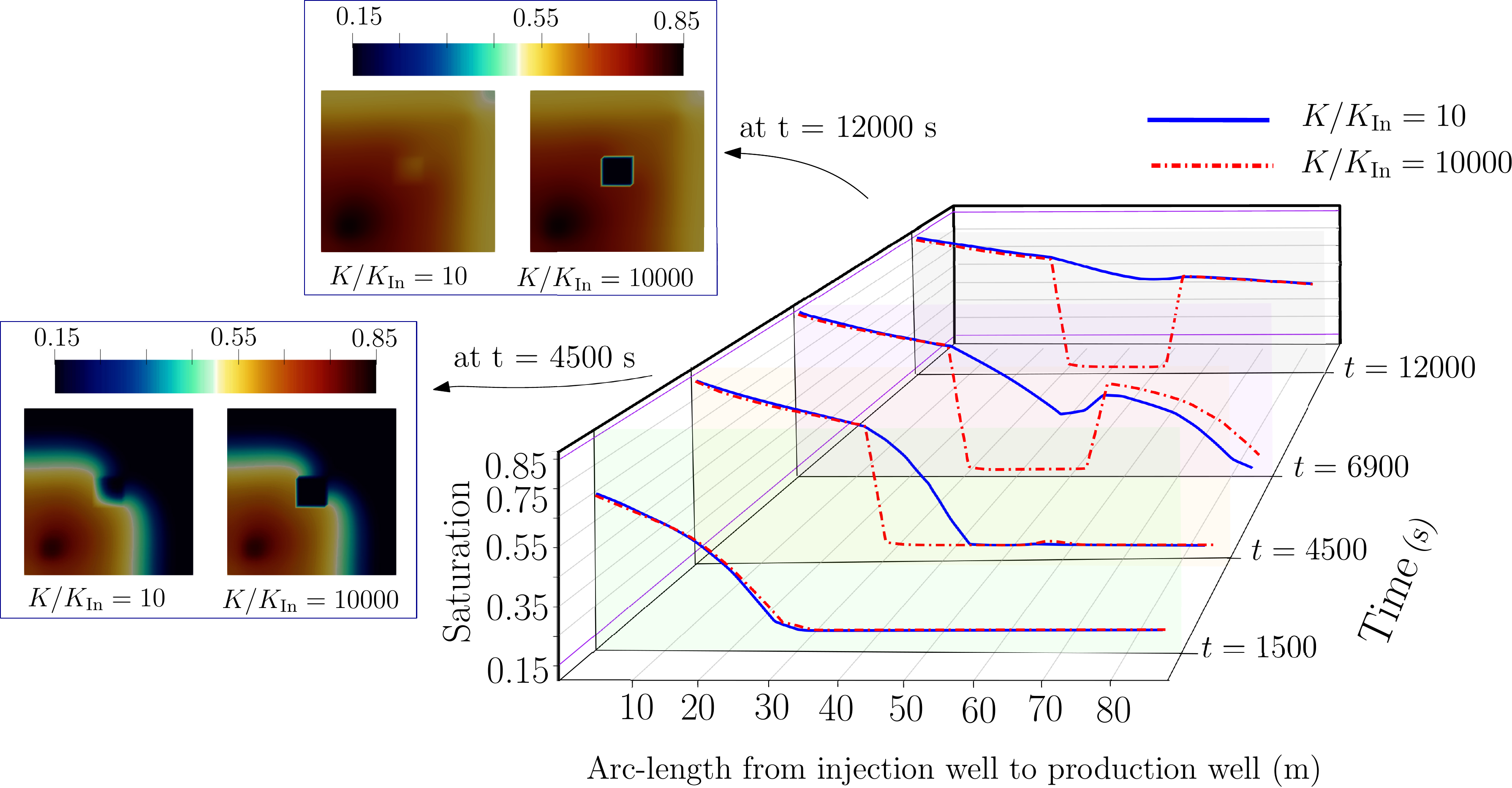}
    \caption{ 
        \textsf{Two-dimensional medium with a low permeability block:}
        This figure shows evolution of the saturation profile along the diagonal line
        and compares the solutions for two cases: (i) domain with a less permeable block 
        (i.e., $K/K_{\mathrm{In}}=10)$ and (ii) domain with almost impermeable block (i.e., $K/K_{\mathrm{In}}=10000$).
        In first case, fluid initially evades the block but as time progresses saturation inside the 
        block starts to increase. However, for the second case, higher permeability difference has made the inclusion 
        impenetrable throughout the simulation.
        We observe that vertex scheme delivers satisfactory results with respect to maximum principle.
    \label{Fig:sat_two_case}}
\end{figure}
\begin{figure}
  \includegraphics[clip,width=0.7\linewidth,trim=0 .4cm 0.5cm 0]{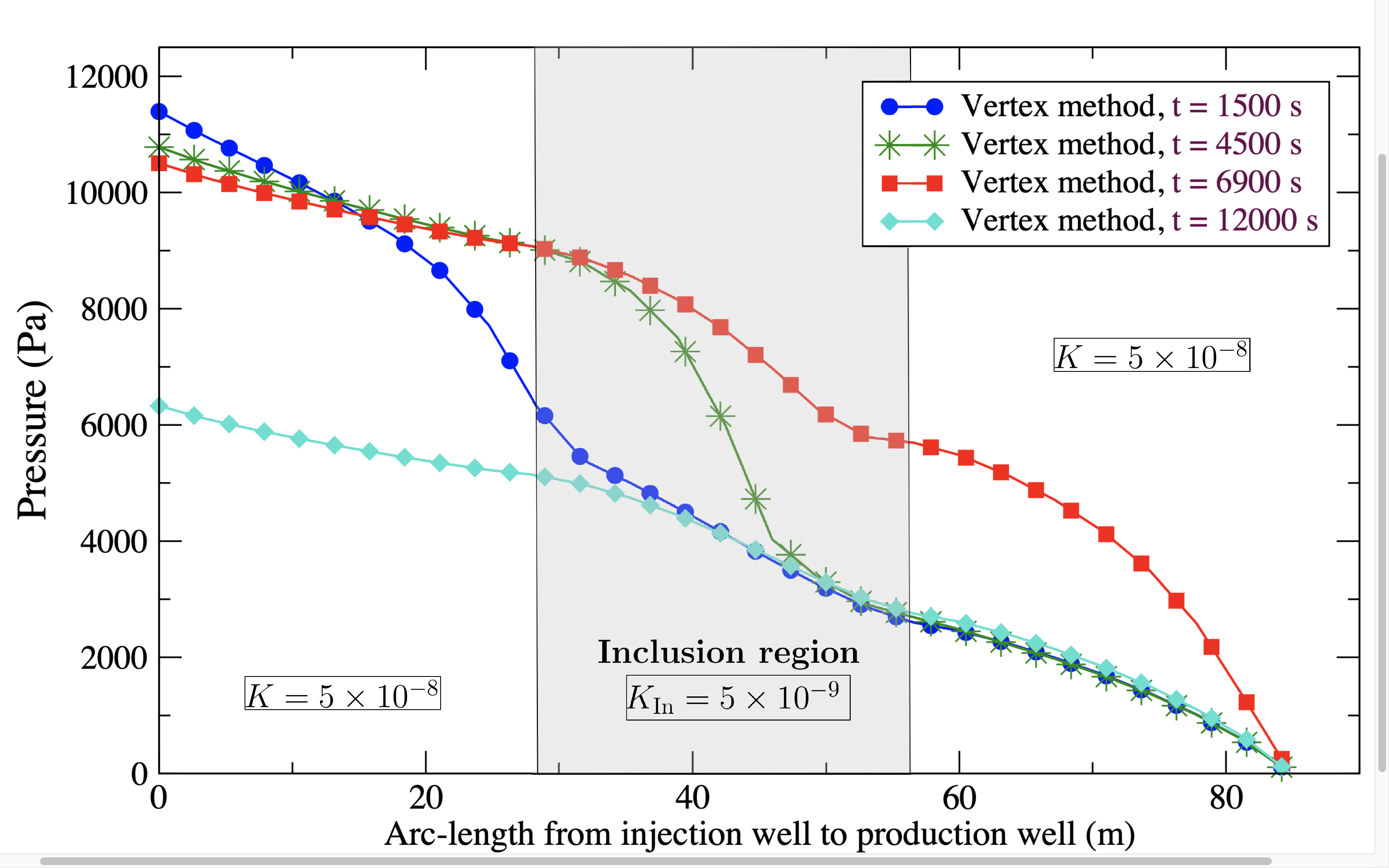}
  \caption{
      \textsf{Two-dimensional medium with a low permeability block:}
      This figure shows the evolution of the numerical pressure.
      Profiles are plotted along the diagonal line from point $(20,20)$~m to $(80,80)$~m. 
      The pressure values drop as we move from the injection well toward the production well.
      The inclusion (or the block) region, which is illustrated with gray-color, slightly perturbs 
      the solutions. 
  \label{Fig:p_heterogen}}
\end{figure}
\subsubsection{Two-dimensional porous medium with highly heterogeneous permeability}
In this example, the domain $\Omega=[0,1000]^2\;\mathrm{m^2}$ is highly heterogeneous because the permeability field
is taken from various horizontal permeability slices from model 2 of the SPE10 benchmark model \citep{christie2001tenth,SPE10}.
This model is characterized by two formations: a shallow-marine Tarbert formation in the top 35 layers, 
where the permeability field is relatively smooth, and a fluivial Upper-Ness permeability in the bottom 50 layers. 
Both formations are characterized by large permeability variations, 8--12 orders of magnitude, but are qualitatively different. 
We choose layer 1 from Tarbert formation and layer 45 and 80 from Upper-Ness formation.
Figure~\ref{Fig:K_layer1},~\ref{Fig:K_layer45}, and~\ref{Fig:K_layer80} show the selected permeability layers.
\begin{figure}
	\subfigure[Layer $1$  \label{Fig:K_layer1}]{
		\includegraphics[clip,scale=0.33,trim=0 0cm 0cm 0]{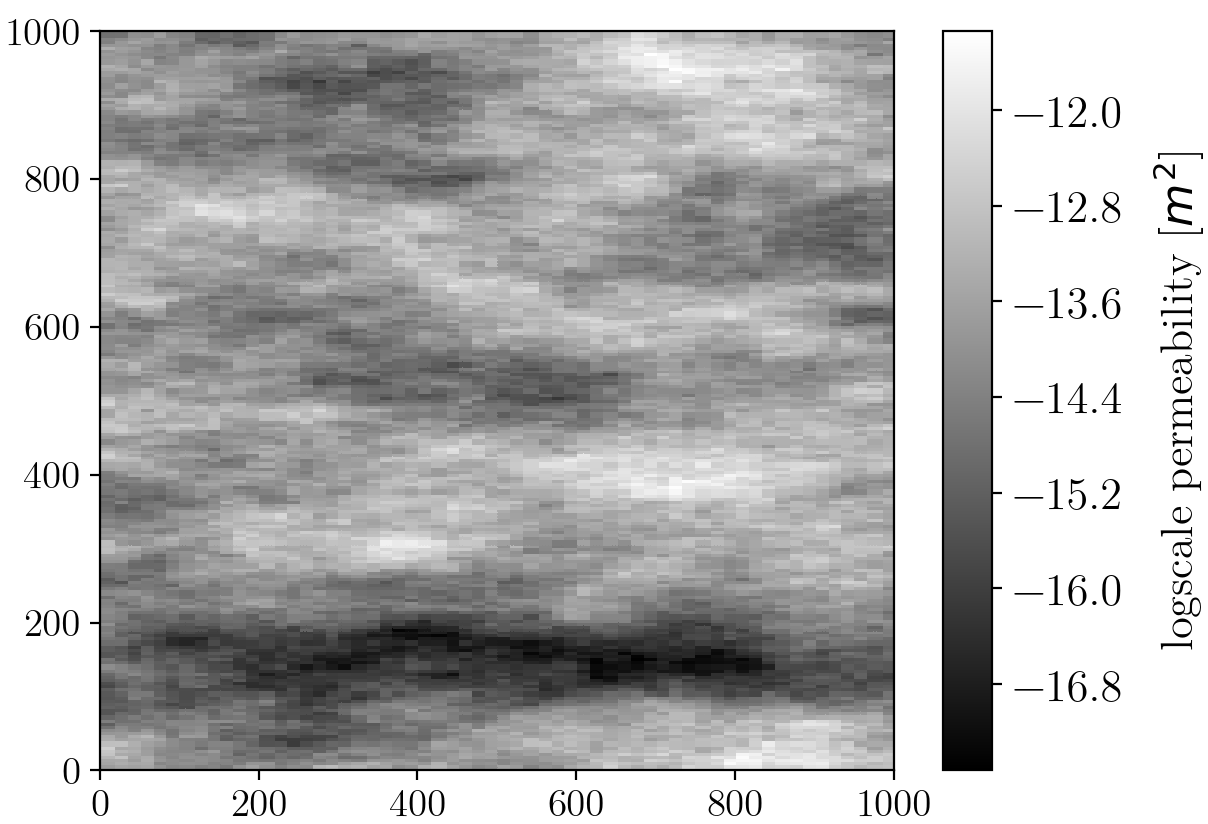}}
	\subfigure[Layer $45$ \label{Fig:K_layer45}]{
		\includegraphics[clip,scale=0.33,trim=0 0cm 0  0]{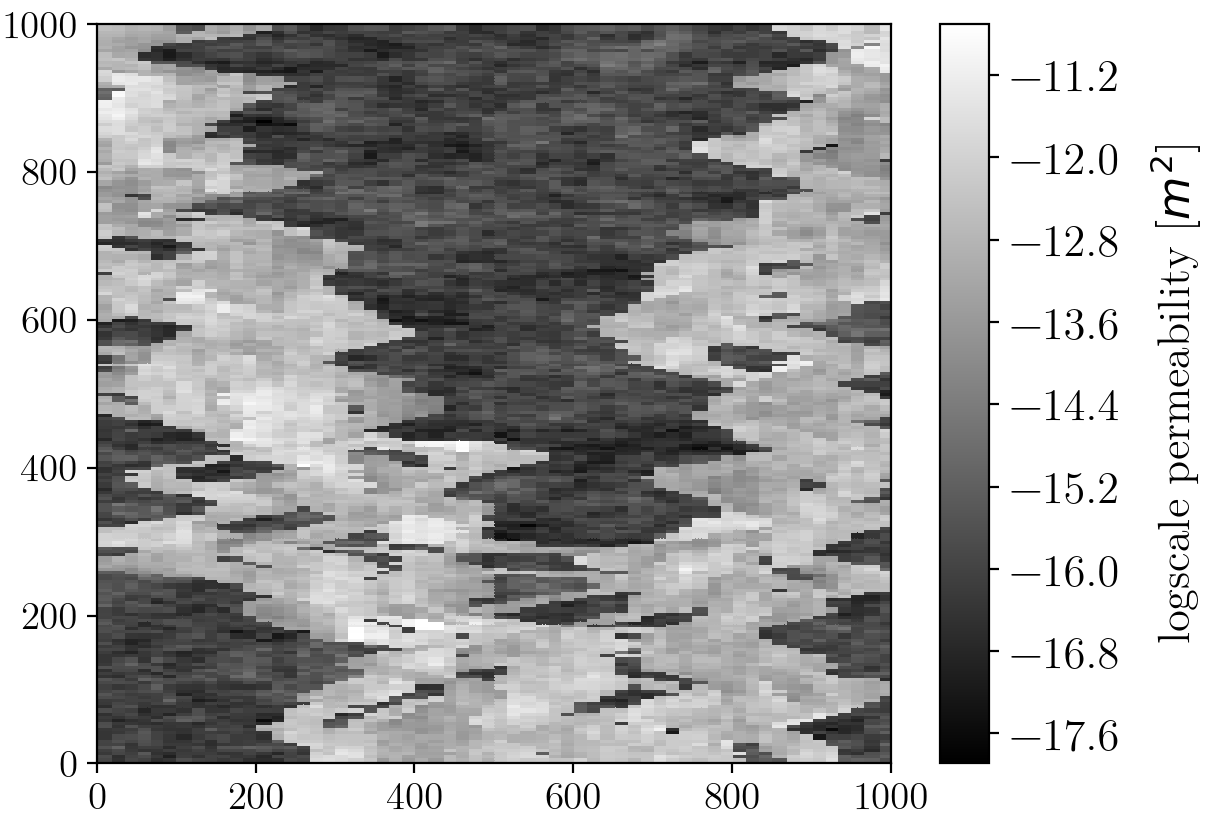}}
	\subfigure[Layer $80$ \label{Fig:K_layer80}]{
		\includegraphics[clip,scale=0.33,trim=0 0cm 0 0]{Figures/Fig10_c.png}}
        \caption{
            \textsf{Two-dimensional porous medium with highly heterogeneous permeability:}
            This figure depicts the permeability fields adopted from three horizontal layers of SPE10 benchmark model.
            Each field is scaled to a resolution of $60 \times 60$ grids.  
            Layer $1$ taken from relatively smooth Tarbert formation, while layer $45$ and $80$ are taken from more rugged 
            Upper-Ness formation.
            Values are displayed in logarithmic scale, since they  vary across a wide range.
        \label{Fig:2D_SPE10_K}}
\end{figure}
These permeability slices are scaled to a $60 \times 60$ grid, instead of the original 60 × 220 grid.
The porosity is set to $0.2$. No flow boundary conditions are prescribed on the entire boundary and as shown 
in Figure~\ref{Fig:Schematic_SPE10} an injection well of size $100\times100$~m is defined on
the center of domain and four production wells with size of $100\times100$~m are located near the corners of the domain.
\begin{figure}
	\subfigure[Two-dimensional domain\label{Fig:2D_SPE10_schematic}]{
		\includegraphics[clip,scale=0.25,trim=0 0cm 0cm 0]{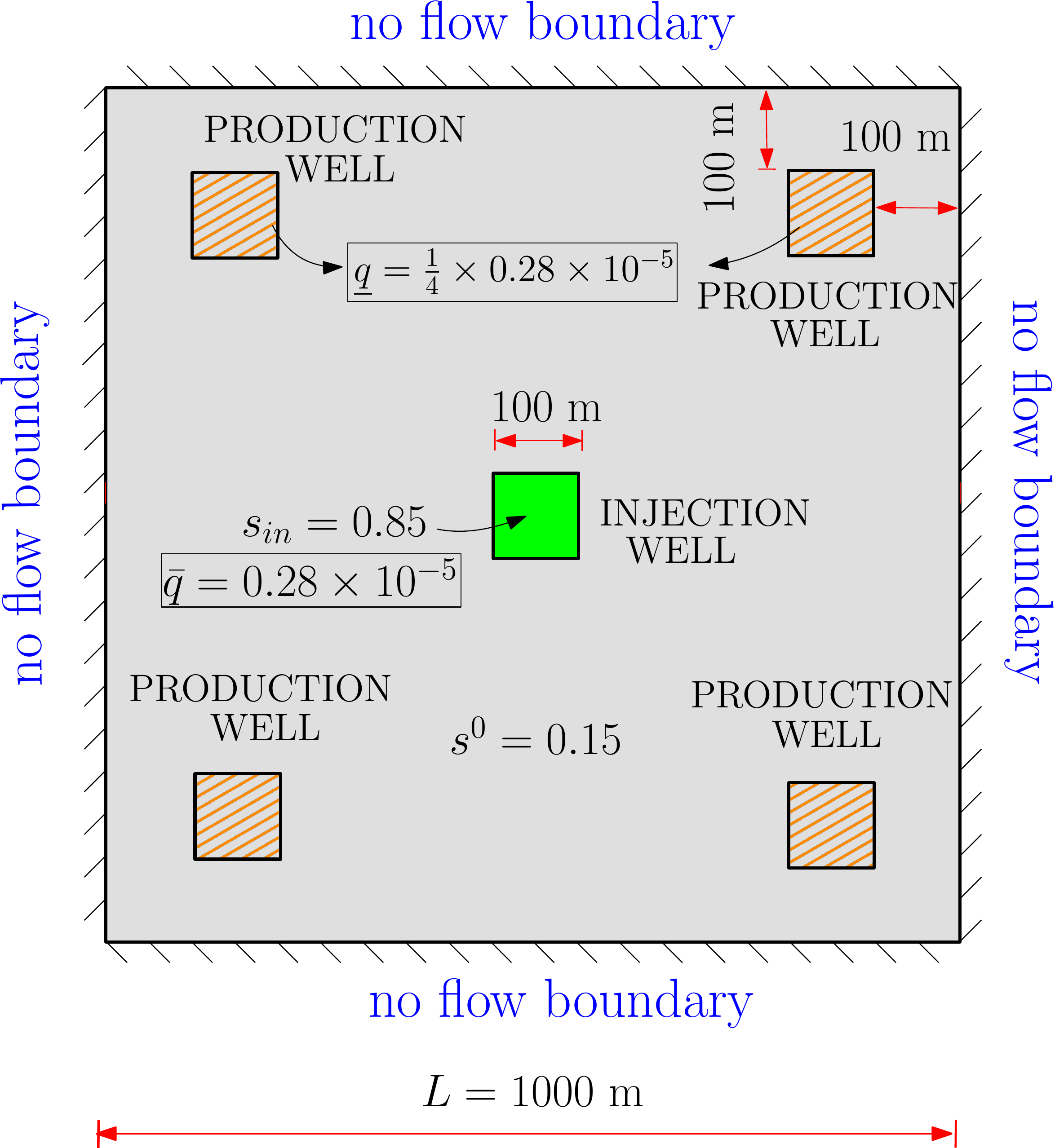}}
        \hspace{0.8cm}
	\subfigure[Three-dimensional domain\label{Fig:3D_SPE10_schemetic}]{
		\includegraphics[clip,scale=0.25,trim=0 0cm 0cm 0]{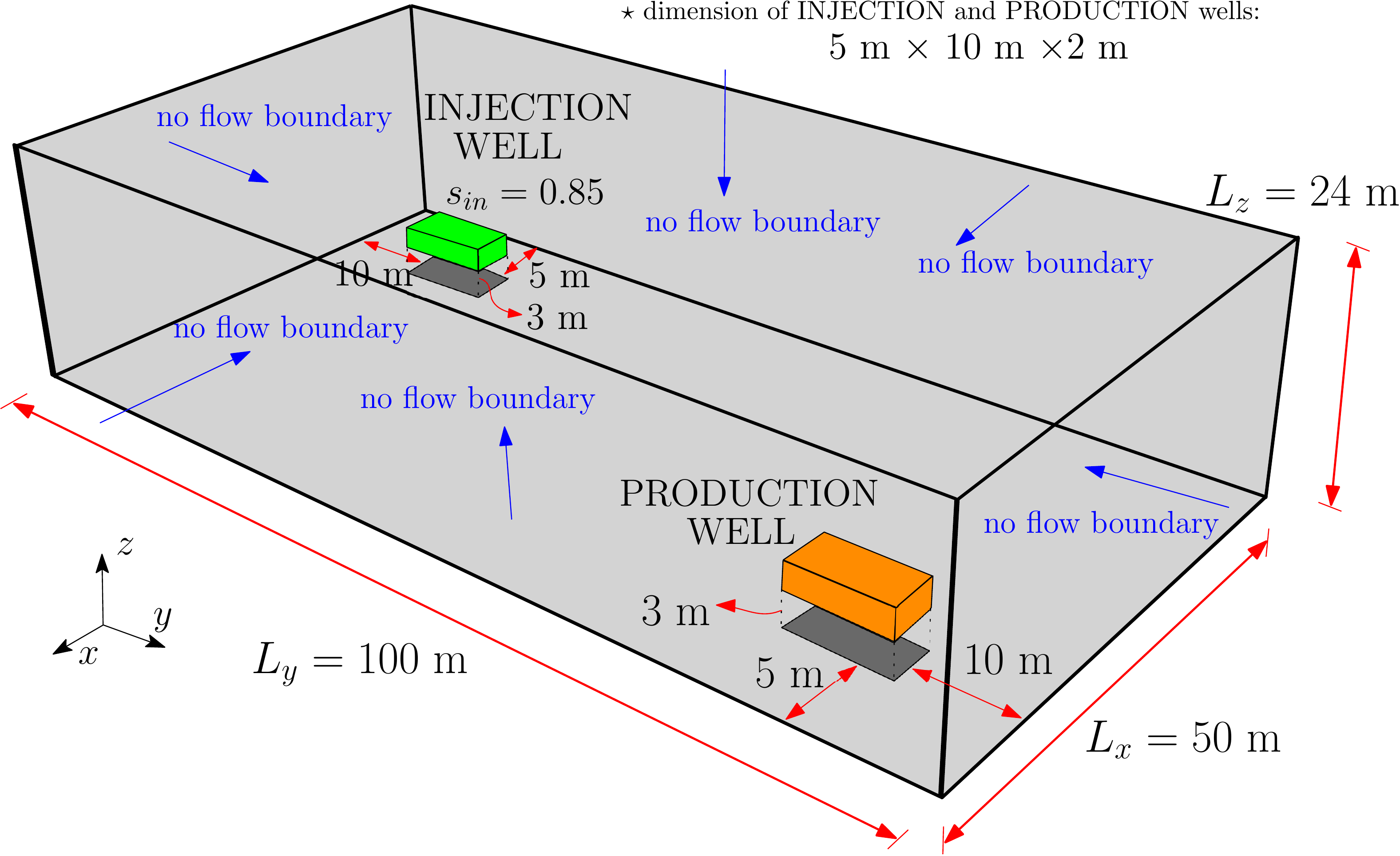}}
        \caption{\textsf{Highly heterogeneous problems:~}
            This figure provides a pictorial description of the boundary value problem and shows
            the computational domain used in numerical simulations. For both two-dimensional and three-dimensional
            problems no flow conditions are assigned across all boundaries. Fluid flow is hence driven by the pressure difference
            prompted by injection and production wells.
        \label{Fig:Schematic_SPE10}}
\end{figure}
Injection and production flow rates are piecewise constant with compact support and are determined by the following constraint:
\begin{align}
    \int_{\Omega}\bar{q} =\int_{\Omega}\underline{q} = 2.8 \times 10^{-1}.
\end{align}
Here we employ the proposed finite element scheme on a structured mesh with $7200$ triangular elements. The simulation runs
to $T=2.5$~days with $600$ time steps, and we provide the solutions at $t=0.4125$, $t=0.725$, and $2.5$~days. 
Saturation contours are depicted in Figure~\ref{Fig:2D_SPE10_sat}.
The wetting phase moves from the injection wells towards the four production wells as expected.
The permeability field determines the pattern of the saturation front throughout the porous media.  
For all three cases, physical instabilities in form of separate finger-like intrusions are generated. 
As expected, the saturation front forms a curve that is less smooth for  porous media 
of Upper Ness types (i.e.,~layers $45$ and $80$).
It is evident that the vertex scheme produces bound-preserving saturations 
and that fronts avoid small regions of lower permeability.
\begin{figure}
	\subfigure[Layer $1$; $t=0.4125$ days \label{Fig:2D_layer1_100}]{
		\includegraphics[clip,scale=0.10,trim=0 0cm 0cm 0]{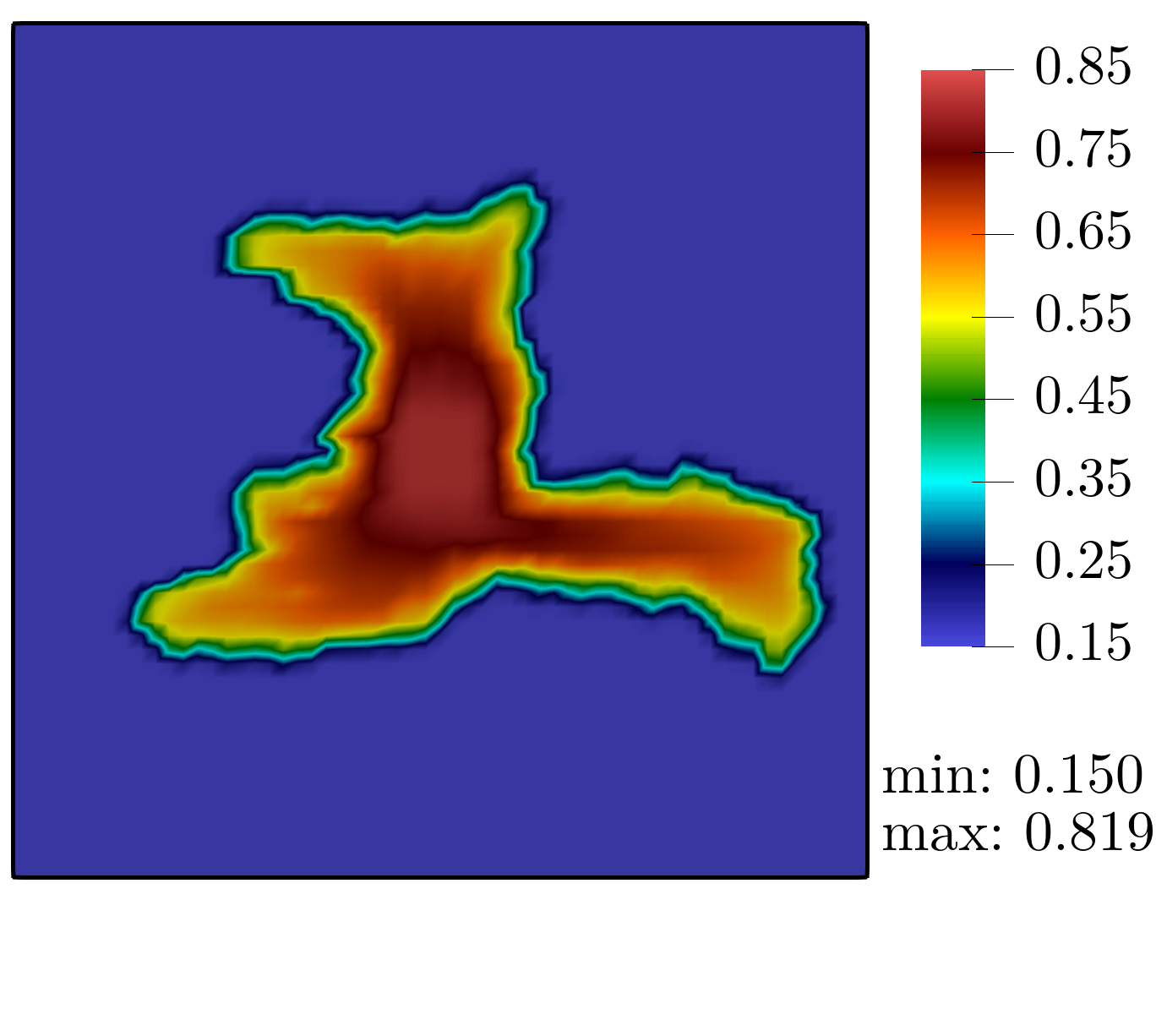}}
	\subfigure[Layer $1$; $t=0.725$ days  \label{Fig:2D_layer1_175}]{
		\includegraphics[clip,scale=0.10,trim=0 0cm 0cm 0]{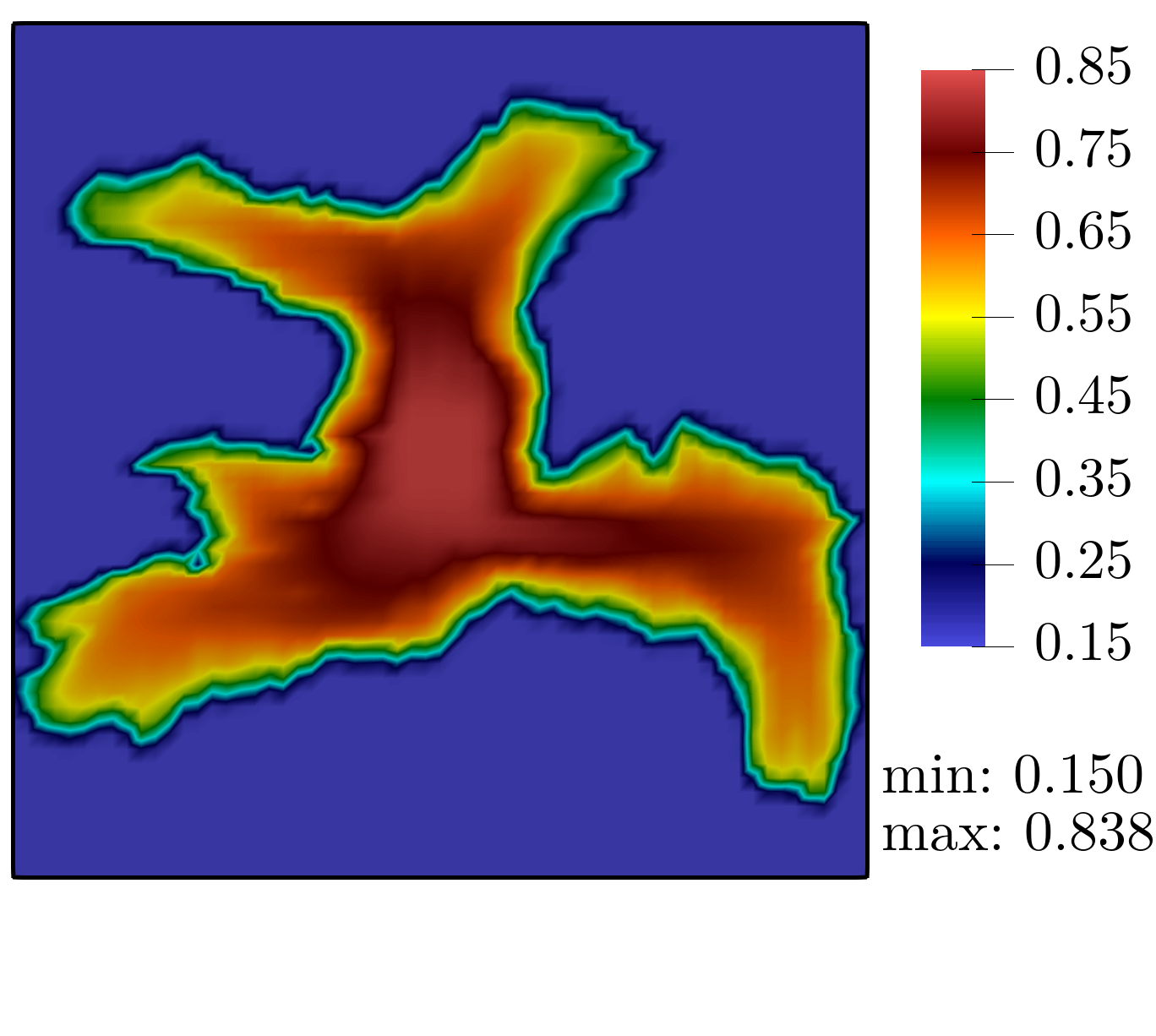}}
	\subfigure[Layer $1$; $t=2.5$ days \label{Fig:2D_layer1_600}]{
		\includegraphics[clip,scale=0.10,trim=0 0cm 0cm 0]{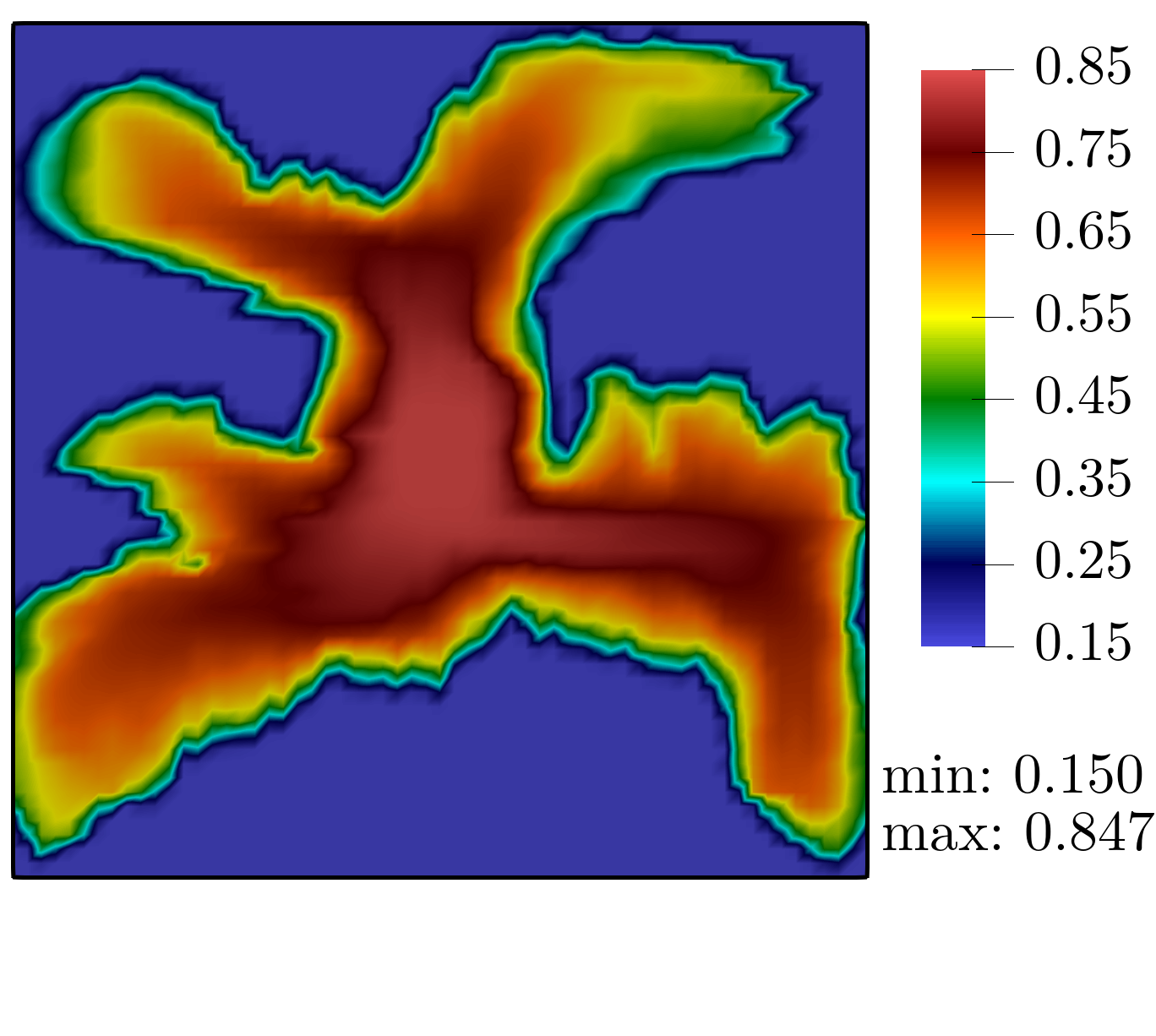}}
        \vspace{0.2cm}
	\subfigure[Layer $45$; $t=0.4125$ days \label{Fig:2D_layer45_100}]{
		\includegraphics[clip,scale=0.10,trim=0 0cm 0cm 0]{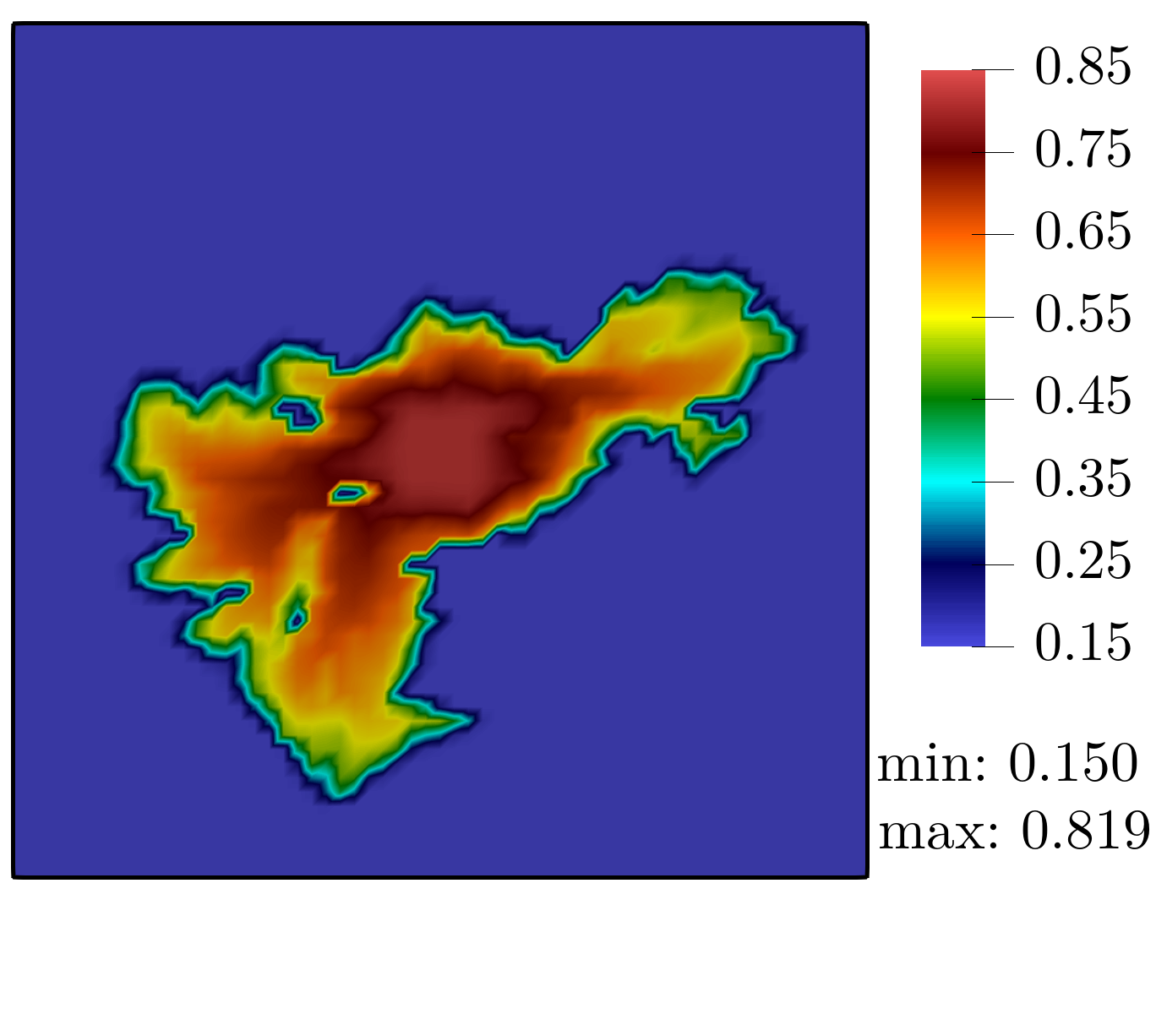}}
	\subfigure[Layer $45$; $t=0.725$ days \label{Fig:2D_layer45_175}]{
		\includegraphics[clip,scale=0.10,trim=0 0cm 0cm 0]{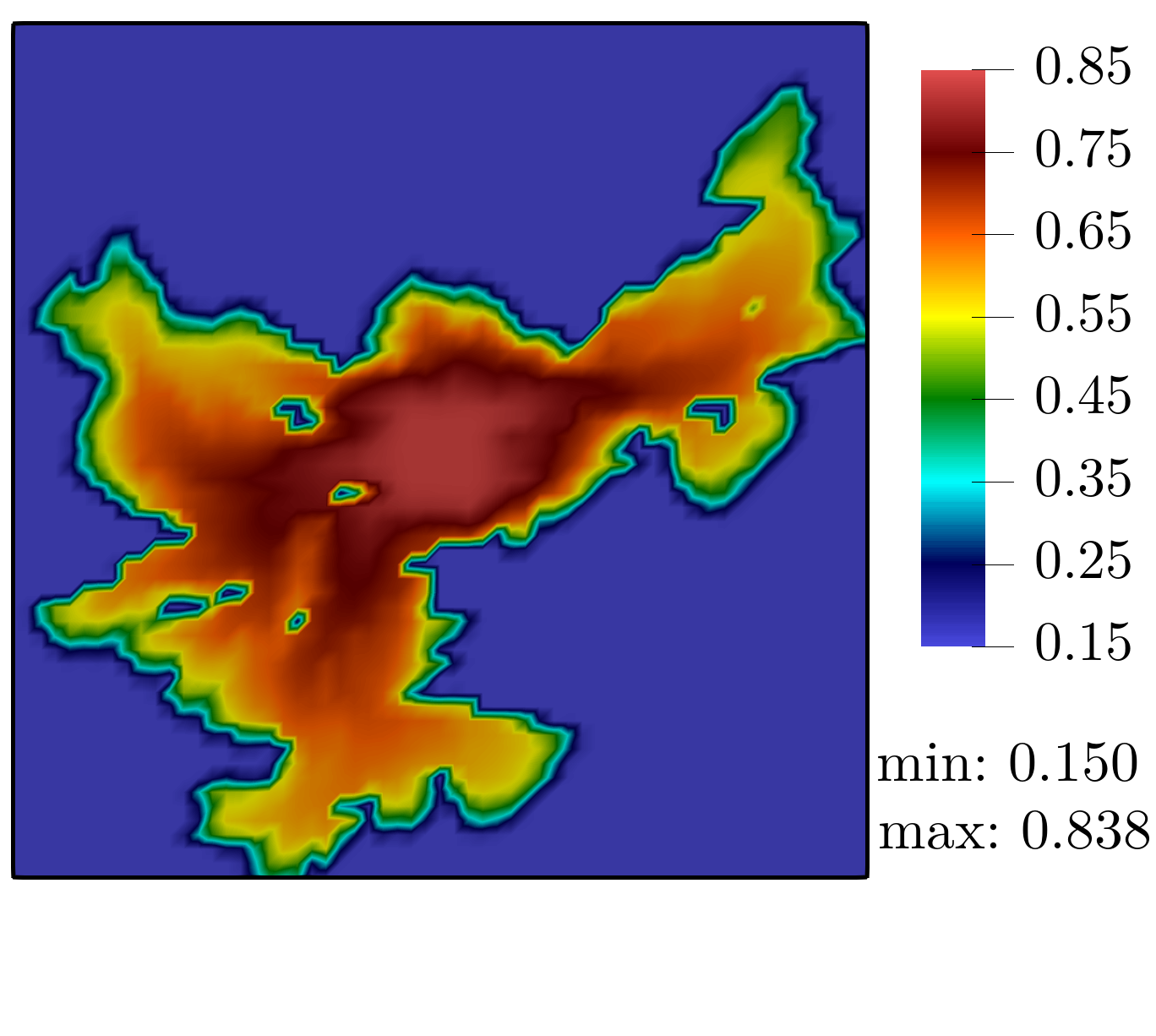}}
	\subfigure[Layer $45$; $t=2.5$ days \label{Fig:2D_layer45_600}]{
		\includegraphics[clip,scale=0.10,trim=0 0cm 0cm 0]{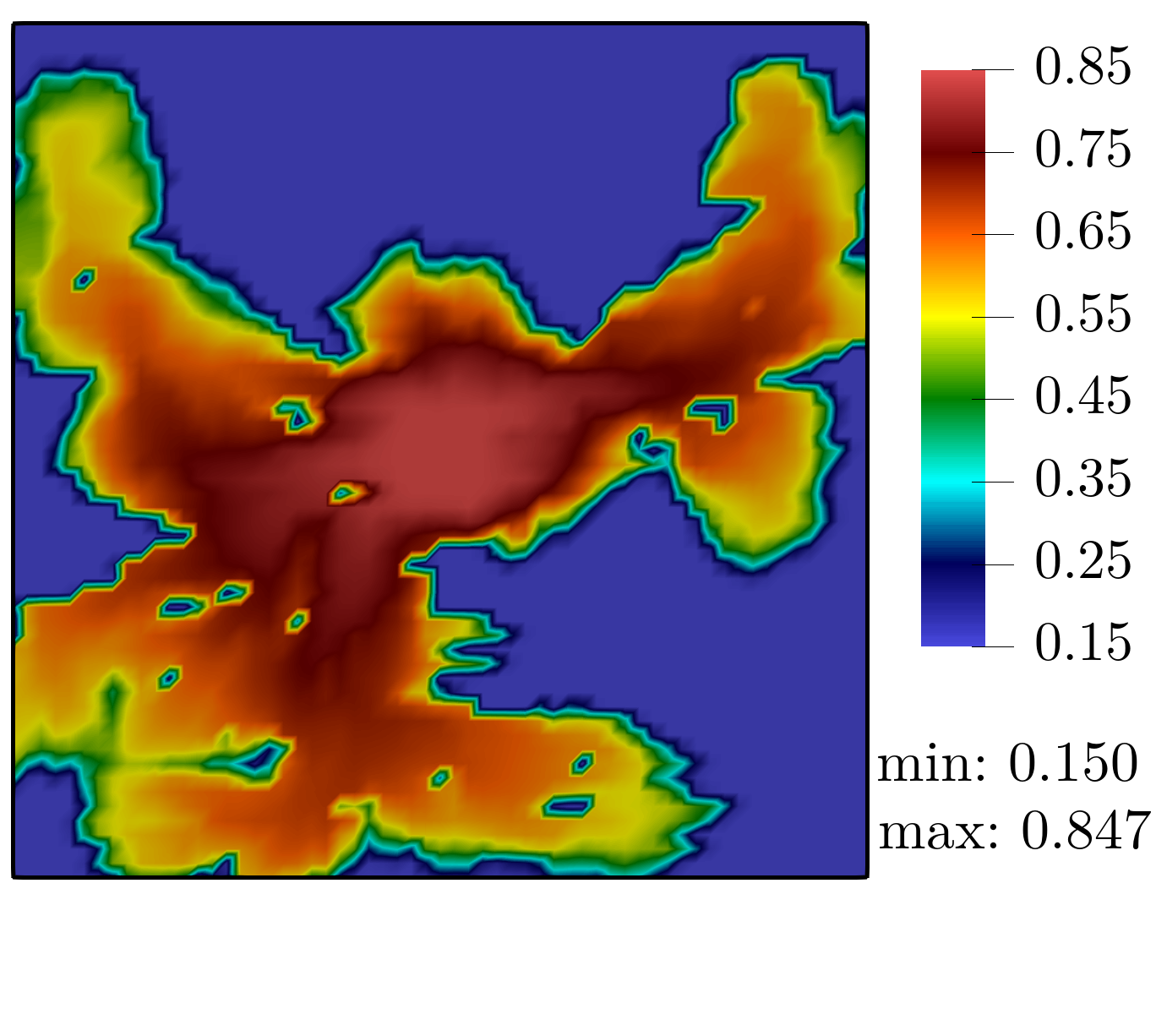}}
        \vspace{0.2cm}
	\subfigure[Layer $80$; $t=0.4125$ days \label{Fig:2D_layer80_100}]{
		\includegraphics[clip,scale=0.10,trim=0 0cm 0cm 0]{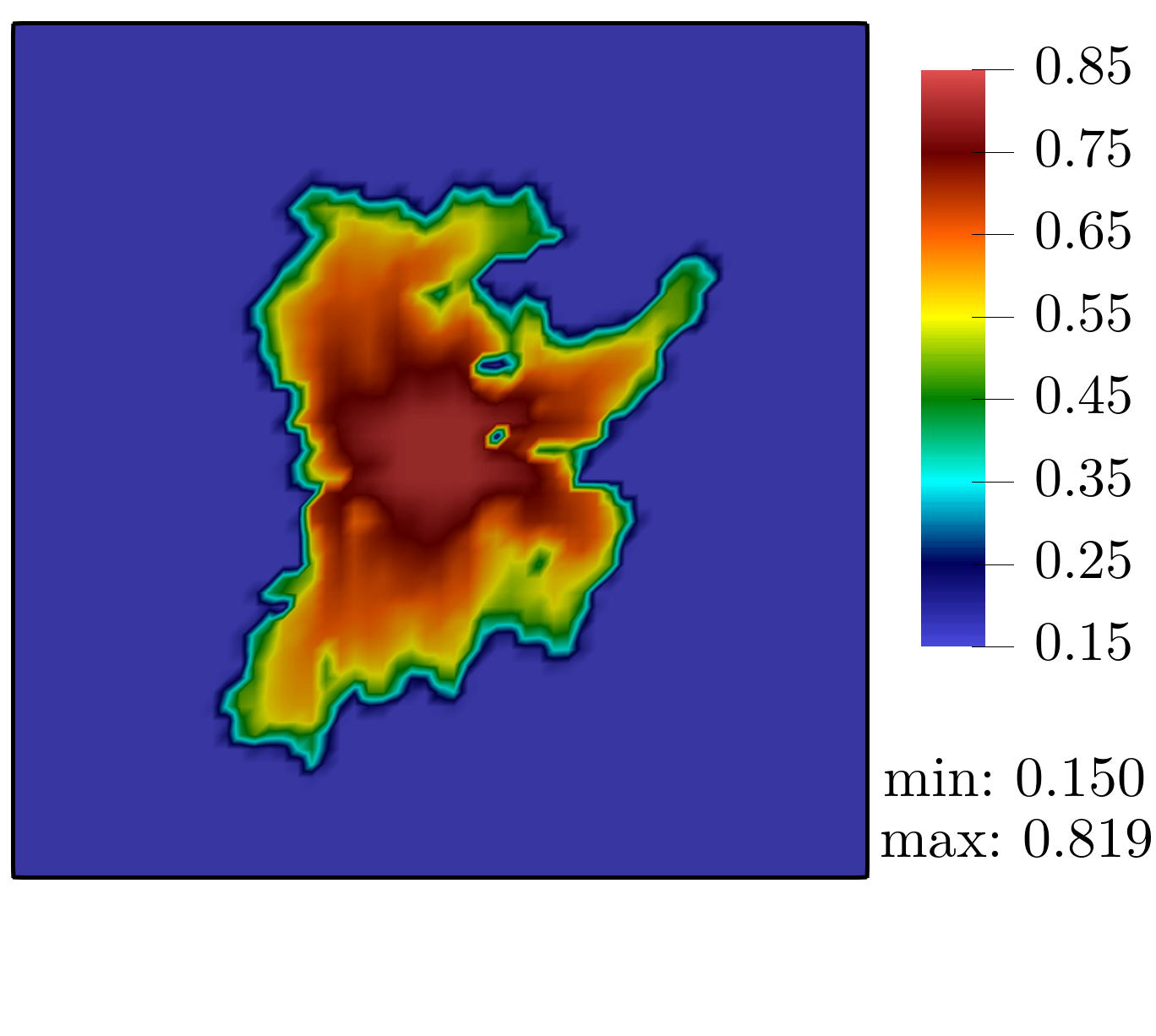}}
	\subfigure[Layer $80$; $t=0.725$ days \label{Fig:2D_layer80_175}]{
		\includegraphics[clip,scale=0.10,trim=0 0cm 0cm 0]{Figures/Fig12_h.png}}
	\subfigure[Layer $80$; $t=2.5$ days \label{Fig:2D_layer80_600}]{
		\includegraphics[clip,scale=0.10,trim=0 0cm 0cm 0]{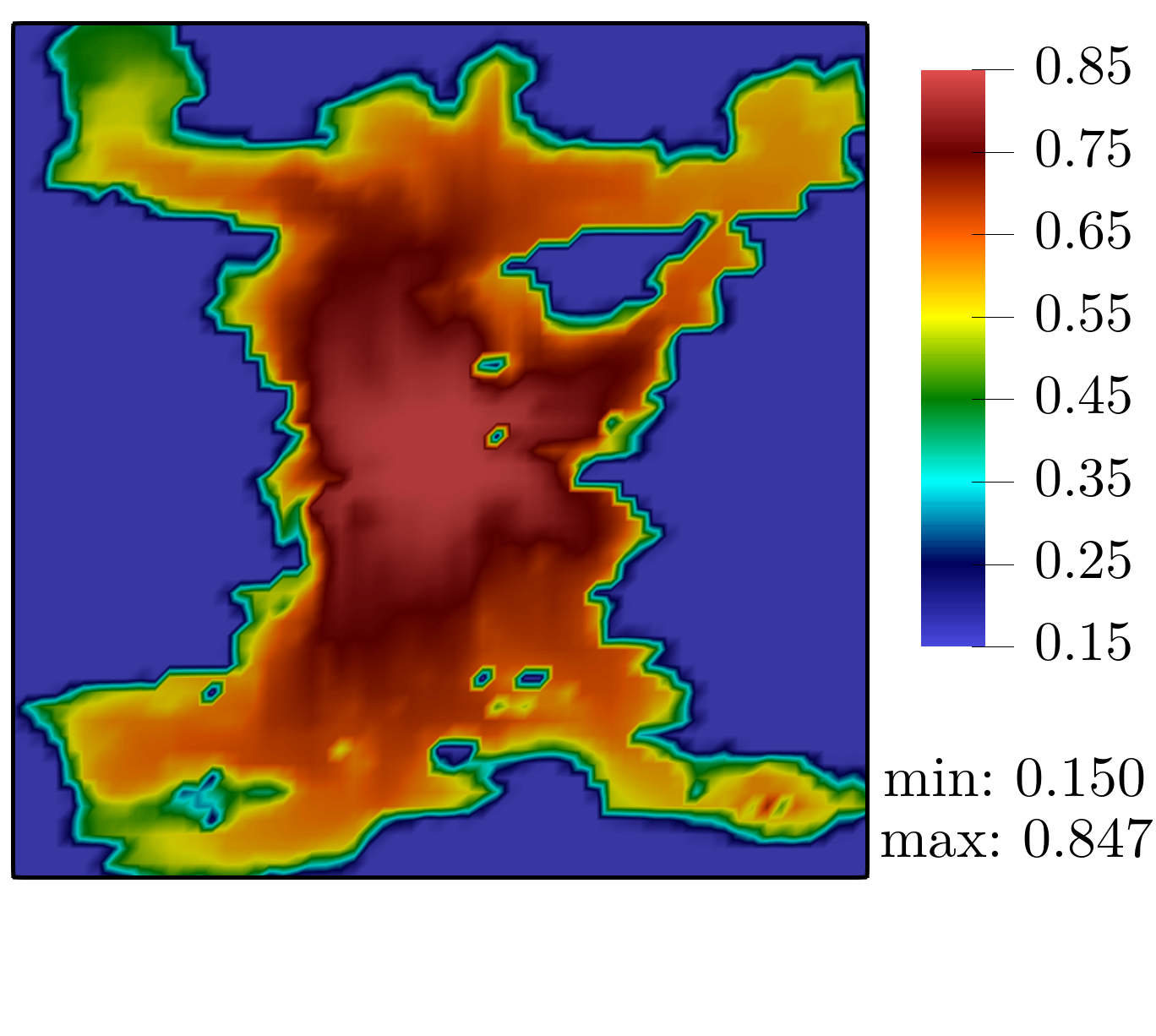}}
        \caption{
            \textsf{Two-dimensional heterogeneous problem:~}
            This figure shows the evolution of the wetting phase saturation for the three chosen layers. 
            For all cases, permeability field determines the pattern fluid flows through porous media.
            Layer $1$ (top figures), compared to other layers, leads to smoother saturation front boundaries. 
            This observation is justified as the permeability field associated with layer $1$ is not as highly varying as for the other layers.  This figure also reiterates that the vertex scheme always produces physical values of saturations, without any overshoots and undershoots,
            even for domains with permeabilites that vary over  many orders of magnitudes.
        \label{Fig:2D_SPE10_sat}}
\end{figure}

\subsubsection{Three-dimensional porous medium}
\label{Sub:SPE10-2D}
Herein we validate our vertex scheme in a three-dimensional set-up.
In particular, we investigate an extension of the numerical experiment performed in Section \ref{sub:2D_homogen}.
The domain is $\Omega=[0,1000]^3~\mathrm{m}^3$ and is partitioned into an unstructured mesh
of tetrahedron elements, as shown in Figure~\ref{Fig:3D_mesh}. The permeability is fixed to $K=5 \times 10^{-8}\;\mathrm{m^2}$.
No flow boundary condition is employed on the entire boundary $\partial \Omega$.
Production and injection wells of size $(10,10,10)$~m with constant flow rates of $\bar{q}=\underline{q}=0.001$ 
are positioned at the opposite corners of the domain.
Figure~\ref{Fig:3D_BVP} shows the computational domain and the boundary conditions for this problem.
The final time is set to $T=2$ days and the time step is set to $\tau=216$ s.
\begin{figure}
	\subfigure[Computational domain \label{Fig:3D_BVP}]{
		\includegraphics[clip,scale=0.26,trim=0 0cm 0 0]{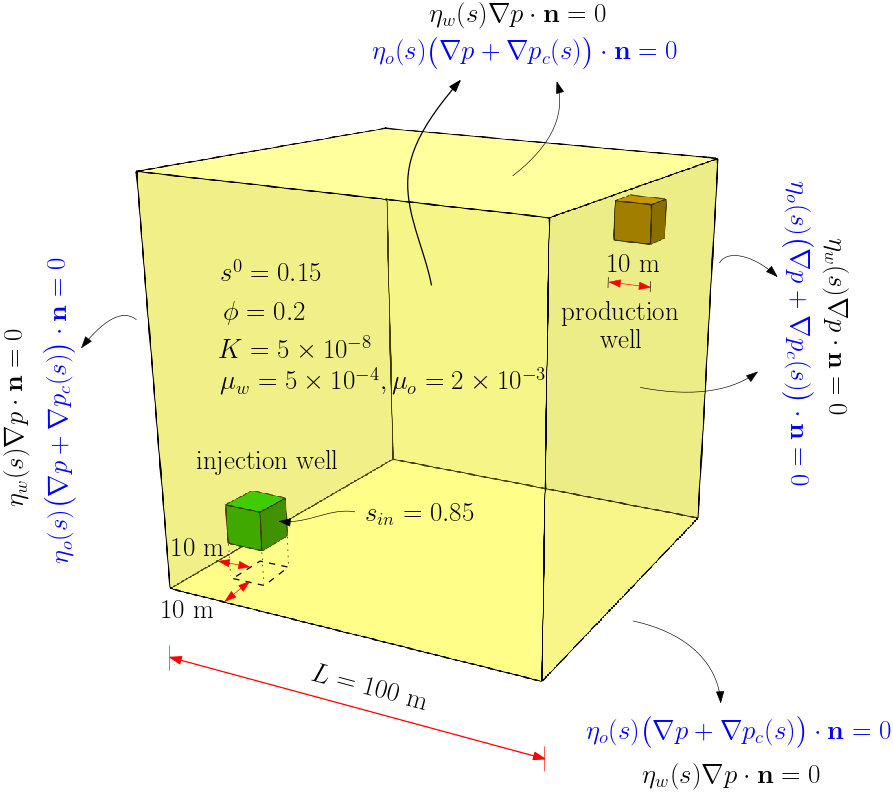}}
            \hfill
        \subfigure[Mesh~\label{Fig:3D_mesh}]{
		\includegraphics[clip,scale=0.33,trim=0 0cm 0 0]{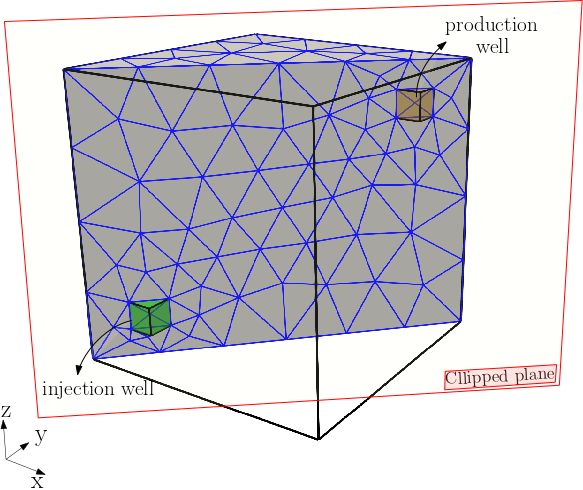}}
        \caption{\textsf{Three-dimensional porous medium with homogeneous permeability}
            Left figure provides a pictorial description of the boundary value problem. no flow conditions prescribed on all boundaries.
            Right figure shows a cross-sectional view of the mesh applied in our numerical experiment.
  \label{Fig:3D_domain_schematic}}
\end{figure}
In Figure~\ref{Fig:3D_sat}, snapshots of the wetting phase saturation are given at times $t=0.125$, $t=0.5$,
$t=1$, and $t=2$ days.
Profile of pressure at the corresponding time steps, 
along the diagonal (from injection well to production well), are exhibited in Figure~\ref{Fig:3D_pressure}.
One can observe that the numerical scheme is robust in three-dimensional domain and the resulting saturation
satisfies the maximum principle.
Only $5$ to $6$ Picard's iterations are needed at each time step for convergence of the vertex scheme. 
This is true for all three-dimensional test cases.
\begin{figure}
	\subfigure[$t=0.125$ days \label{Fig:BoM_vertex_1800}]{
		\includegraphics[clip,scale=0.135,trim=0 0cm 18cm 0]{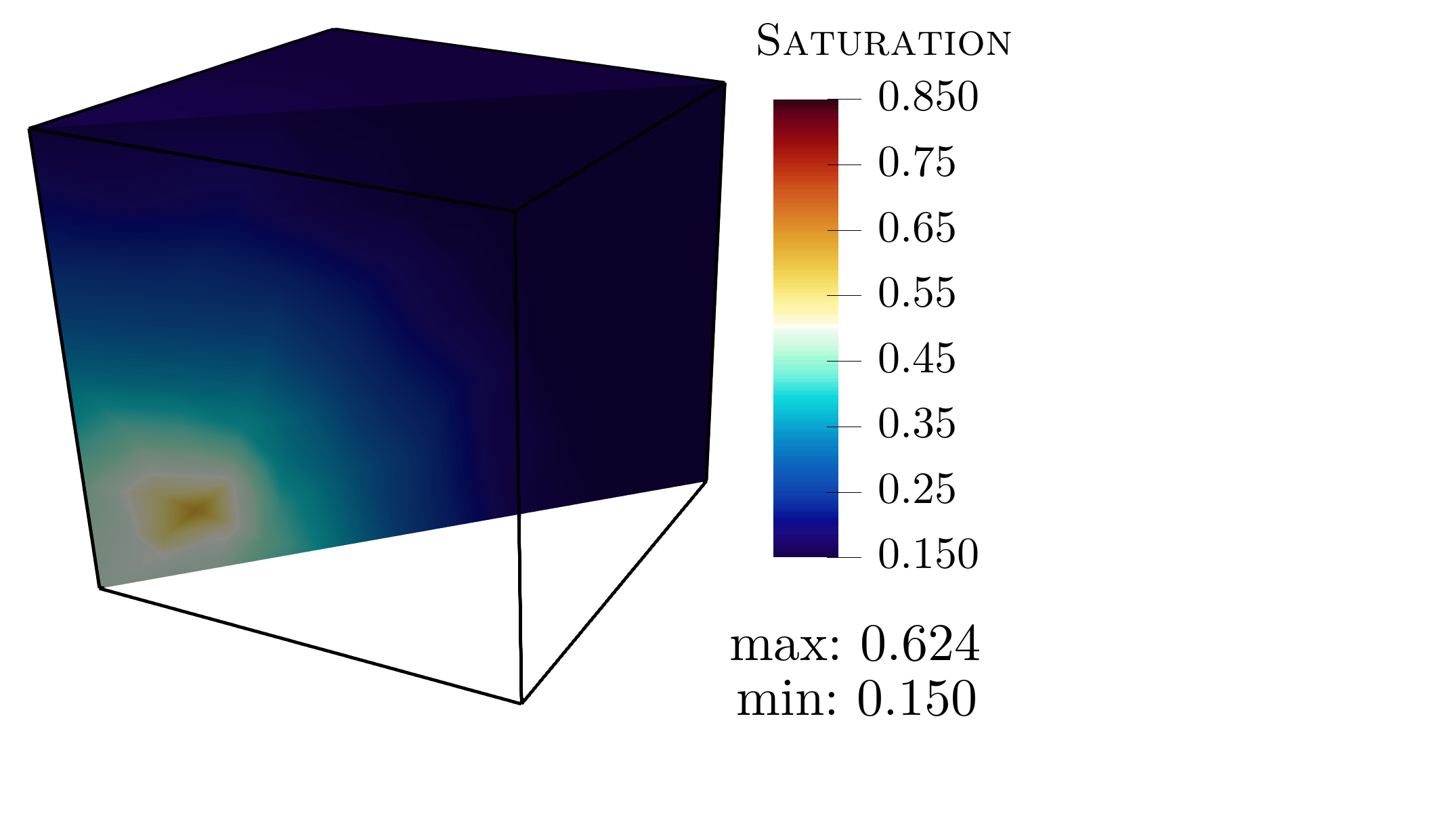}}
	\subfigure[$t=0.5$ days \label{Fig:BoM_vertex_1800}]{
		\includegraphics[clip,scale=0.135,trim=0 0cm 18cm 0]{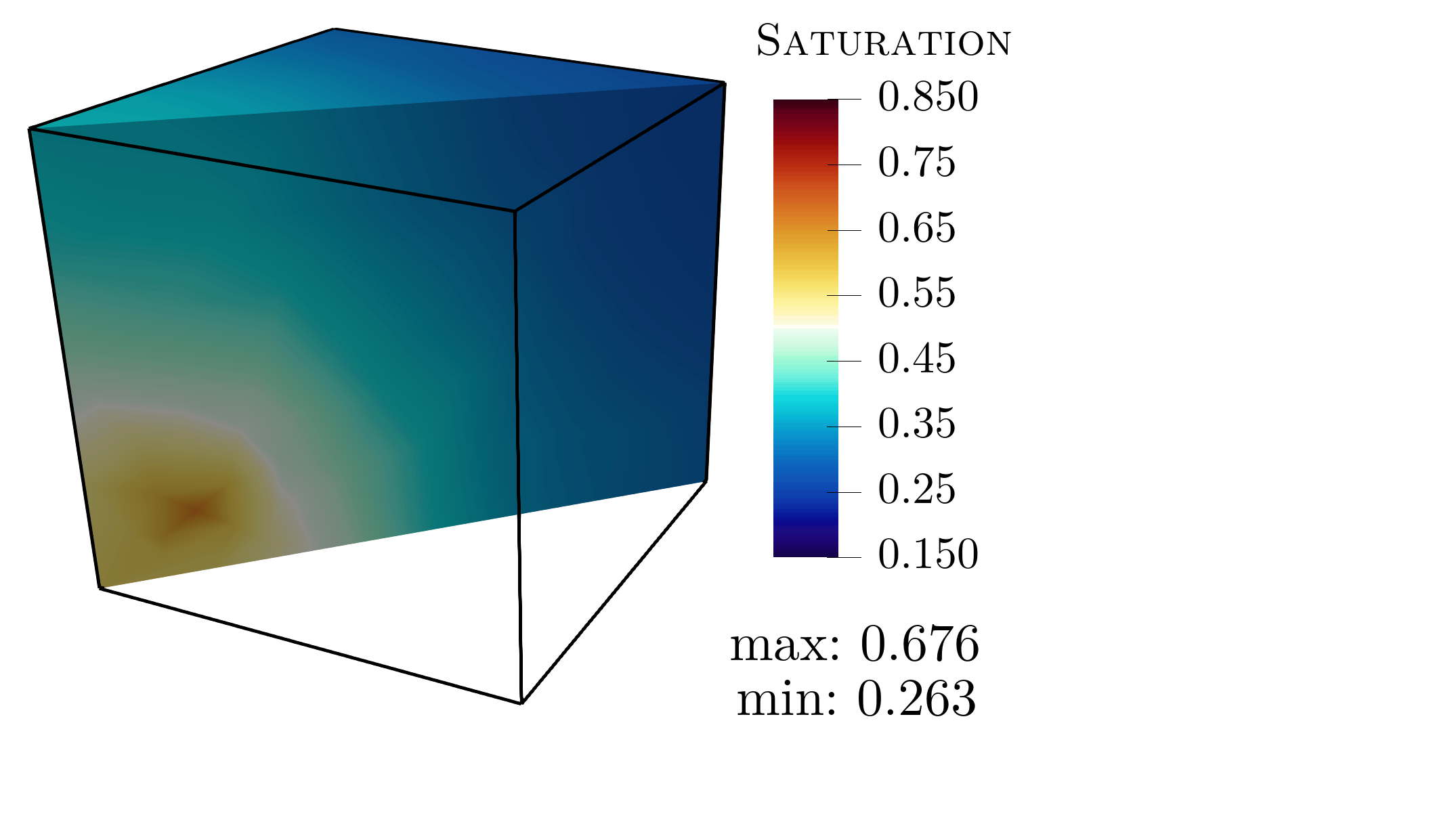}}
        \vspace{0.5cm}
	\subfigure[$t=1.0$ days \label{Fig:BoM_vertex_1800}]{
		\includegraphics[clip,scale=0.135,trim=0 0cm 18cm 0]{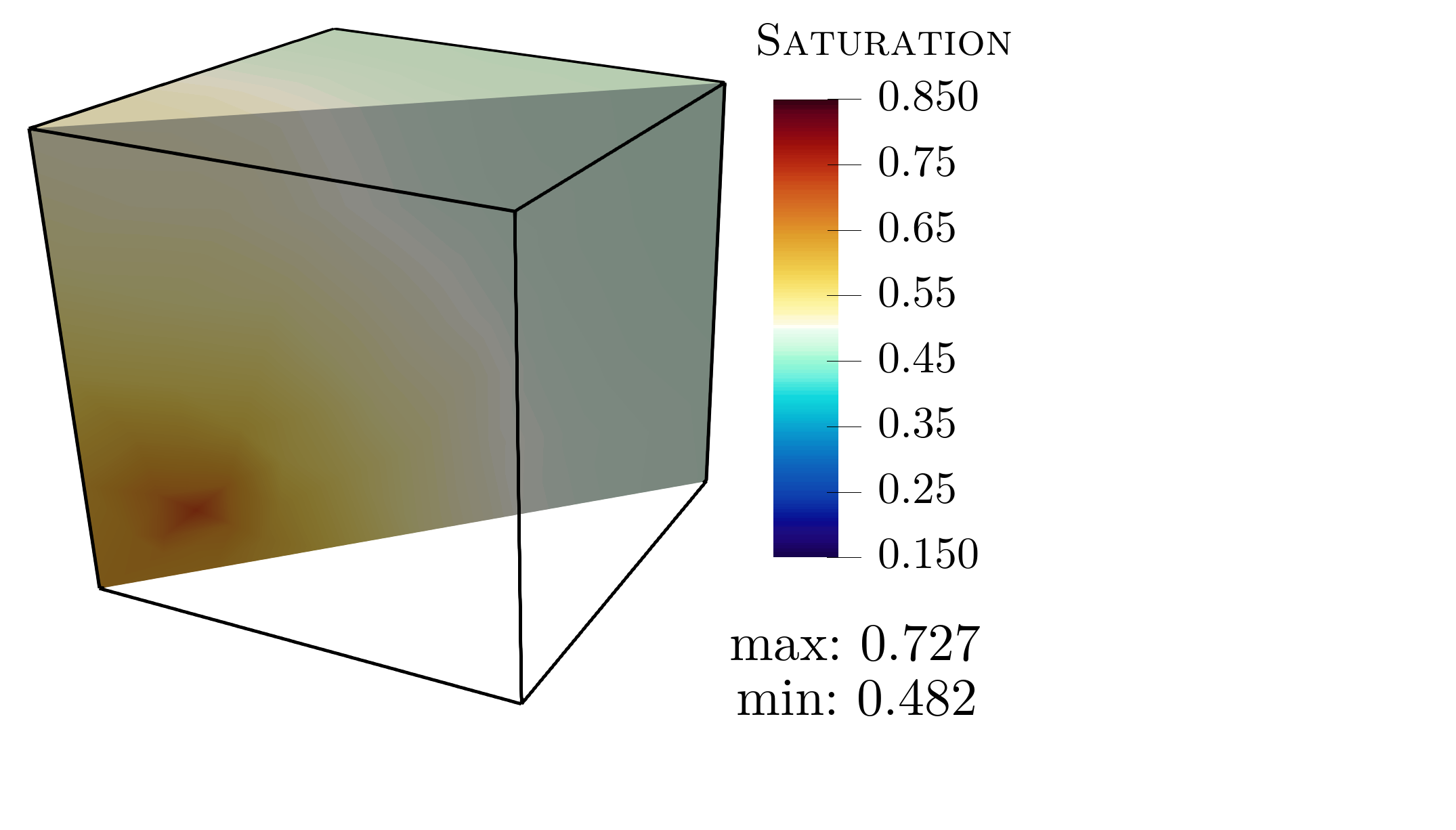}}
	\subfigure[$t=2.0$ days \label{Fig:BoM_vertex_1800}]{
		\includegraphics[clip,scale=0.135,trim=0 0cm 18cm 0]{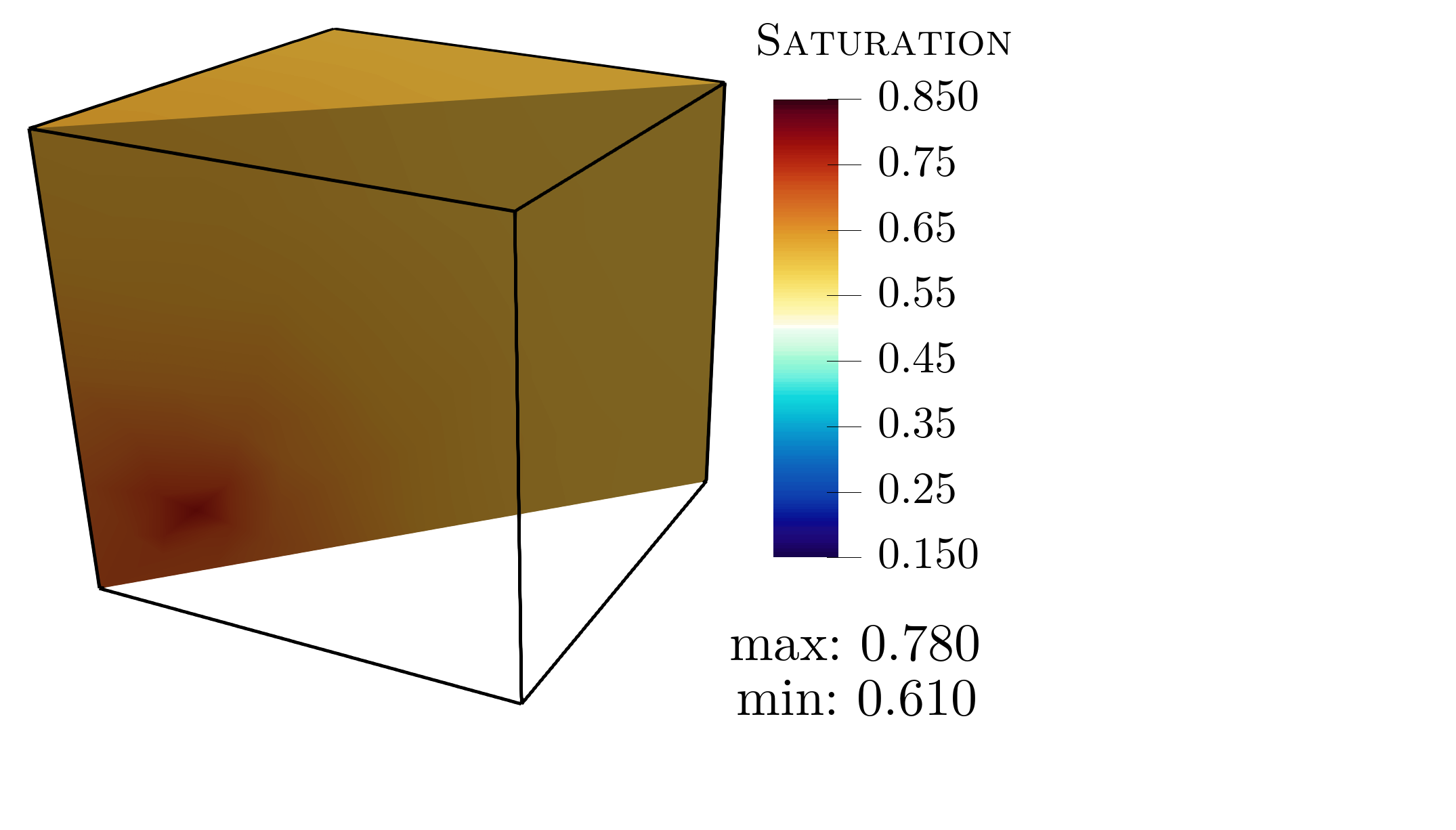}}
        \caption{
            \textsf{Three-dimensional porous medium:}
            This figure shows saturation contours in a homogeneous domain with $K=5 \times 10^{-8} \; \mathrm{m^2}$.
        Unstructured tetrahedron mesh (Figure~\ref{Fig:3D_mesh}) is used for this experiment.
        Similar to two-dimensional problems, the scheme exhibits satisfactory results with respect to maximum principle. 
        This means that solutions always remain between $0.15$ and $0.85$. 
        In this figure a cutaway view of solutions is provided for better visualization.
        \label{Fig:3D_sat}}
\end{figure}
\begin{figure}
  \includegraphics[clip,width=0.8\linewidth,trim=0 3.5cm 0 3cm]{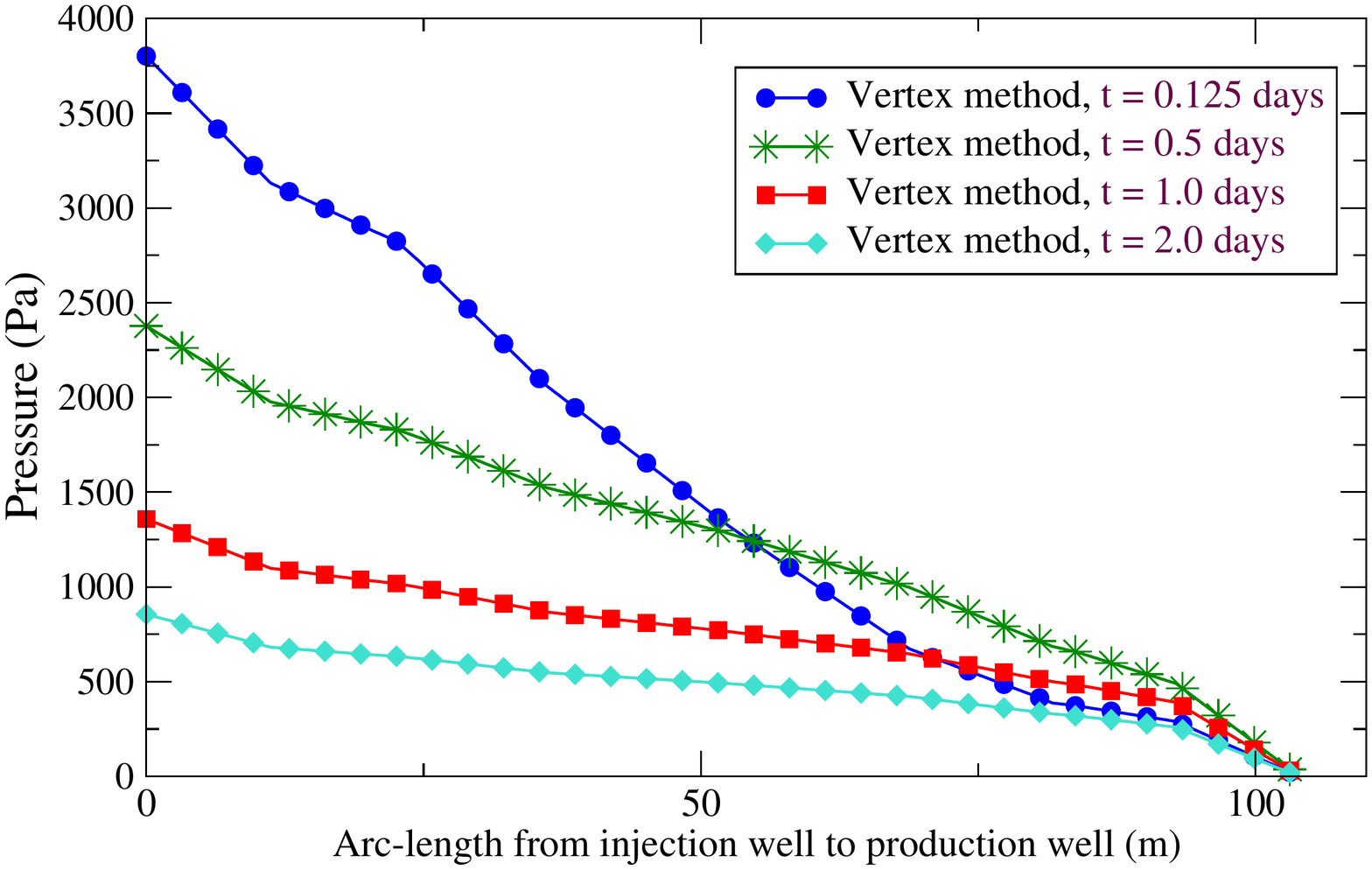}
  \caption{
        \textsf{Three-dimensional porous medium:}
        This figure shows pressure solutions obtained from the vertex scheme on a homogeneous domain with $K=5 \times 10^{-8}
        \;\mathrm{m^2}$.
        Profiles are plotted on the diagonal line from point $(20, 20, 20)\;\mbox{m}$ to point $(80,80,80)\;\mbox{m}$.
        Vertex scheme correctly predicts the response. Maximum and minimum pressures are detected on the injection well and 
        production well, respectively. As time advances, more fluid reaches the production well and the pressure 
        difference decreases.
  \label{Fig:3D_pressure}}
\end{figure}

\subsubsection{Three-dimensional porous medium with highly heterogeneous permeability}
In realistic problems, heterogeneities in three-dimensional media have a large impact on the propagation of the fluid phases.
We now examine a three-dimensional problem analogous to the 2D numerical experiment carried out in Section~\ref{Sub:SPE10-2D}.
The aim of this boundary value problem is to show that the proposed finite element method can perform satisfactorily in highly 
heterogeneous three-dimensional domains. 
The domain is $\Omega=[0,50]\times[0,100]\times[0,24]~\mathrm{m}^3$. 
As shown in Figure~\ref{Fig:SPE3D_K}, we adopt a sample permeability field of size $32\;\mbox{m}\times64\;\mbox{m}\times 12\;\mbox{m}$ 
from the SPE10 benchmark problem~\citep{SPE10}. The coordinates of injection and production wells are  
$(7.5,15,4)$~m and $(42.5,85,4)$~m, respectively (see Figure~\ref{Fig:3D_SPE10_schemetic}). 
The size of both wells are $5\;\mbox{m}\times10\;\mbox{m}\times2\;\mbox{m}$ with $\int_{\Omega}\bar{q}=\int_{\Omega}\underline{q}=0.1$.
The mesh is made of $19200$ structured tetrahedral elements. Total time is set to $T=3$~days and 
the time step is~$\tau = 259.2$~s. 
\begin{figure}
  \includegraphics[clip,width=0.6\linewidth]{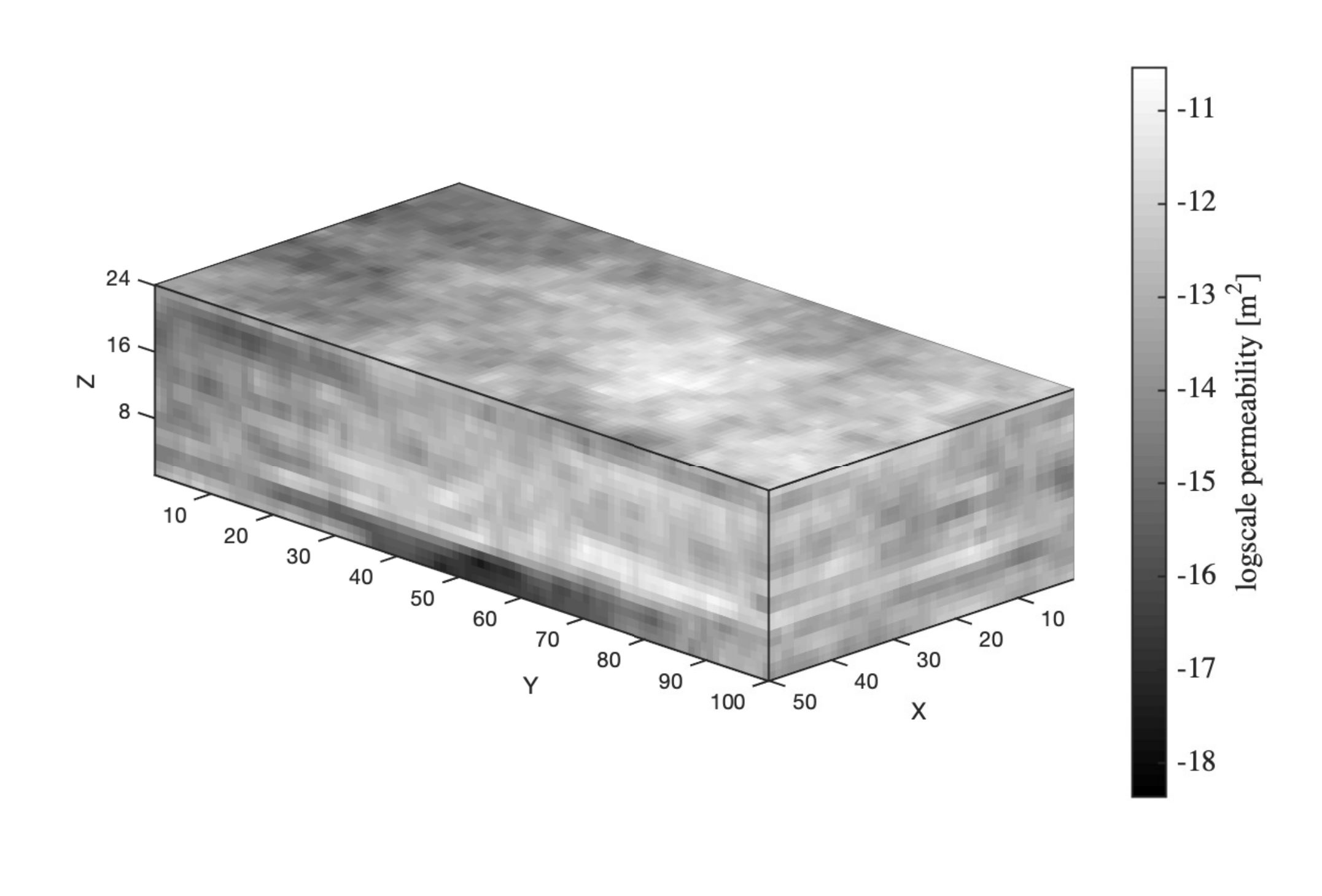}
  \caption{
  \textsf{Three-dimensional porous medium with highly heterogeneous permeability:}
  This figure shows the permeability field extracted from SPE10 benchmark model.
  The resolution of this field is $32\times64\times 12$ grids. 
  Values are displayed in logarithmic scale, since they vary across a wide range.
  \label{Fig:SPE3D_K}}
\end{figure}
We apply the proposed finite element method and plot the saturation contours at different time steps in Figure~\ref{Fig:SPE10_sat}.
The wetting phase reaches the production well by sweeping the regions with highest permeability value.
Evidently, the numerical saturation remains within physical bounds and no undershoots and overshoots are observed.
This experiments reinforces that the proposed  scheme satisfies the  maximum-principle for two-phase incompressible flow 
and remains robust for highly heterogeneous three-dimensional media.

\begin{figure}
	\subfigure[$t=0.15$ days  \label{Fig:BoM_vertex_1800}]{
		\includegraphics[clip,scale=0.16,trim=0 7.5cm 0 7.5cm]{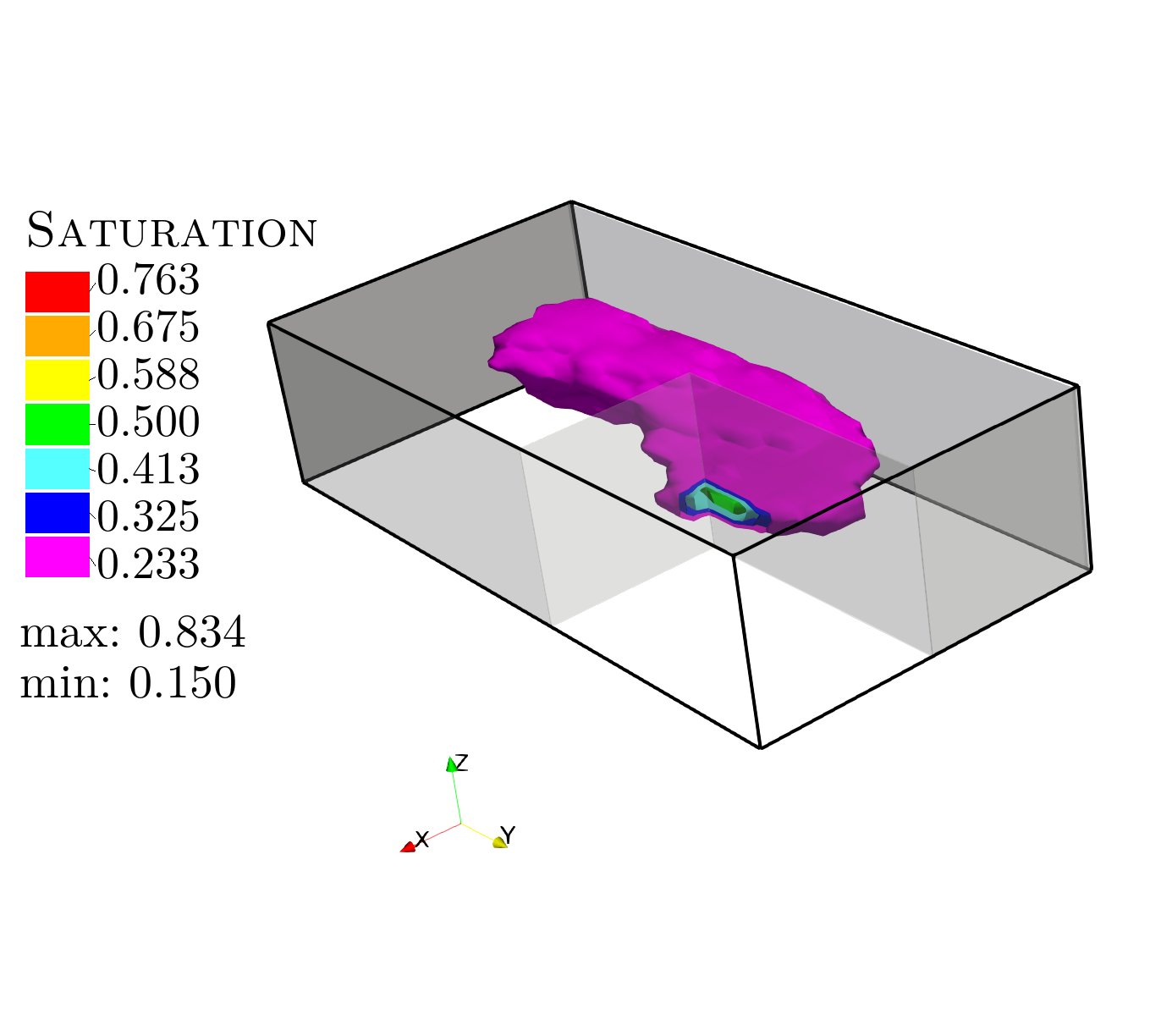}}
        \hfill 
	\subfigure[$t=0.39$ days  \label{Fig:BoM_vertex_1800}]{
		\includegraphics[clip,scale=0.16,trim=0 7.5cm 0 7.5cm]{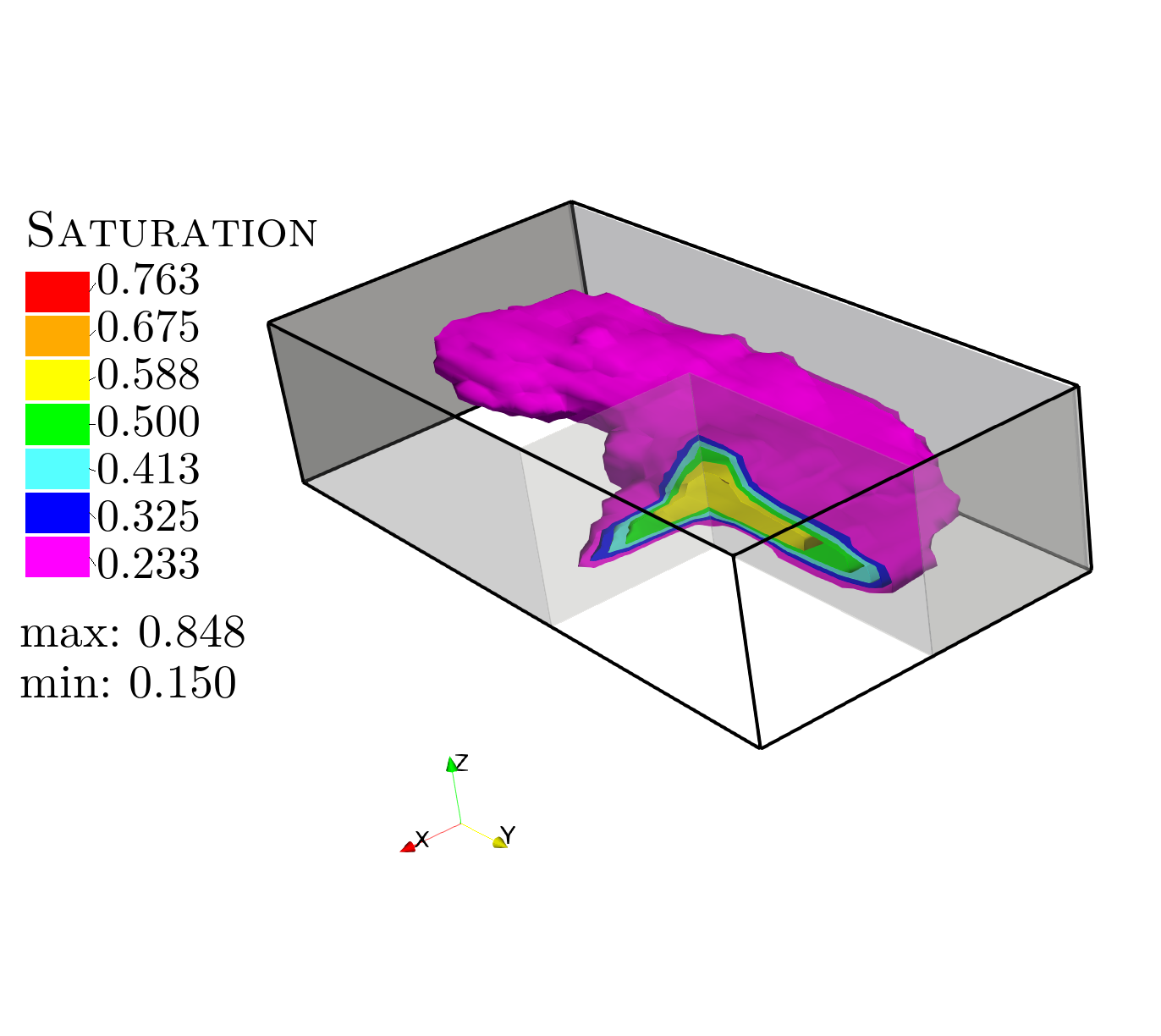}}
        \vspace{0.5cm}
	\subfigure[$t=0.75$ days  \label{Fig:BoM_vertex_1800}]{
		\includegraphics[clip,scale=0.16,trim=0 7.5cm 0 7.5cm]{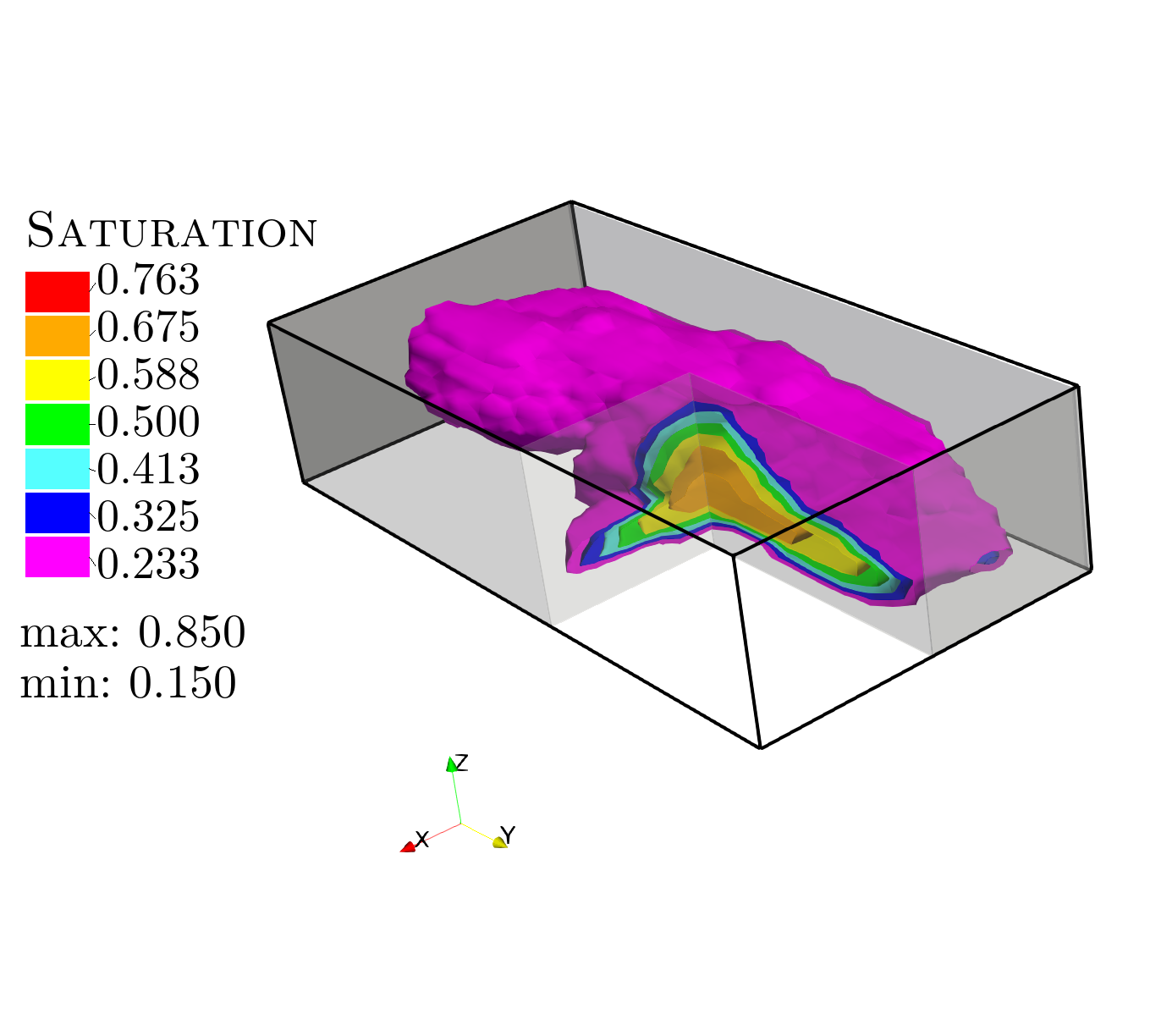}}
        \hfill
	\subfigure[$t=3$ days  \label{Fig:BoM_vertex_1800}]{
		\includegraphics[clip,scale=0.16,trim=0 7.5cm 0 7.5cm]{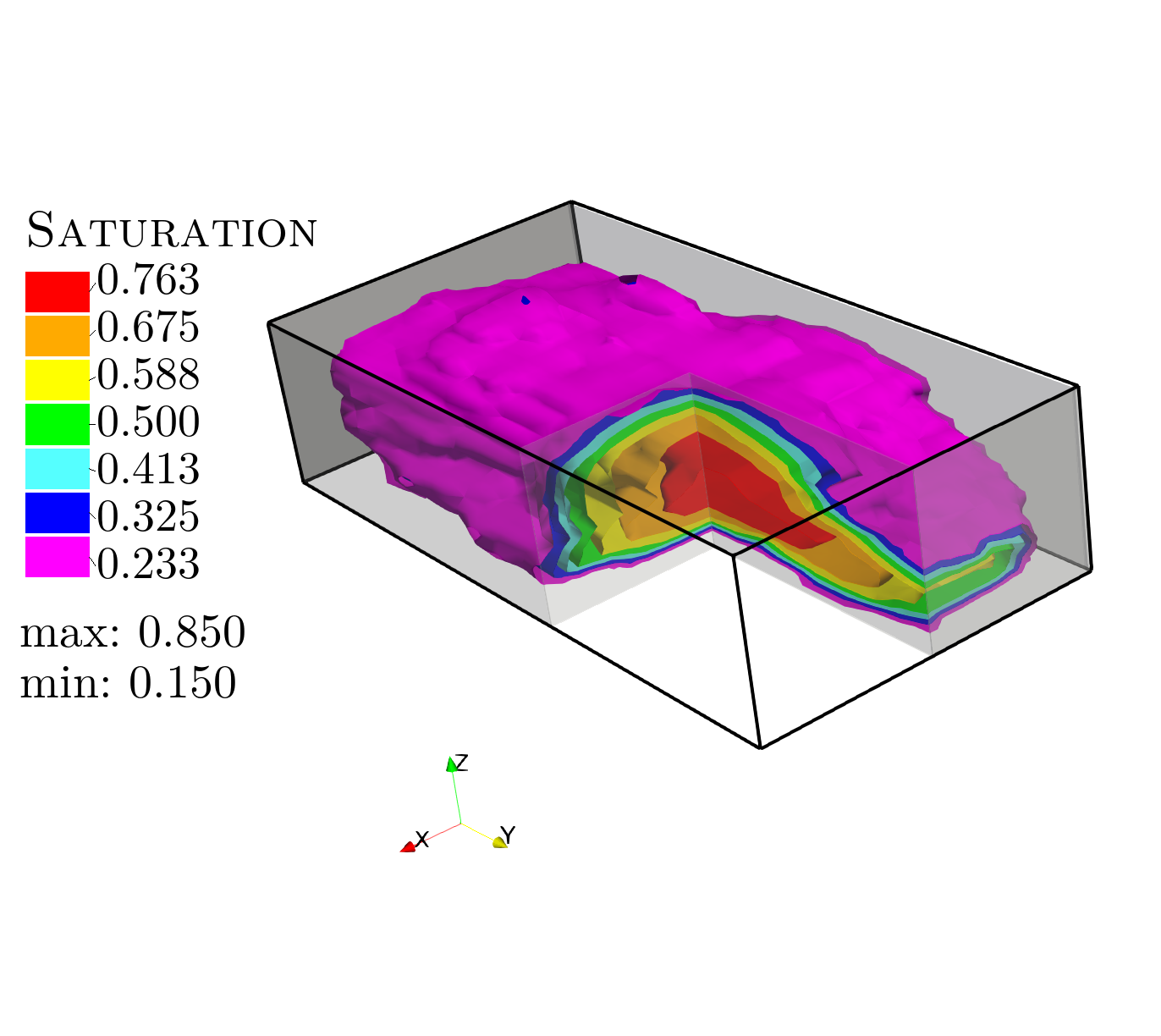}}
        \caption{
            \textsf{Three-dimensional SPE10 problem:}
            This figure shows contours of the saturation solutions obtained under the proposed vertex scheme on a 
            highly heterogeneous domain (i.e., SPE10 permeability field).
            The numerical experiment performed on a structured tetrahedron mesh with $19200$ elements.
            The main inference from this figure is that 
            (i) the proposed scheme generates robust and accurate results. It can be seen that the wetting phase fluid flows 
            through the most permeable pore-networks from injection well towards production well.
            (ii) The scheme always respects maximum principle, since no undershoots and overshoots has been detected 
            throughout the simulation.
             (For interpretation of the references to color in this figure legend, the reader is referred to the web 
             version of this article.)
        \label{Fig:SPE10_sat}}
\end{figure}
\section{CONCLUSION}
\label{Sec:Conclusion}
We have developed a new first-order finite element method with mass-lumping and flux upwinding, which we refer to as vertex scheme, 
to solve the immiscible two-phase flow problem in porous media. We show optimal convergence rates for manufactured solutions.
Numerical examples in two and three dimensions pinpoint that the method is accurate, and robust, even in the case of realistic 
discontinuous highly varying permeability.
Furthermore, we show that the proposed method is locally mass-conservative and the resulting solutions satisfy the maximum principle.
The method is mesh-independent and does not require penalization or any external bound-preserving mechanism.

\clearpage
\appendix
\section*{APPENDIX}
\label{Sec:Appendix}
Below, we have provided the PETSc command-line options for the Schur complement approach discussed in Section~\ref{Sec:Solver}. 
(see online version for color-coded terms).
\vspace{-\baselineskip}
\begin{lstlisting}[language=Python,
caption=PETSc command-line options for the proposed Schur complement approach, abovecaptionskip=5pt, belowcaptionskip=5pt,
label=Code:ex_fields,frame=single]
# Outer solver
PETScOptions.set("ksp_type", "fgmres")
PETScOptions.set("ksp_rtol", 1e-8)

# Schur complement with full factorization
PETScOptions.set('pc_type', 'fieldsplit')
PETScOptions.set('pc_fieldsplit_type', 'schur')
PETScOptions.set('pc_fieldsplit_schur_fact_type', 'full')

# Diagonal mass lumping
PETScOptions.set('pc_fieldsplit_schur_precondition', 'selfp')

# Automatically determine fields based on zero diagonal entries
PETScOptions.set('pc_fieldsplit_detect_saddle_point')

# Single sweep of ILU(0) for the mass matrix
PETScOptions.set('fieldsplit_0_ksp_type', 'preonly')
PETScOptions.set('fieldsplit_0_pc_type', 'ilu')

# Single sweep of multi-grid for the Schur complement
PETScOptions.set('fieldsplit_1_ksp_type', 'preonly')
PETScOptions.set('fieldsplit_1_pc_type', 'hypre')
\end{lstlisting}

\bibliographystyle{plainnat}
\bibliography{./Master_References/Master_References,./Master_References/Books}
\end{document}